\documentclass[11pt]{article}
\usepackage{amssymb}
\usepackage{latexsym}
\usepackage{amsfonts}
\usepackage{amsmath}

\newtheorem{thrm}{Theorem}[section]
\newtheorem{lemma}[thrm]{Lemma}

\newtheorem{cor}[thrm]{Corollary}
\newtheorem{remark}[thrm]{Remark}
\newtheorem{example}[thrm]{Example}
\newtheorem{assumption}{Assumption}
\numberwithin{equation}{section}

\usepackage[dvips]{color}
\usepackage{showkeys}

\def\cL{{\cal L} }
\def\P{\mathbb{P} }
\def\V{\mathbb{V} }
\def\E{\mathbb{E} }
\def\N{\mathcal{N}}
\def\R{\mathcal{R}}

\def\C{\mathcal{C}}
\def\A{\mathcal{A}}
\def\I{\mathbb{I}}

\begin{document}
\allowdisplaybreaks
\begin{doublespace}
\title{\Large\bf Central Limit Theorems for Supercritical Branching Nonsymmetric Markov Processes}
\author{ \bf  Yan-Xia Ren\footnote{The research of this author is supported by NSFC (Grant No.  11271030 and 11128101) and Specialized Research Fund for the
Doctoral Program of Higher Education.\hspace{1mm} } \hspace{1mm}\hspace{1mm}
Renming Song\thanks{Research supported in part by a grant from the Simons
Foundation (208236).} \hspace{1mm}\hspace{1mm} and \hspace{1mm}\hspace{1mm}
Rui Zhang\footnote{Supported by the China Scholarship Council.}
\hspace{1mm} }
\date{}
\maketitle

\begin{abstract}
In this paper, we establish a spatial central limit theorem for a large class of supercritical
branching, not necessarily symmetric, Markov processes with spatially dependent
branching mechanisms satisfying a second moment condition.
This central limit theorem generalizes and unifies all the central limit theorems obtained recently in \cite{RSZ2}
for supercritical branching symmetric Markov processes.
To prove our central limit theorem, we have to carefully develop the spectral theory
of nonsymmetric strongly continuous semigroups which should be of independent interest.
\end{abstract}

\medskip
\noindent {\bf AMS Subject Classifications (2000)}: Primary 60F05, 60J80;
Secondary 60J25, 60J35

\medskip

\noindent{\bf Keywords and Phrases}: Central limit theorem,
branching Markov process, supercritical, martingale.

\bigskip

\baselineskip=6.0mm

\section{Introduction}

Central limit theorems for supercritical branching processes were initiated by Kesten and Stigum
in \cite{KS, KS66}.
In these two papers, they established central limit theorems for
supercritical multi-type Galton-Watson processes by using the
Jordan canonical form of the expectation matrix $M$.
Then in \cite{Ath69a, Ath69, Ath71}, Athreya proved central limit
theorems for supercritical multi-type continuous time branching processes,
using the Jordan canonical form and the eigenvectors of the matrix $M_t$, the mean matrix at time $t$.
Asmussen and Keiding \cite{AK} used martingale central limit theorems to prove central limit theorems for
supercritical multi-type branching processes. In \cite{AH83}, Asmussen and Hering established
spatial central limit theorems for general supercritical branching Markov processes under a certain
condition. However, the condition in \cite{AH83} is not easy to check and essentially the only
examples given in \cite{AH83} of branching Markov processes satisfying this condition are branching
diffusions in bounded smooth domains. We note that the limit normal random variables in \cite{AH83}
may be degenerate.

 The recent study of spatial central limit theorem for branching Markov processes started with
\cite{RP}.
In this paper, Adamczak and Mi{\l}o\'{s} proved some central limit theorems for supercritical branching
Ornstein-Uhlenbeck processes with binary branching mechanism.
We note that branching Ornstein-Uhlenbeck processes do not satisfy the condition in \cite{AH83}.
In \cite{Mi}, Mi{\l}o\'{s} proved
some central limit theorems for supercritical super Ornstein-Uhlenbeck processes with branching
mechanisms satisfying a fourth moment condition.
Similar to the case of \cite{AH83}, the limit normal random variables in \cite{RP, Mi}
may be degenerate.
In \cite{RSZ}, we established central limit
theorems for supercritical super Ornstein-Uhlenbeck processes with
branching mechanisms satisfying only a second moment condition. More importantly, the central limit
theorems in \cite{RSZ} are more satisfactory since our limit normal random variables are non-degenerate.
In \cite{RSZ2}, we obtained central limit theorems for a large class
of general supercritical branching symmetric Markov processes with  spatially dependent branching mechanisms satisfying
only a second moment condition. In \cite{RSZ3}, we obtained central limit theorems for a large class of
general supercritical superprocesses with symmetric spatial motions and with spatially dependent branching mechanisms satisfying
only a second moment condition. Furthermore, we also obtained the covariance structure of the limit
Gaussian field in \cite{RSZ3}.

Compared with \cite{Ath69a, Ath69, Ath71, KS, KS66}, the spatial processes in \cite{RP, Mi, RSZ, RSZ2, RSZ3}
are assumed to be symmetric. The reason for this assumption is that
one of the main tools in \cite{RP, Mi, RSZ, RSZ2, RSZ3} is the well-developed spectral theory
of self-adjoint operators.

The main purpose of this  paper is to establish central limit theorems
for general supercritical branching, not necessarily symmetric,
Markov processes  with  spatially dependent branching mechanisms satisfying
only a second moment condition. To accomplish this, we need to carefully develop
the spectral theory of not necessarily symmetric strongly continuous semigroups.
We believe these spectral results are of independent interest and should be very
useful in studying non-symmetric Markov processes.

In this paper, $\mathbb{R}$ and $\mathbb{C}$ stand for the sets of real and complex numbers respectively,
all vectors in $\mathbb{R}^n$ or $\mathbb{C}^n$ will be understood
as column vectors. For any $z\in \mathbb{C}$, we use $\Re(z)$ and $\Im(z)$ to denote
real and imaginary parts of $z$ respectively. For a matrix $A$, we use $\overline{A}$
and $A^T$ to denote the conjugate and transpose of $A$ respectively.

\subsection{Spatial process}\label{subs:sp}

In this subsection, we spell out our assumptions on the spatial Markov process.
Throughout this paper, $E$ stands for
a locally compact separable metric
space, $m$ is a $\sigma$-finite Borel measure
on $E$ with full support and $\partial$ is a separate point not contained
in $E$. $\partial$ will be interpreted as the cemetery point.
We will use $E_{\partial}$ to denote $E\cup\{\partial\}$.
Every function $f$ on $E$ is automatically extended to $E_{\partial}$ by setting $f(\partial)=0$.
We will assume that $\xi=\{\xi_t,\Pi_x\}$ is a Hunt process on $E$ and $\zeta:=
 \inf\{t>0: \xi_t=\partial\}$ is the lifetime of $\xi$.
We will use $\{P_t:t\geq 0\}$ to denote the semigroup of $\xi$.
Our standing assumption on $\xi$ is that there exists a family of
continuous strictly positive functions $\{p(t,x,y):t>0\}$ on $E\times E$ such that,
for any $t>0$ and nonnegative function $f$ on $E$,
$$
  P_tf(x)=\int_E p(t,x,y)f(y)m(dy).
$$
For $p\ge 1$, we define
$L^p(E,m;\mathbb{C}):=\{f:E\to\mathbb{C}:\int_{E}|f(x)|^p\,m(dx)<\infty\}$
and $L^p(E,m):=\{f\in L^p(E,m;\mathbb{C}): f \mbox{ is real}\}$ .
We also define
$$a_t(x):=\int_E p(t,x,y)^2\,m(dy),\qquad \widehat{a}_t(x):=\int_E p(t,y,x)^2\,m(dy).$$
In this paper, we assume that
\begin{assumption}
\begin{description}
  \item[(a)]
  For all $t>0$ and  $x\in E$, $\int_E p(t,y,x)\,m(dy)\le 1$.
  \item[(b)] For any $t>0$,
  $a_t$ and $\widehat{a}_t$ are continuous functions in $E$ and they belong to $L^1(E,m)$.
  \item[(c)] There exists $t_0>0$ such that $a_{t_0},\widehat{a}_{t_0}\in L^2(E,m)$.
\end{description}
\end{assumption}

It is easy to see that
\begin{equation}\label{1.1}
  p(t+s,x,y)=\int_E p(t,x,z)p(s,z,y)\,m(dz)\le (a_t(x))^{1/2}(\widehat{a}_s(y))^{1/2},
\end{equation}
which implies
\begin{eqnarray}\label{1.2}
  a_{t+s}(x)\le \int_E \widehat{a}_s(y)\,m(dy)a_t(x)\quad \mbox{and}\quad \widehat{a}_{t+s}(x)\le \int_E a_s(y)\,m(dy)\widehat{a}_t(x).
\end{eqnarray}
So condition $(c)$ above is equivalent to
 \begin{description}
   \item[(c$'$)] There exists $t_0>0$ such that for all $t\ge t_0$, $a_{t},\widehat{a}_t\in L^2(E,m)$.
 \end{description}
 It is well known and easy to check that, for $p\in [1, \infty)$, $\{P_t:t\ge 0\}$ is a
strongly continuous contraction semigroup on $L^p(E, m;\mathbb{C})$.
We claim that the function $t\rightarrow \int_E a_{t}(x)\,m(dx)$ is decreasing.
In fact, by Fubini's theorem and  H\"{o}lder's inequality, we get
\begin{eqnarray*}
  a_{t+s}(x) &=& \int_E p(t+s,x,y)\int_{E}p(t,x,z)p(s,z,y)\,m(dz)\,m(dy)\\
&=& \int_E p(t,x,z)\int_{E}p(t+s,x,y)p(s,z,y)\,m(dy)\,m(dz)\\
&\le&  a_{t+s}(x)^{1/2}\int_E p(t,x,z)a_s(z)^{1/2}\,m(dz)
\end{eqnarray*}
which implies
\begin{equation}\label{8.9}
  a_{t+s}(x)\le \left(\int_E p(t,x,z)a_s(z)^{1/2}\,m(dz)\right)^2\le \int_E p(t,x,z)a_s(z)\,m(dz).
\end{equation}
Thus, by Fubini's theorem and condition $(a)$, we get
\begin{equation}\label{8.10}
  \int_E a_{t+s}(x)\,m(dx)\le \int_E a_{s}(z)\int_E p(t,x,z)\,m(dx)\,m(dz)\le \int_E a_{s}(z)\,m(dz).
\end{equation}
Therefore, the function $t\rightarrow \int_E a_{t}(x)\,m(dx)$ is decreasing.

Now we give some examples of non-symmetric Markov processes satisfying the above assumptions.
The purpose of these examples is to show that the above assumptions are satisfied by many
Markov processes. We will not try to give the most general examples possible. For examples
of symmetric Markov processes satisfying the above assumptions, see \cite{RSZ2}.

\begin{example}\label{examp0}
{\rm
Suppose that $E$ consists of finitely many points.
If $X=\{X_t: t\ge 0\}$ is an irreducible conservative Markov process in $E$,
then $X$ satisfies Assumption 1 for some finite measure $m$ on $E$ with  full support.
}
\end{example}

\begin{example}\label{examp1}
{\rm Suppose that $\alpha\in (0, 2)$ and that $Y=\{Y_t: t\ge0\}$ is a strictly
$\alpha$-stable process in $\mathbb{R}^d$. Suppose that, in the case $d\ge 2$, the spherical
part $\eta$ of the L\'evy measure $\mu$ of $Y$ satisfies the following assumption: there
exist a positive function $\Phi$ on the unit sphere $S$ in $\mathbb{R}^d$ and $\kappa>1$
such that
$$
\Phi=\frac{d\eta}{d\sigma} \quad \mbox{and} \quad
\kappa^{-1}\le \Phi(z)\le \kappa \quad \mbox{on } S
$$
where $\sigma$ is the surface measure on $S$. In the case $d=1$, we assume that the L\'evy
measure of $Y$ is given by
$$
\mu(dx)=c_1x^{-1-\alpha}1_{\{x>0\}}+ c_2|x|^{-1-\alpha}1_{\{x<0\}}
$$
with $c_1, c_2>0$. Suppose that $D$ is an open set in $\mathbb{R}^d$ of finite Lebesgue measure.
Let $X$ be the process in $D$ obtained by killing $Y$ upon exiting $D$. Then $X$ satisfies Assumption 1
with $E=D$ and $m$ being the Lebesgue measure. For details, see \cite[Example 4.1]{KiSo09}.
}
\end{example}

\begin{example}\label{examp2}
{\rm Suppose that $\alpha\in (0, 2)$ and that $Z=\{Z_t: t\ge 0\}$ is a truncated
strictly $\alpha$-stable process in $\mathbb{R}^d$, that is, $Z$ is a L\'evy
process with L\'evy measure given by
$$
\widetilde{\mu}(dx)=\mu(dx)1_{\{|x|<1\}},
$$
where $\mu$ is the L\'evy measure of the process $Y$ in the previous example.
Suppose that $D$ is a connected open set in $\mathbb{R}^d$ of finite Lebesgue measure.
Let $X$ be the process in $D$ obtained by killing $Z$ upon exiting $D$. Then $X$ satisfies Assumption 1
with $E=D$ and $m$ being the Lebesgue measure. For details,
see \cite[Example 4.2 and Proposition 4.4]{KiSo09}.
}
\end{example}

\begin{example}\label{examp3}
{\rm
Suppose $\alpha\in (0, 2)$, $Y=\{Y_t: t\ge0\}$ is a strictly
$\alpha$-stable process in $\mathbb{R}^d$ satisfying the assumptions in Example
\ref{examp1} and that $B$ is an independent Brownian motion in $\mathbb{R}^d$. Let $W$ be the
process defined by $W_t=Y_t+B_t$. Suppose that $D$ is an open set in $\mathbb{R}^d$ of finite Lebesgue measure.
Let $X$ be the process in $D$ obtained by killing $W$ upon exiting $D$. Then $X$ satisfies Assumption 1
with $E=D$ and $m$ being the Lebesgue measure. For details, see \cite[Example 4.5 and
Lemma 4.6]{KiSo09}.
}
\end{example}

\begin{example}\label{examp4}
{\rm
Suppose $\alpha\in (0, 2)$, $Z=\{Z_t: t\ge0\}$ is a truncated strictly
$\alpha$-stable process in $\mathbb{R}^d$ satisfying the assumptions in Example
\ref{examp2} and that $B$ is an independent Brownian motion in $\mathbb{R}^d$. Let $V$ be the
process defined by
$V_t=Z_t+B_t$.
Suppose that $D$ is
a connected open set in $\mathbb{R}^d$ of finite Lebesgue measure.
Let $X$ be the process in $D$ obtained by killing $V$ upon exiting $D$. Then $X$ satisfies Assumption 1
 with $E=D$ and $m$ being the Lebesgue measure. For details, see \cite[Example 4.7 and
Lemma 4.8]{KiSo09}.
}
\end{example}

\begin{example}\label{examp5}
{\rm
Suppose $d\ge 3$ and that $\mu=(\mu^1, \cdots, \mu^d)$, where each $\mu^j$ is a signed measure
on
$\mathbb{R}^d$ such that
$$
 \lim_{r\to 0}\sup_{x\in\mathbb{R}^d}\int_{B(x, r)}\frac{|\mu^j|(dy)}{|x-y|^{d-1}}=0.
$$
Let $Y=\{Y_t: t\ge 0\}$ be a Brownian motion with drift $\mu$ in $\mathbb{R}^d$, see \cite{KiSo06}.
Suppose that $D$ is a bounded connected open set in $\mathbb{R}^d$
and suppose $K>0$ is a constant such that
$D\subset B(0, K/2)$. Put $B=B(0, K)$. Let $G_B$ be the Green function of $Y$ in $B$ and define
$H(x):=\int_BG_B(x, y)dy$. Then $H$ is a strictly positive continuous function on $B$.
Let $X$ be the
process obtained by killing $Y$ upon exiting $D$.
Then $X$ satisfies Assumption 1
with $E=D$ and $m$ being the measure defined by $m(dx)=H(x)dx$.
For details, see \cite[Example 4.6]{ZLS} or \cite{KiSo08, KiSo08c}.
}
\end{example}

\begin{example}\label{examp6}
{\rm
Suppose $d\ge 2$,
$\alpha\in (1, 2)$,
and that $\mu=(\mu^1, \cdots, \mu^d)$, where each $\mu^j$ is a signed measure
on $\mathbb{R}^d$ such that
$$
 \lim_{r\to 0}\sup_{x\in\mathbb{R}^d}\int_{B(x, r)}\frac{|\mu^j|(dy)}{|x-y|^{d-\alpha+1}}=0.
$$
Let $Y=\{Y_t: t\ge 0\}$ be an $\alpha$-stable process with drift $\mu$ in $\mathbb{R}^d$, see \cite{KiSo13}.
Suppose that $D$ is a bounded open set in $\mathbb{R}^d$ and suppose $K>0$ is such that
$D\subset B(0, K/2)$. Put $B=B(0, K)$. Let $G_B$ be the Green function of $Y$ in $B$ and define
$H(x):=\int_BG_B(x, y)dy$. Then $H$ is a strictly positive continuous function on $B$.
Let $X$ be the
process obtained by killing $Y$ upon exiting $D$.
Then $X$ satisfies Assumption 1
with $E=D$ and $m$ being the measure defined by $m(dx)=H(x)dx$.
For details, see \cite[Example 4.7]{ZLS} or \cite{CKS}.
}
\end{example}

 \subsection{Branching Markov Processes}
 The branching Markov process $\{X_t: t\ge 0\}$ on $E$ we are going to work with is
determined by three parameters: a spatial motion $\xi=\{\xi_t, \Pi_x\}$ on $E$ satisfying the
assumptions at the beginning of the previous subsection,
 a branching rate function $\beta(x)$ on $E$ which is a non-negative bounded measurable function
and an offspring distribution $\{p_n(x): n=0, 1,, 2, \dots\}$ satisfying the assumption
\begin{equation}\label{1.16}
  \sup_{x\in E}\sum_{n=0}^\infty n^2p_n(x)<\infty.
\end{equation}
We denote the generating function of the offspring distribution by
$$
  \varphi(x,z)=\sum_{n=0}^\infty p_n(x)z^n,\quad x\in E,\quad |z|\leq 1.
$$

Consider a branching system on $E$ characterized by the following properties:
(i) each individual has a random birth and death time;
(ii) given that an individual is born at $x\in E$, the conditional distribution of its
path is determined by $\Pi_x$;
(iii) given the path  $\xi$ of an individual up to time $t$
and given that the particle is alive at time $t$ , its
probability of dying in the interval $[t, t +dt)$ is $\beta(\xi_t)dt + o(dt)$;
(iv) when an individual dies at $x\in E$,
it splits into $n$ individuals all positioned at $x$, with probability $p_n(x)$;
(v) when an individual reaches $\partial$, it disappears from the system;
(vi) all the individuals, once born, evolve independently.

Let $\mathcal{M}_a(E)$ be the space of finite integer-valued atomic measures on $E$,
and let $\mathcal{B}_b(E)$  be the set of bounded real-valued Borel measurable functions on $E$.
Let $X_t(B)$ be the number of particles alive at time $t$ located in $B\in \mathcal{B}(E)$.
Then $X=\{X_t,t\geq 0\}$ is an $\mathcal{M}_a(E)$-valued Markov process.
For any $\nu\in\mathcal{M}_a(E)$, we denote the law of $X$ with initial configuration $\nu$ by $\P_\nu$.
As usual, $\langle f,\nu\rangle:=\int_E f(x)\,\nu(dx)$.
For $0\le f\in \mathcal{B}_b(E)$, let
$$
  \omega(t,x):=\P_{\delta_x}e^{-\langle f,X_t\rangle},
$$
then $\omega(t,x)$ is the unique positive solution to the equation
\begin{equation}\label{1.3}
  \omega(t,x)=\Pi_x\int_0^t \psi(\xi_s,\omega(t-s,\xi_s))\,ds+\Pi_x(e^{-f(\xi_t)}),
\end{equation}
where $\psi(x,z)=\beta(x)(\varphi(x,z)-z),x\in E, z\in [0,1],$ while $\psi(\partial,z)=0, z\in [0,1]$.
By the branching property, we have
$$
\P_{\nu}e^{-\langle f,X_t\rangle}=e^{\langle\log\omega(t,\cdot),\nu\rangle}.
$$
Define
\begin{equation}\label{e:alpha}
\alpha(x):=\frac{\partial\psi}{\partial z}(x,1)=\beta(x)\left(\sum_{n=1}^\infty np_n(x)-1\right)
\end{equation}
and
\begin{equation}\label{e:A}
A(x):=\frac{\partial^2\psi}{\partial z^2}(x,1)=\beta(x)\sum_{n=2}^\infty (n-1)n p_n(x).
\end{equation}
By \eqref{1.16}, there exists $K>0$, such that
\begin{equation}\label{1.5}
  \sup_{x\in E}\left(|\alpha(x)|+A(x)\right)\le K.
\end{equation}
 For any $f\in\mathcal{B}_b(E)$ and $(t, x)\in (0, \infty)\times E$, define
\begin{equation}\label{1.26}
   T_tf(x):=\Pi_x \left[e^{\int_0^t\alpha(\xi_s)\,ds}f(\xi_t)\right].
\end{equation}
It is well known that $T_tf(x)=\P_{\delta_x}\langle f,X_t\rangle$ for every $x\in E$.

It is elementary to show that, see \cite[Lemma 2.1]{RSZ4}, that there exists a function
$q(t, x, y)$ on $(0, \infty)\times E\times E$ which is continuous in $(x, y)$ for each $t>0$
such that
\begin{equation}\label{comp}
 e^{-Kt}p(t,x,y) \le q(t,x,y)\le e^{Kt}p(t,x,y), \quad (t, x, y)\in
 (0, \infty)\times E\times E
\end{equation}
and that for any bounded Borel function $f$ on $E$ and $(t, x)\in (0, \infty)\times E$,
$$
T_tf(x)=\int_Eq(t, x, y)f(y)m(dy).
$$

Define
$$b_t(x):=\int_E q(t,x,y)^2\,m(dy),\qquad \widehat{b}_t(x):=\int_E q(t,y,x)^2\,m(dy).$$
The functions $x\to b_t(x)$ and $x\to\widehat{b}_t(x)$ are continuous. In fact, by \eqref{1.1},
\begin{equation}\label{1.4}
 q(t,x,y)\le e^{Kt}p(t,x,y)\le e^{Kt}a_{t/2}(x)^{1/2}\widehat{a}_{t/2}(y)^{1/2}.
\end{equation}
Since $q(t,\cdot,y)$ and $a_{t/2}$ are continuous, by the dominated convergence theorem, we get $b_t$ is continuous.
Similarly, $\widehat{b}_t$ is also continuous.
Thus, it follows from \eqref{1.4} and the assumptions (b) and (c$'$) in the previous subsection that $b_t$ and
$\widehat{b}_t$ enjoy the following properties.
\begin{description}
  \item[(i)] For any $t>0$, we have $b_t\in L^1(E,m)$. Moreover, $b_t(x)$ and $\widehat{b}_t(x)$ are continuous
  in $x\in E$;

   \item[(ii)] There exists $t_0>0$ such that for all $t\ge t_0$, $b_{t},\widehat{b}_{t}\in L^2(E,m)$.
\end{description}

\subsection{Preliminaries}

For $p\ge 1$, $\{T_t: t\ge 0\}$ is a
strongly continuous semigroup on $L^p(E, m;\mathbb{C})$. In fact,
by \eqref{comp}, we get $|T_tf(x)|\le e^{Kt}P_t|f|(x).$
Thus,
\begin{equation}\label{Lp}
  \|T_tf\|_p\le e^{Kt}\|P_t|f|\|_p\le e^{Kt}\|f\|_p.
\end{equation}
For $f,g \in L^2(E,m;\mathbb{C})$, define
$$\langle f,g\rangle_m:=\int_E f(x)\overline{g(x)}\,m(dx).$$
Let $\{\widehat{T}_t, t>0\}$ be the adjoint semigroup of $\{T_t: t\ge 0\}$
on $L^2(E, m;\mathbb{C})$,
that is, for $f,g\in L^2(E, m;\mathbb{C})$,
\begin{equation}\label{adjiont}
  \langle T_tf, g\rangle_m=\langle f, \widehat{T}_tg\rangle_m.
\end{equation}
Thus,
$$
\widehat{T}_tg(x)=\int_E q(t,y,x)g(y)\,m(dy).
$$
It is well known, see for instance \cite[Corollary 1.10.6, Lemma 1.10.1]{Pa}, that
$\{\widehat{T}_t: t\ge 0\}$ is a strongly continuous semigroup on $L^2(E, m;\mathbb{C})$
and that
\begin{equation}\label{1.66}
 \|\widehat{T}_t\|_2=\|T_t\|_2\le e^{Kt}.
\end{equation}

For all $t>0$ and $f\in L^2(E,m;\mathbb{C})$, $T_tf$  and $\widehat{T}_tf$ are continuous.
In fact,
since $q(t,x, y)$ is continuous, by \eqref{1.4} and Assumption 1(b),
using the dominated convergence theorem, we get $T_tf$  and $\widehat{T}_tf$ are continuous.

 It follows from (i) above that, for any $t>0$, $T_t$ and $\widehat{T}_t$ are compact operators on $L^2(E,m;\mathbb{C})$.
Let $\mathcal{A}$  and $\widehat{\mathcal{A}}$ be the infinitesimal generators of $\{T_t:t\geq 0\}$
and $\{\widehat{T}_t:t\geq 0\}$ in $L^2(E, m;\mathbb{C})$ respectively.
Let $\sigma(\mathcal{A})$ and $\sigma(\widehat{\mathcal{A}})$ be the spectra of
$\mathcal{A}$  and $\widehat{\mathcal{A}}$ respectively.
It follows from \cite[Theorem 2.2.4 and Corollary 2.3.7]{Pa} that both
$\sigma(\mathcal{A})$ and $\sigma(\widehat{\mathcal{A}})$ consist of
eigenvalues only, and that $\mathcal{A}$ and $\widehat{\mathcal{A}}$ have the same number, say $N$, of eigenvalues.
Of course $N$ might be finite or infinite.
Let $\mathbb{I}=\{1,2,\ldots,N\}$, when $N<\infty$; otherwise $\mathbb{I}=\{1,2,\ldots\}$.
Under the assumptions of Subsection \ref{subs:sp}, using \eqref{comp} and Jentzsch's theorem (\cite[Theorem V.6.6
on page 337]{Sch}, we know that the common value $\lambda_1=\sup\Re(\sigma(\mathcal{A}))=\sup\Re(\sigma(\widehat{\mathcal{A}}))$ is an eigenvalue of multiplicity one
for both $\mathcal{A}$ and $\widehat{\mathcal{A}}$, and that an eigenfunction $\phi_1$ of $\mathcal{A}$
associated with $\lambda _1$ can be chosen to be strictly positive almost everywhere with $\|\phi_1\|_2=1$
and an eigenfunction $\psi_1$ of $\widehat{\mathcal{A}}$
associated with $\lambda _1$ can be chosen to be strictly positive almost everywhere
with $\langle \phi_1,\psi_1\rangle_m=1$.
We list the eigenvalues $\{-\lambda_k,k\in\I \}$ of $\mathcal{A}$ in an order so
$\lambda_1<\Re(\lambda_2)\le\Re(\lambda_3)\le\cdots$. Then $\{-\overline{\lambda}_k,k\in\I \}$
are the eigenvalues of $\widehat{\mathcal{A}}$.
For convenience, we define, for any positive integer $k$ not belong to $\I$, $\lambda_k=\overline{\lambda}_k=\infty$.
For $k\in \I$, we write
$\Re_k:=\Re(\lambda_k)$ and $\Im_k:=\Im(\lambda_k)$.
We use the convention $\Re_\infty=\infty$.

Let $\sigma(T_t)$ be the spectrum of $T_t$ in $L^2(E,m;\mathbb{C})$.
It follows from \cite[Theorem 2.2.4]{Pa} that
$\sigma(T_t)\setminus\{0\}:=\{e^{-\lambda_k t}:k\in \I\}.$
In particular,
$\sigma(T_1)\setminus\{0\}=\{e^{-\lambda_k},k\in \I\}$ .
\begin{remark}\label{rek4}
It is easy to see that, there exists $t^*$ such that, for any $k\ne j$, $e^{-\lambda_kt^*}\ne e^{-\lambda_jt^*}.$
So without lose of generality, we assume that, for $k\ne j$, $e^{-\lambda_k}\ne e^{-\lambda_j}.$ Otherwise, we can consider $T_{t^*}$ instead of $T_1$ in the following arguments.
\end{remark}

Now we recall some basic facts about spectral theory, for more details, see \cite[Chapter 6]{BP}.
For any
$k\in\I$,
we define $\N_{k,0}:=\{0\}$ and for $n\ge 1$,
$$
\N_{k,n}:=\N((e^{-\lambda_k}I-T_1)^n)=\{f\in L^2(E,m;\mathbb{C}): (e^{-\lambda_k}I-T_1)^nf=0\}
$$
and
$$
\R_{k,n}:=\R((e^{-\lambda_k}I-T_1)^n)=(e^{-\lambda_k}I-T_1)^n(L^2(E,m;\mathbb{C})).
$$
For each
$k\in\I$,
there exists an integer $\nu_k\ge 1$ such that
$$
  \N_{k,n} \varsubsetneqq \N_{k,n+1}, \quad n=0,1,\cdots,\nu_k-1;\quad
  \N_{k,n} = \N_{k,n+1},\quad n\ge \nu_k
$$
and
$$
  \R_{k,n} \supsetneqq \R_{k,n+1}, \quad n=0,1,\cdots,\nu_k-1;\quad
  \R_{k,n} = \R_{k,n+1},\quad n\ge \nu_k.
$$

For all $k\in \I$ and $n\ge 0$, $\N_{k,n}$ is a finite dimensional linear subspace of $L^2(E,m;\mathbb{C})$.
$\N_{k,n}$ and $\R_{k,n}$ are  invariant subspaces of $T_t$. In fact,
for any $f\in \N_{k,n}$,
$$
(e^{-\lambda_k}I-T_1)^n(T_tf)=T_t(e^{-\lambda_k}I-T_1)^nf=0,
$$
which implies that $T_tf\in \N_{k,n}$.
If $f=(e^{-\lambda_k}I-T_1)^ng$, then
$T_tf=T_t(e^{-\lambda_k}I-T_1)^ng=(e^{-\lambda_k}I-T_1)^nT_tg\in\R_{k,n}$.
Thus, $\{T_t|_{\N_{k,\nu_k}},t>0\}$ is a semigroup
on $\N_{k,\nu_k}$.
We denote the corresponding infinitesimal generator as $\mathcal{A}_{k}$.
By \cite[Theorem 6.7.4]{BP}, $\sigma(T_1|_{\N_{k,\nu_k}})=\{e^{-\lambda_k}\}$.
Since $\sigma(\A_k)\subset \sigma(\A)$, we have $\sigma(\A_k)=\{-\lambda_k\}$.
Define $n_k:=\dim(\N_{k,\nu_k}) $ and $r_k:=\dim(\N_{k,1})$.
Then from linear algebra we know that
there exists a basis $\{\phi_j^{(k)},j=1,2,\cdots,n_k\}$ of $\N_{k,\nu_k}$ such that
\begin{eqnarray}\label{Jordan}
  \A_k(\phi_1^{(k)},\phi_2^{(k)},\cdots,\phi_{n_k}^{(k)})&=&(\phi_1^{(k)},\phi_2^{(k)},\cdots,\phi_{n_k}^{(k)})\left(
  \begin{array}{cccc}
    J_{k,1} &  & &0 \\
     & J_{k,2} & &\\
     &  & \ddots&\\
     0&  &     &J_{k,r_k}\\
  \end{array}
\right)\nonumber\\
&:=&(\phi_1^{(k)},\phi_2^{(k)},\cdots,\phi_{n_k}^{(k)})D_k,
\end{eqnarray}
where
\begin{equation}\label{J_kj}
  J_{k,j}=
  \left(
  \begin{array}{ccccc}
    -\lambda_k & 1          &       &          &0 \\
               &  -\lambda_k& 1     &          &\\
               &            & \ddots& \ddots   &\\
               &            &       &-\lambda_k&1\\
    0          &            &       &          &-\lambda_k\\
  \end{array}
\right), \quad \mbox{a }  d_{k,j}\times d_{k,j} \quad \mbox{matrix}
\end{equation}
with $\sum_{j=1}^{r_k}d_{k,j}=n_k.$  $D_k$ is uniquely determined by the dimensions of
$\N_{k,n},n=1,2,\cdots, \nu_k$ (see \cite[Section 7.8]{meyer} for more details).
Here and in the remainder of this paper we use the convention that when an operator, like
$\A$ or $\A_k$ or $T_t$, acts on a vector-valued function, it acts componentwise.
For convenience, we define the following $\mathbb{C}^{n_k}$-valued functions:
$$
\Phi_k(x):=\left(\phi_1^{(k)}(x),\phi_2^{(k)}(x),\cdots,\phi_{n_k}^{(k)}(x)\right)^T.
$$
Thus, we have,
for a.e. $x\in E$,
\begin{eqnarray}
  T_t(\Phi_k)^T(x)&=&e^{-\lambda_k t}(\Phi_k(x))^T\left(
  \begin{array}{cccc}
    J_{k,1}(t) &  & &0 \\
     & J_{k,2}(t) & &\\
     &  & \ddots&\\
     0&  &     &J_{k,r_k}(t)\\
  \end{array}
\right)\nonumber\\
&:=&e^{-\lambda_k t}(\Phi_k(x))^TD_k(t),\label{T-Jordan}
\end{eqnarray}
where $J_{k,j}(t)$ is a $d_{k,j}\times d_{k,j}$ matrix given by
\begin{equation}\label{J_kjt}
  J_{k,j}(t)=
  \left(
  \begin{array}{ccccc}
    1 & t &   t^2/2!   &      \cdots    &t^{d_{k,j}-1}/(d_{k,j}-1)! \\
             0  &  1& t     &       t^2/2!   &\cdots\\
               &            & \ddots& \ddots   &\\
               &            &       &1&t\\
    0          &            &       &          &1\\
  \end{array}
\right).
\end{equation}
More details can be found in \cite[p. 609]{meyer}.
Under our assumptions, $T_t(\Phi_k)^T(x)$ is continuous.
Thus, by \eqref{T-Jordan}, we can choose $\Phi_k$ to be continuous,
which implies \eqref{T-Jordan} holds for all $x\in E$.
We note that here the matrix $D_k(t)$ satisfies the semigroup property,
that is, for $t,s>0$, $D_{k}(t+s)=D_k(t)D_k(s)$ and $D_k(t)$ is invertible with $D_k(t)^{-1}=D_k(-t)$.

For any vector $a=(a_1,\cdots, a_n)^T\in \mathbb{C}^n$,
we define the $L^{p}$ norm  of $a$ by
$|a|_p:=\left(\sum_{j=1}^n|a_j|^p\right)^{1/p}$ when $1\le p<\infty$ and
$|a|_\infty:=\max_i(|a_i|)$ when $p=\infty$.

By H\"{o}lder's inequality, $|T_{t}(\phi_j^{(k)})(x)|\le b_{t}(x)^{1/2}.$
By \eqref{T-Jordan}, we get $(\Phi_k)^T=e^{\lambda_kt}T_{t}(\Phi_k)^T(D_k(t))^{-1}.$
Thus,
\begin{equation}\label{phi}
  |\Phi_k(x)|_\infty\le c(t,k)b_{t}(x)^{1/2},
\end{equation}
where $c(t, k)$ does not depend on $x$. When we choose $t=t_0$, we get that $\phi_j^{(k)}\in L^2(E,m;\mathbb{C})\cap L^4(E,m;\mathbb{C})$.

Now we consider the corresponding formula for $\widehat{T}_t$.
We know that
$\sigma(\widehat{T}_1)\setminus\{0\}=\{e^{-\overline{\lambda_k}},k\in\I\}$.
Define
$$\widehat{\N}_{k,n}:=\N((e^{-\overline{\lambda_k}}I-\widehat{T}_1)^n)=\{f\in L^2(E,m;\mathbb{C}): (e^{-\overline{\lambda_k}}I-\widehat{T}_1)^nf=0\}.$$
We have
\begin{equation}\label{7.4}
  (e^{-\lambda_k}I-T_1)^n=e^{-n\lambda_k}I-\sum_{j=1}^ne^{-(n-j)\lambda_k}T_1^j.
\end{equation}
Since $\sum_{j=1}^ne^{-(n-j)\lambda_k}T_1^j$ is also a compact operator,
by \cite[Theorem 6.6.13]{BP}, $\widehat{\N}_{k,n}$ is of the same dimension as $\N_{k,n}$.
In particular, $\dim(\widehat{\N}_{k,\nu_k})=\dim(\N_{k,\nu_k})=n_k.$
Thus we have
$$
  \widehat{\N}_{k,n} \varsubsetneqq \widehat{\N}_{k,n+1}, \quad n=0,1,\cdots,\nu_k-1;\quad
  \widehat{\N}_{k,n} = \widehat{\N}_{k,n+1},\quad n\ge \nu_k.
$$
Similarly, we can get, for all $k\in\I$ and $n\ge 0$, $\widehat{\N}_{k,n}$ is an invariant subspace of $\widehat{T}_t$. Hence, $\{\widehat{T}_t|_{\widehat{\N}_{k,\nu_k}},t>0\}$ is a semigroup on $\widehat{\N}_{k,\nu_k}$ with
infinitesimal generator $\widehat{\A}_{k}$.

Let $\{\widehat{\psi}_1^{(k)},\widehat{\psi}_2^{(k)},\cdots,\widehat{\psi}_{n_k}^{(k)}\}$ be a basis of $\widehat{\N}_{k,\nu_k}$ such that
\begin{equation}\label{1.14}
  \widehat{T}_t(\widehat{\psi}_1^{(k)},\widehat{\psi}_2^{(k)},\cdots,\widehat{\psi}_{n_k}^{(k)})
  =(\widehat{\psi}_1^{(k)},\widehat{\psi}_2^{(k)},\cdots,\widehat{\psi}_{n_k}^{(k)})\widehat{D}_k(t),
\end{equation}
where $\widehat{D}_k(t)$ is an $n_k\times n_k$ invertible matrix.
Since $\widehat{T}_t(\widehat{\psi}_1^{(k)},\widehat{\psi}_2^{(k)},\cdots,\widehat{\psi}_{n_k}^{(k)})(x)$ is continuous, we can choose $(\widehat{\psi}_1^{(k)},\widehat{\psi}_2^{(k)},\cdots,\widehat{\psi}_{n_k}^{(k)})$ to be continuous.
We define an $n_k\times n_k$ matrix $\widetilde{A}_k$ by
\begin{equation}\label{A_k}
  (\widetilde{A}_k)_{j, l}:=\langle \phi_j^{(k)},
\widehat{\psi}_l^{(k)}\rangle_m.
\end{equation}

\begin{lemma}\label{lemma1.1}
For each $k\in\I$,
\begin{equation}\label{1.15}
  L^2(E,m;\mathbb{C})=\N_{k,\nu_k}\oplus(\widehat{\N}_{k,\nu_k})^\bot=\widehat{\N}_{k,\nu_k}\oplus(\N_{k,\nu_k})^\bot.
\end{equation}
Morover, the matrix $\widetilde{A}_k$ defined in \eqref{A_k} is invertible.
\end{lemma}
\textbf{Proof:}
By \cite[Theorem 6.6.7]{BP}, we have
$L^2(E,m;\mathbb{C})=\N_{k,\nu_k}\oplus\R_{k,\nu_k}$.
It follows from  \cite[Theorem 6.6.14]{BP} that
$\R_{k,\nu_k}=(\widehat{\N}_{k,\nu_k})^\bot.$
Thus,
$L^2(E,m;\mathbb{C})=\N_{k,\nu_k}\oplus(\widehat{\N}_{k,\nu_k})^\bot$.
Similarly, we have $L^2(E,m;\mathbb{C})=\widehat{\N}_{k,\nu_k}\oplus(\N_{k,\nu_k})^\bot$.

For any vector $a=(a_1,\cdots, a_{n_k})^T\in \mathbb{C}^{n_k}$, we have
$$
\widetilde{A}_ka=(\langle \phi_1^{(k)},h\rangle_m,\langle \phi_2^{(k)},h\rangle_m,\cdots,\langle \phi_{n_k}^{(k)},h\rangle_m)^T,
$$
where
$h=(\widehat{\psi}_1^{(k)},\widehat{\psi}_2^{(k)},\cdots,\widehat{\psi}_{n_k}^{(k)})\bar{a}\in \widehat{\N}_{k,\nu_k} .$

If $\widetilde{A}_ka=0$, then $h\in (\N_{k,\nu_k})^\bot.$
Since $\widehat{\N}_{k,\nu_k}\cap (\N_{k,\nu_k})^\bot=\{0\}$,  we have $h=0,$ which implies $a=0.$
Therefore, $\widetilde{A}_k$ is invertible.
\hfill$\Box$

\begin{lemma}\label{lemma T*}
For any $k\in\I$, define
$$
(\Psi_k(x))^T:=\left(\psi_1^{(k)}(x),\psi_2^{(k)}(x),\cdots,\psi_{n_k}^{(k)}(x)\right)
:=\left(\widehat{\psi}_1^{(k)}(x),\widehat{\psi}_2^{(k)}(x),\cdots,\widehat{\psi}_{n_k}^{(k)}(x)\right)
\overline{\widetilde{A}_k^{-1}}.
$$
Then $\{\psi_1^{(k)},\psi_2^{(k)},\cdots,\psi_{n_k}^{(k)}\}$ is a basis of $\widehat{\N}_{k,\nu_k}$ such
that the $n_k\times n_k$ matrix $A_k:=(\langle \phi_j^{(k)}, \psi_l^{(k)}\rangle_m)$ satisfies
\begin{equation}\label{A}
  A_k=I
\end{equation}
and for any $x\in \E$,
\begin{equation}\label{T^*}
  \widehat{T}_t(\Psi_k)(x)=e^{-\overline{\lambda_k} t}D_k(t)\Psi_k(x).
  \end{equation}
Moreover, the basis of $\widehat{\N}_{k,\nu_k}$ satisfying \eqref{A} is unique.
\end{lemma}
\textbf{Proof:}
For any $\mathbb{C}^{n}$-valued functions $(f_1(x), f_2(x), \cdots f_{n}(x))^T$ and
$(g_1(x), g_2(x), \cdots g_{n}(x))^T$, we use $\langle (f_1, f_2, \cdots f_{n}), (g_1, g_2, \cdots g_{n}) \rangle_m$ to denote the $n\times n$ matrix $(\langle f_j, g_l\rangle_m)$.
Since $\overline{\widetilde{A}_k^{-1}}$ is invertible, $\{\psi_1^{(k)},\psi_2^{(k)},\cdots,\psi_{n_k}^{(k)}\}$ is a basis of $\widehat{\N}_{k,\nu_k}$.
By \eqref{T-Jordan} and \eqref{1.14}, we get
\begin{eqnarray*}
  &&e^{-\lambda_k t}(D_k(t))^T\widetilde{A}_k=\left\langle T_t\left(\phi_1^{(k)},\phi_2^{(k)},\cdots,\phi_{n_k}^{(k)}\right),
\left(\widehat{\psi}_1^{(k)},\widehat{\psi}_2^{(k)},\cdots,\widehat{\psi}_{n_k}^{(k)}\right)\right\rangle_m \\
&=& \left\langle \left(\phi_1^{(k)},\phi_2^{(k)},\cdots,\phi_{n_k}^{(k)}\right),
\widehat{T}_t\left(\widehat{\psi}_1^{(k)},\widehat{\psi}_2^{(k)},\cdots,\widehat{\psi}_{n_k}^{(k)}\right)\right\rangle_m
  =\widetilde{A}_k\overline{\widehat{D_k}(t)}.
\end{eqnarray*}
Since $D_k(t)$ is a real matrix, we have
\begin{equation}\label{1.17}
  e^{-\overline{\lambda_k} t}\overline{\widetilde{A}_k^{-1}}(D_k(t))^T=\widehat{D_k}(t)\overline{\widetilde{A}_k^{-1}}.
\end{equation}
By \eqref{1.14} and \eqref{1.17}, we have
\begin{eqnarray*}
&&\widehat{T}_t\left(\psi_1^{(k)},\psi_2^{(k)},\cdots,\psi_{n_k}^{(k)}\right)
=\left(\widehat{\psi}_1^{(k)},\widehat{\psi}_2^{(k)},\cdots,\widehat{\psi}_{n_k}^{(k)}\right)\widehat{D_k}(t)
\overline{\widetilde{A}_k^{-1}}\\
&=&e^{-\overline{\lambda_k} t}\left(\widehat{\psi}_1^{(k)},\widehat{\psi}_2^{(k)},\cdots,\widehat{\psi}_{n_k}^{(k)}\right)
\overline{\widetilde{A}_k^{-1}}(D_k(t))^T
=e^{-\overline{\lambda_k} t}\left(\psi_1^{(k)},\psi_2^{(k)},\cdots,\psi_{n_k}^{(k)}\right)(D_k(t))^T.
\end{eqnarray*}
Assume that there exists another basis $\widetilde{\Psi}_k(x)$ of $\widehat{\N}_{k,\nu_k}$ satisfying \eqref{A}.
Then there exists  matrix $B$ such that
$(\widetilde{\Psi}_k(x))^T=(\Psi_k(x))^TB$.
Thus,
$$I=\langle (\Phi_k)^T,(\widetilde{\Psi}_k)^T\rangle_m=\langle(\Phi_k)^T,(\Psi_k)^T\rangle_m\overline{B}=\overline{B},$$
which implies $B=I$. Thus, we get $\widetilde{\Psi}_k(x)=\Psi_k(x)$.
The proof is now complete.
\hfill$\Box$

\begin{remark}\label{rek5}
We know that $T_t(\overline{\Phi_k^T})(x)=e^{-\overline{\lambda_k}t}\overline{\Phi^T_j(x)}D_k(t)$. Thus
$e^{-\overline{\lambda_k}t}$ is also a eigenvalue of $T_t$.
Hence there exists a unique $k'$ such that $\lambda_{k'}=\overline{\lambda_k}$. It is obvious that $D_k(t)=D_{k'}(t)$
and we can choose $\Phi_{k'}(x)=\overline{\Phi_{k}(x)}$.
By Lemma \ref{lemma T*}, we have $\Psi_{k'}(x)=\overline{\Psi_{k}(x)}$.
In particular, if $\lambda_k$ is  real, then $k'=k$.
\end{remark}

\begin{lemma}\label{lemma1.2}
For $j,k\in\I$ and $j\ne k$ ,
we have

\begin{equation}\label{1.18}
  \N_{j,\nu_j}\subset\R_{k,\nu_k}=(\widehat{\N}_{k,\nu_k})^\bot.
\end{equation}
In particular, $\N_{j,\nu_j}\cap \N_{k,\nu_k}=\{0\}$.
\end{lemma}
\textbf{Proof:}
Assume $f\in \N_{j,\nu_j}$, then $(e^{-\lambda_j}I-T_1)^{\nu_j}f=0$.
Since $\nu_j\ge 1$, we can define $g=(e^{-\lambda_j}I-T_1)^{\nu_j-1}f$.
Thus
$e^{-\lambda_j}g=T_1g.$
Hence, $(e^{-\lambda_k}I-T_1)g=(e^{-\lambda_k}-e^{-\lambda_j})g$, which implies
$$
(e^{-\lambda_k}I-T_1)^{\nu_k}g=(e^{-\lambda_k}-e^{-\lambda_j})^{\nu_k}g.
$$
Therefore
$g=(e^{-\lambda_k}-e^{-\lambda_j})^{-\nu_k}(e^{-\lambda_k}I-T_1)^{\nu_k}g\in \R_{k,\nu_k}.$

Assume $f=f_1+f_2$ with $f_1\in\N_{k,\nu_k}$ and $f_2\in\R_{k,\nu_k}.$
Then  $(e^{-\lambda_j}I-T_1)^{\nu_j-1}f_1\in\N_{k,\nu_k}$.
On the other hand,
$(e^{-\lambda_j}I-T_1)^{\nu_j-1}f_1=g-(e^{-\lambda_j}I-T_1)^{\nu_j-1}f_2\in\R_{k,\nu_k}.$
Thus
$(e^{-\lambda_j}I-T_1)^{\nu_j-1}f_1=0.$

If $\nu_j=1$, then $f=g\in\R_{k,\nu_k}$.
If $\nu_j>1$ and $f_1\ne 0$, then $e^{-\lambda_j}\in \sigma(T_1|_{\N_{k,\nu_k}})$.
By \cite[Theorem 6.7.4]{BP},
$\sigma(T_1|_{\N_{k,\nu_k}})=\{e^{-\lambda_k}\}$.
This is a contradiction. Thus, $f_1=0$, which implies $f=f_2\in\R_{k,\nu_k}$.
Therefore $ \N_{j,\nu_j}\subset\R_{k,\nu_k}.$
\hfill$\Box$

By Lemma \ref{lemma1.2}, for $k\in\I$, we can define
$$
\mathcal{M}_k:=\N_{1,\nu_1}\oplus\N_{2,\nu_2}\oplus\cdots\oplus\N_{k,\nu_k}
\mbox{ and }
\widehat{\mathcal{M}}_k:=\widehat{\N}_{1,\nu_1}\oplus\widehat{\N}_{2,\nu_2}\oplus\cdots\oplus\widehat{\N}_{k,\nu_k}.
$$

\begin{cor}\label{cor2}
For any $k\in\I$,
\begin{eqnarray}\label{2.1}
  L^2(E,m;\mathbb{C})=\mathcal{M}_k\oplus(\widehat{\mathcal{M}}_k)^\bot=\widehat{\mathcal{M}}_k\oplus(\mathcal{M}_k)^\bot.
\end{eqnarray}
\end{cor}
\textbf{Proof:} By \eqref{1.15}, \eqref{2.1} holds for $k=1$.
Assume that \eqref{2.1} holds for $k-1$.
Then
\begin{equation}\label{2.9}
  L^2(E,m;\mathbb{C})=\mathcal{M}_{k-1}\oplus(\widehat{\mathcal{M}}_{k-1})^\bot.
\end{equation}
For any $f\in(\widehat{\mathcal{M}}_{k-1})^\bot$, by \eqref{1.15},
we have $f=f_3+f_4$, where $f_3\in \N_{k,\nu_k}$ and $f_4\in(\widehat{\N}_{k,\nu_k})^\bot$.
By \eqref{1.18},
$f_3\in \bigcap_{j=1}^{k-1}(\widehat{\N}_{j,\nu_j})^\bot=(\widehat{\mathcal{M}}_{k-1})^\bot$,
which implies $f_4=f-f_3\in (\mathcal{M}_{k-1})^\bot$.
Thus, we obtain
$$
f_4\in (\widehat{\N}_{k,\nu_k})^\bot \cap(\mathcal{M}_{k-1})^\bot=(\mathcal{M}_k)^\bot.
$$
Hence
$$(\mathcal{M}_{k-1})^\bot=\N_{k,\nu_k}\oplus (\mathcal{M}_{k})^\bot.$$
Therefore, by induction, the first part of \eqref{2.1} holds for all $k\in\I$.

The proof of  $L^2(E,m;\mathbb{C})=\widehat{\mathcal{M}}_k\oplus(\mathcal{M}_k)^\bot$
is similar.
\hfill$\Box$

\begin{remark}\label{rek2}
Since $-\lambda_1$ is simple, which means $n_1=r_1=\nu_1=1$,
we know that $\Phi_1(x)=\phi_1(x)$ and $\Psi_1(x)=\psi_1(x)$.
Moreover, since $T_t\phi_1(x)=e^{-\lambda_1 t}\phi_1(x)$
and $\widehat{T}_t\psi_1(x)=e^{-\lambda_1t}\psi_1(x)$ for every $x$, $\phi_1$ and $\psi_1$ are continuous and  strictly positive.
It is easy to see that $D_1(t)\equiv1$.

By Lemma \ref{lemma1.2},  $\{\phi_{l}^{(j)},j=1,\cdots,k,l=1,\cdots,n_j\}$ is a basis of $\mathcal{M}_k$
and $\{\psi_{l}^{(j)},j=1,\cdots,k,l=1,\cdots,n_j\}$ is a basis of $\widehat{\mathcal{M}}_k$.
By \eqref{1.18} and \eqref{A}, we get
$\langle\phi_l^{(j)},\psi_n^{(k)}\rangle_m=1$, when $j=k$ and $l=n$;
otherwise $\langle\phi_l^{(j)},\psi_n^{(k)}\rangle_m=0.$
\end{remark}

In this paper, we always assume that the branching Markov process $X$ is supercritical, that is,
\begin{assumption}
 $\lambda_1<0$.
 \end{assumption}

We will use $\{{\cal F}_t: t\ge0\}$ to denote the filtration of $X$, that is
${\cal F}_t=\sigma(X_s: s\in [0, t])$.
Using the expectation formula of $\langle \phi_1, X_t\rangle$ and the Markov property of $X$,
it is easy to show that (see Lemma \ref{thrm1}),
for any nonzero $\nu\in {\cal M}_a(E)$, under $\P_{\nu}$,
the process $W_t:=e^{\lambda_1 t}\langle \phi_1, X_t\rangle$ is a positive martingale.
Therefore it converges:
$$
  W_t \to W_\infty,\quad \P_{\nu}\mbox{-a.s.} \quad \mbox{ as }t\to \infty.
$$
Using the assumption \eqref{1.16}
we can show that, as $t\to \infty$, $W_t$ also converges in $L^2(\P_{\nu})$,
so $W_\infty$ is non-degenerate and
 the second moment is finite. Moreover, we have $\P_{\nu}(W_\infty)=\langle\phi_1, \nu\rangle$.
Put $\mathcal{E}=\{W_\infty=0\}$, then $\P_{\nu}(\mathcal{E})<1$.
It is clear that  $\mathcal{E}^c\subset\{X_t(E)>0,\forall t\ge 0\}$.

\subsection{Main results}

For any $k\in\I$,
every function $f\in L^2(E,m;\mathbb{C})$ can be written uniquely as the sum of a function $f_k\in  \mathcal{M}_k$
and a function in $(\widehat{\mathcal{M}}_k)^\bot$. Similarly, every function $f\in L^2(E,m;\mathbb{C})$ can be written uniquely as the sum of a function $\widehat{f}_k\in  \widehat{\mathcal{M}}_k$
and a function in $(\mathcal{M}_k)^\bot$. Using \eqref{A}, we can easily get that
\begin{equation}\label{fk}
  f_k(x)=\sum_{j=1}^k(\Phi_j(x))^T\langle f,\Psi_j\rangle_m\in \mathcal{M}_k\quad \mbox{and}\quad
\widehat{f}_k(x)=\sum_{j=1}^k(\Psi_j(x))^T\langle f,\Phi_j\rangle_m\in \widehat{\mathcal{M}}_k ,
\end{equation}
where
$$
\langle f,\Psi_j\rangle_m:=(\langle f,\psi^{(j)}_1\rangle_m,\langle f,\psi^{(j)}_2\rangle_m,\cdots,\langle f,\psi^{(j)}_{n_j}\rangle_m)^T
$$
and
$$
\langle f,\Phi_j\rangle_m:=(\langle f,\phi^{(j)}_1\rangle_m,\langle f,\phi^{(j)}_2\rangle_m,\cdots,\langle f,\phi^{(j)}_{n_j}\rangle_m)^T.
$$
For any $f\in L^2(E,m;\mathbb{C})$, we define
$$
\gamma(f):=\inf\{j\in\I:\langle f,\Psi_j\rangle_m\ne 0\},
$$
where we use the usual convention that $\inf\varnothing=\infty$.
If $\gamma(f)<\infty$, define
$$
\zeta(f):=\sup\{j\in\I:\Re_j=\Re_{\gamma(f)}\}.
$$
For each $j\in\I$, every component of the function $t:\rightarrow D_j(t)\langle f,\Psi_j\rangle_m$ is a polynomial of $t$.
Denote the degree of the $l$-th component of $D_j(t)\langle f,\Psi_j\rangle_m$ by  $\tau_{j,l}(f).$
We define
$$\tau(f):=\sup\{\tau_{j,l}(f): \gamma(f)\le j\le \zeta(f),1\le l\le n_j\}.$$
Then for any $j$ with $\Re_j=\Re_{\gamma(f)}$,
\begin{equation}\label{1.20}
F_{f,j}:=\lim_{t\to\infty} t^{-\tau(f)}D_j(t)\langle f,\Psi_j\rangle_m
\end{equation}
exists and there exists a $j$ such that $F_{f,j}\ne 0$.

Note that if $g\in L^2(E,m)$, then for any $j\in\I$,
$$\overline{\langle g,\Psi_j\rangle_m}=\langle g,\overline{\Psi_j}\rangle_m=\langle g,\Psi_{j'}\rangle_m,$$
where $j'$ is defined in Remark \ref{rek5}.
For $g(x)=\sum_{k:\lambda_1\ge 2\Re_k} (\Phi_k(x))^Tb_k$, we get
$b_k=\langle g,\Psi_j\rangle_m.$
Thus, if $g(x)$ is real, we get $\overline{b_k}=b_{k'}$.
The following three subsets of $L^2(E, m)$ will be needed in the statement of our main result:
$$
\C_l:=\left\{g(x)=\sum_{k\in\I:\lambda_1>2\Re_k} (\Phi_k(x))^Tb_k: b_k\in \mathbb{C}^{n_k} \mbox{ with } \overline{b_k}=b_{k'}\right\},
$$
$$
\C_c:=\left\{g(x)=\sum_{k\in\I:\lambda_1=2\Re_k} (\Phi_k(x))^Tb_k: b_k\in \mathbb{C}^{n_k} \mbox{ with } \overline{b_k}=b_{k'}\right\},
$$
and
$$
\C_s:=\left\{g\in L^2(E,m)\cap L^4(E,m):\lambda_1< 2\Re_{\gamma(g)}\right\}.
$$

\subsubsection{Some basic law of large numbers}

For any $k\in \I$, we define an $n_k$-dimensional random vector $H^{(k)}_t$ as follows:
$$
H_{t}^{(k)}:=e^{\lambda_kt}(\langle \phi_1^{(k)}, X_t\rangle, \cdots,\langle\phi_{n_k}^{(k)},X_t\rangle) (D_k(t))^{-1}.
$$
One can show (see Lemma \ref{thrm1} below) that, if $\lambda_1>2\Re_k$,
then, for any $\nu\in {\cal M}_a(E)$ and $b\in\mathbb{C}^{n_k}$, $H_t^{(k)}b$ is a martingale under $\P_{\nu}$ and bounded in
 $L^2(\P_{\nu})$.
Thus the limit $H_\infty^{(k)}:=\lim_{t\to\infty}H_t^{(k)}$ exists $\P_{\nu}$-a.s. and in $L^2(\P_{\nu})$.

\begin{thrm}\label{thrm2}
If $f\in L^2(E,m;\mathbb{C})\cap L^4(E,m;\mathbb{C})$ with $\lambda_1>2\Re_{\gamma(f)}$,
then for any nonzero $\nu\in {\cal M}_a(E)$, as $t\to\infty$,
$$
t^{-\tau(f)}e^{\Re_{\gamma(f)}t}\langle f,X_t\rangle-\sum_{j=\gamma(f)}^{\zeta(f)}e^{-i\Im_j t}H_\infty^{(j)}F_{f,j}\to 0,
\quad \mbox{ in } L^2(\P_{\nu}).
$$
\end{thrm}
\begin{remark}\label{rem:large}
Suppose $f\in L^2(E,m;\mathbb{C})\cap L^4(E,m;\mathbb{C})$ with $\gamma(f)=1$.
Then $\zeta(f)=1$.
Since $D_1(t)\equiv1$, $\tau(f)=0$.
Thus $H_t^{(1)}$ reduces to $W_t$ and $H_\infty^{(1)}=W_\infty$.
Therefore by Theorem \ref{thrm2} and the fact that $F_{f,1}=\langle f,\psi_1\rangle_m$,
we get that for any nonzero $\nu\in{\cal M}_a(E)$,
\begin{equation*}
  e^{\lambda_1 t}\langle f,X_t\rangle\to \langle f,\psi_1\rangle_m W_\infty, \quad \mbox{in }
  L^2(\P_{\nu}),
\end{equation*}
as $t\to\infty$.
It is obvious that the convergence also holds in $\P_{\nu}$-probability.

In particular, if  $f$ is non-zero and non-negative, then $\langle f,\psi_1\rangle_m\ne 0$
which implies $\gamma(f)=1$.\hfill$\Box$
\end{remark}

\subsubsection{Main result}
For $f\in \C_s$, define
\begin{equation}\label{e:sigma}
\sigma_{f}^2:=\int_0^\infty e^{\lambda_1 s}\langle A|T_s f|^2,\psi_1\rangle_m \,ds+\langle |f|^2,\psi_1\rangle_m.
\end{equation}
For $h=\sum_{k:\lambda_1=2\Re_k} (\Phi_k(x))^Tb_k\in\C_c$, define
\begin{equation}\label{e:rho}
\rho_{h}^2:=(1+2\tau(h))^{-1}\left\langle AF_{h},\psi_1\right\rangle_m,
\end{equation}
where
$F_{h}(x):=\sum_{k:\lambda_1=2\Re_k}\left|(\Phi_k(x))^TF_{h,k}\right|^2$.
For $g(x)=\sum_{k: \lambda_1>2\Re_k}\left(\Phi_k(x)\right)^Tb_k\in \C_l$, define
$$
I_sg(x):=\sum_{k: \lambda_1>2\Re_k}e^{\lambda_ks}\Phi_k(x)^TD_k(s)^{-1}b_k, \quad \beta_g^2
:=\int_{0}^\infty e^{-\lambda_1u}\langle A\left|I_ug\right|^2,\psi_1\rangle_m\,du-\langle  g^2,\psi_1\rangle_m
$$
and
$$
E_t(g):=\sum_{k: \lambda_1>2\Re_k}\left(e^{-\lambda_kt}H^{(k)}_\infty D_k(t)b_k\right).
$$

\begin{thrm}\label{The:1.3}
If $f\in \C_s$, $h\in\C_c$ and
 $g\in\C_l$, then
$\sigma_f^2$, $\rho_h^2$ and $\beta_g^2$ all belong to $(0, \infty)$.
Furthermore, it holds that, under $\P_{\nu}(\cdot\mid \mathcal{E}^c)$, as $t\to\infty$,
\begin{eqnarray}\label{result}
   &&\left(e^{\lambda_1 t}\langle \phi_1, ,X_t\rangle,
 ~\frac{ \langle g,X_t\rangle-E_t(g)}
  {\sqrt{\langle \phi_1,X_t\rangle}},
  ~\frac{\langle h , X_t\rangle}{\sqrt{t^{1+2\tau(h)}\langle \phi_1,X_t\rangle}},
  ~\frac{\langle f , X_t\rangle}{\sqrt{\langle \phi_1,X_t\rangle}} \right) \nonumber\\
 && \stackrel{d}{\rightarrow}(W^*,G_3(g),G_2(h),~G_1(f)),
\end{eqnarray}
where $W^*$ has the same distribution as $W_\infty$ conditioned on $\mathcal{E}^c$,
$G_3(g)\sim \mathcal{N}(0,\beta_g^2)$, $G_2(h)\sim \mathcal{N}(0,\rho_h^2)$
and $G_1(f)\sim \mathcal{N}(0,\sigma_f^2)$. Moreover, $W^*$, $G_3(g)$, $G_2(h)$ and $G_1(f)$ are independent.
\end{thrm}

Whenever $f\in \C_s$, we will use $G_1(f)$ to denote a normal random variable $\mathcal{N}(0,\sigma_f^2)$.
For $f_1,f_2\in \C_s$, define
$$
\sigma(f_1,f_2):=\int_0^\infty e^{\lambda_1 s}\langle A(T_s f_1)(T_sf_2),\psi_1\rangle_m \,ds+\langle f_1f_2,\psi_1\rangle_m.
$$
\begin{cor}\label{Cor:1}
If $f_1,f_2\in \C_s$,
then, under $\P_{\nu}(\cdot\mid \mathcal{E}^c)$,
$$
 \left(\frac{\langle f_1 , X_t\rangle}{\sqrt{\langle \phi_1,X_t\rangle}},
 \frac{\langle f_2, X_t\rangle}{\sqrt{\langle \phi_1,X_t\rangle}} \right)
 \stackrel{d}{\rightarrow}(G_1(f_1),G_1(f_2)), \quad t\to\infty,
$$
and $(G_1(f_1),G_1(f_2))$ is a bivariate normal random variable with covariance
\begin{equation}\label{sigma(fg)}
  \mbox{\rm Cov}(G_1(f_1),G_1(f_2))=\sigma(f_1,f_2).
\end{equation}
\end{cor}
\textbf{Proof:}
Using the convergence of the fourth component in Theorem \ref{The:1.3}, we get
\begin{eqnarray*}
  &&{\P}_\nu \left(\exp\left\{i\theta_1 \frac{\langle f_1,X_t\rangle}{\sqrt{\langle \phi_1,X_t\rangle}}
  +i\theta_2 \frac{\langle f_2,X_t\rangle}{\sqrt{\langle \phi_1,X_t\rangle}}\right\}\mid \mathcal{E}^c\right)\\
  &=&{\P}_\nu \left(\exp\left\{i\frac{\langle \theta_1 f_1+\theta_2 f_2,X_t\rangle}{\sqrt{\langle \phi_1,X_t\rangle}}\right\} \mid \mathcal{E}^c\right) \\
   &\to& \exp\left\{ -\frac{1}{2}\sigma_{(\theta_1f_1+\theta_2f_2)}^2\right\},\quad \mbox{as}\quad t\to\infty,
\end{eqnarray*}
where
\begin{eqnarray*}
 \sigma_{(\theta_1f_1+\theta_2f_2)}^2 &=& \int_0^\infty e^{\lambda_1 s}\langle A(T_s (\theta_1f_1+\theta_2f_2))^2,\psi_1\rangle_m \,ds
+\langle (\theta_1f_1+\theta_2f_2)^2,\psi_1\rangle_m \\
   &=&  \theta_1^2\sigma_{f_1}^2+2\theta_1\theta_2\sigma(f_1,f_2)+\theta_2^2\sigma_{f_2}^2.
\end{eqnarray*}
Now \eqref{sigma(fg)} follows immediately.\hfill$\Box$

Whenever $h\in \C_c$, we will use $G_2(h)$ to denote a normal random variable $\mathcal{N}(0,\rho_h^2)$.
For $h_1, h_2\in\C_c$, define
\begin{equation}\label{rho2}
  \rho(h_1,h_2):=(1+\tau(h_1)+\tau(h_2))^{-1}\left\langle AF_{h_1,h_2},\psi_1\right\rangle_m,
\end{equation}
where
\begin{equation}\label{e:Ffg}
F_{h_1,h_2}(x):=\sum_{j:\lambda_1=2\Re_j}\Phi_j(x)^TF_{h_1,j}\Phi_{j'}(x)^TF_{h_2,j'}
=\sum_{j:\lambda_1=2\Re_j}\Phi_j(x)^TF_{h_1,j}\overline{\Phi_{j}(x)^TF_{h_2,j}}.
\end{equation}

\begin{cor}\label{cor:2}
If $h_1, h_2\in\C_c$,
then we have, under $\P_{\nu}(\cdot\mid \mathcal{E}^c)$,
$$
 \left(\frac{\langle h_1 , X_t\rangle}{\sqrt{t^{1+2\tau(h_1)}\langle \phi_1,X_t\rangle}}, \frac{\langle h_2, X_t\rangle}{\sqrt{t^{1+2\tau(h_2)}
\langle \phi_1,X_t\rangle}} \right)\stackrel{d}{\rightarrow}(G_2(h_1),G_2(h_2)), \quad t\to\infty,
$$
and $(G_2(h_1),G_2(h_2))$ is a bivariate normal random variable with covariance
$$
\mbox{\rm Cov}(G_2(h_1),G_2(h_2))=\rho(h_1,h_2).
$$
\end{cor}

Whenever $g\in \C_l$, we will use $G_3(g)$ to denote a normal random variable $\mathcal{N}(0,\beta_g^2)$.
For $g_1(x) ,g_2(x)\in\C_l$,
define
$$
\beta(g_1,g_2):=\int_0^\infty e^{-\lambda_1s}\langle A(I_sg_1)(I_sg_2),\psi_1\rangle_m\,ds-\langle g_1g_2,\psi_1\rangle_m.
$$
Using the convergence of the second component in Theorem \ref{The:1.3}
and an argument similar to that in the proof of Corollary \ref{Cor:1}, we get
\begin{cor}\label{cor:3}
If $g_1(x) ,g_2(x)\in\C_l$,
then we have, under $\P_{\nu}(\cdot\mid \mathcal{E}^c)$,
\begin{eqnarray*}
&&\left(\frac{\langle g_1,X_t\rangle -E_t(g_1)}
{\sqrt{\langle \phi_1,X_t\rangle}},
\frac{\langle g_2,X_t\rangle -E_t(g_2)}
{\sqrt{\langle \phi_1,X_t\rangle}}\right)\stackrel{d}{\rightarrow}(G_3(g_1),G_3(g_2)),
\end{eqnarray*}
and $(G_3(g_1),G_3(g_2))$ is a bivariate normal random variable with covariance
$$
\mbox{\rm Cov}(G_3(g_1),G_3(g_2))=\beta(g_1,g_2).
$$
\end{cor}

For any
$f\in L^2(E, m)\cap L^4(E,m)$, define
\begin{eqnarray*}
f_{(s)}(x)&:=&\sum_{j: 2\Re_j<\lambda_1}(\Phi_j(x))^T\langle f,\Psi_j\rangle_m,
\\
f_{(c)}(x)&:=&\sum_{j:2\Re_j=\lambda_1} (\Phi_j(x))^T\langle f,\Psi_j\rangle_m,
\\
f_{(l)}(x)&:=&f(x)-f_{(s)}(x)-f_{(l)}(x).
\end{eqnarray*}
Then $f_{(s)}\in \C_l$, $f_{(c)}\in \C_c$ and $f_{(l)}\in \C_s$.

\begin{remark}\label{r:critical}
If $f\in L^2(E,m)\cap L^4(E,m)$ with $\lambda_1=2\Re_{\gamma(f)}$, then $f=f_{(c)}+f_{(l)}$.
Using the convergence of the fourth component in Theorem \ref{The:1.3} for $f_{(l)}$, it holds under
$\P_{\nu}(\cdot\mid \mathcal{E}^c)$ that
$$
\frac{\langle f_{(l)},X_{t}\rangle}{\sqrt{t^{1+2\tau(f)}\langle \phi_1,X_t\rangle}} \stackrel{d}{\to} 0,\quad t\to \infty.
$$
Thus using the convergence of the first and third components in Theorem \ref{The:1.3}, we get, under $\P_{\nu}(\cdot\mid \mathcal{E}^c)$,
$$
  \left(e^{\lambda_1 t}\langle \phi_1, X_t\rangle, ~\frac{\langle f , X_t\rangle}
  {\sqrt{t^{1+2\tau(f)}\langle \phi_1,X_t\rangle}} \right)
  \stackrel{d}{\rightarrow}(W^*,~G_2(f_{(c)})), \quad t\to\infty,
$$
where $W^*$ has the same distribution as $W_\infty$ conditioned on $\mathcal{E}^c$
and $G_2(f_{(c)})\sim \mathcal{N}(0,\rho_{f_{(c)}}^2)$. Moreover, $W^*$ and $G_2(f_{(c)})$ are independent.
\end{remark}

\begin{remark}\label{large}
Assume $f\in L^2(E,m)\cap L^4(E,m)$ satisfies
$\lambda_1>2\Re_{\gamma(f)}$.

If $f_{(c)}=0$, then $f=f_{(l)}+f_{(s)}$.
Using the convergence of the first, second and fourth components in Theorem \ref{The:1.3},
we get for any nonzero $\nu\in {\cal M}_a(E)$, it holds under $\P_{\nu}(\cdot\mid \mathcal{E}^c)$ that, as $t\to\infty$,
$$
  \left(e^{\lambda_1 t}\langle \phi_1,X_t\rangle,~\frac
  {\left(\langle f,X_t\rangle-\sum_{k: 2\Re_k<\lambda_1}
  e^{-\lambda_kt}H^{(k)}_\infty D_k(t)\langle f, \Psi_k\rangle_m\right)}
  {\langle \phi_1,X_t\rangle^{1/2}} \right)
  \stackrel{d}{\rightarrow}(W^*,~G_1(f_{(l)})+G_3(f_{(s)})),
$$
where $W^*$, $G_3(f_{(s)})$ and $G_1(f_{(l)})$ are the same as those in Theorem \ref{The:1.3}.
Since $G_3(f_{(s)})$ and $G_1(f_{(l)})$ are independent,
$G_1(f_{(l)})+G_3(f_{(s)})\sim\mathcal{N}\left(0,\sigma_{f_{(l)}}^2+\beta_{f_{(s)}}^2\right)$.

If $f_{(c)}\ne 0$, then as $t\to\infty$,
$$
\frac{\left(\langle f_{(l)}+f_{(s)},X_t\rangle-\sum_{k: 2\Re_k<\lambda_1}
  e^{-\lambda_kt}H^{(k)}_\infty\langle f, \Psi_k\rangle_m\right)}
  {\sqrt{t^{1+2\tau(f)}\langle \phi_1,X_t\rangle}}
  \stackrel{d}{\rightarrow} 0.
  $$
Then using the convergence of the first and third components in Theorem \ref{The:1.3}, we get
$$
\left(e^{\lambda_1 t}\langle \phi_1,X_t\rangle,~\frac
  {\left(\langle f,X_t\rangle-\sum_{k: 2\Re_k<\lambda_1}
  e^{-\lambda_kt}H^{(k)}_\infty D_k(t)\langle f, \Psi_k\rangle_m\right)}
  {\sqrt{t^{1+2\tau(f)}\langle \phi_1,X_t\rangle}} \right)
  \stackrel{d}{\rightarrow}(W^*,~G_2(f_{(c)})),
$$
where $W^*$ and $G_2(f_{(c)})$ are the same as those in Remark \ref{r:critical}.
Thus \cite[Theorem 1.13]{RSZ} is a consequence of Theorem \ref{The:1.3}.
\end{remark}

\section{Estimates on the moments of $X$}
In the remainder of this paper we will use the  following notation:
for two positive functions $f(t,x)$ and $g(t,x)$, $f(t,x)\lesssim g(t,x)$ means that there exists a constant $c>0$ such that
$f(t,x)\le cg(t,x)$ for all $t,x$.

\subsection{Estimates on the first moment of $X$}
\begin{lemma}\label{lemma2.1}
For each $k\in\I$,
if $a<\Re_{k+1}$, there exists a  constant $c(k,a)>0$ such that for all $t>0$,
\begin{equation}\label{2.4}
  \|T_t|_{(\widehat{\mathcal{M}}_k)^\bot}\|_2\le c(k,a)e^{-at}\quad
  \mbox{and}\quad \|\widehat{T}_t|_{(\mathcal{M}_k)^\bot}\|_2\le c(k,a)e^{-at}.
\end{equation}
\end{lemma}
\textbf{Proof:}
Since $(\mathcal{M}_k)^\bot$ is invariant for $\widehat{T}_t$,
$\{\widehat{T}_t|_{(\mathcal{M}_k)^\bot}:t>0\}$ is a semigroup on $(\mathcal{M}_k)^\bot$ .
By \cite[Theorem 6.7.5]{BP}, we have
$\sigma(\widehat{T}_1|_{(\mathcal{M}_k)^\bot})=\{e^{-\overline{\lambda_{j}}},k+1\le j\in\I\}\cup \{0\}.$
Thus,
if $k+1\in\I$,
 the spectral radius of $\widehat{T}_1|_{(\mathcal{M}_k)^\bot}$ is
$r(\widehat{T}_t|_{(\mathcal{M}_k)^\bot})=e^{-\Re_{k+1}}<e^{-a}$.
If $k+1$ does not belong to $\I$, then $r(\widehat{T}_t|_{(\mathcal{M}_k)^\bot})=0<e^{-a}$.

By \cite[Theorem 6.3.10]{BP},
$r(\widehat{T}_1|_{(\mathcal{M}_k)^\bot})=\lim_{n\to\infty}(\|\widehat{T}_n|_{(\mathcal{M}_k)^\bot}\|_2)^{1/n},$
thus there exists a constant $n_1$, such that
\begin{equation}\label{2.3}
  \|\widehat{T}_{n_1}|_{(\mathcal{M}_k)^\bot}\|_2\le e^{-an_1}.
\end{equation}
By \eqref{1.66}, we have
\begin{equation}\label{2.2}
  \sup_{0\le t\le n_1}\|\widehat{T}_t|_{(\mathcal{M}_k)^\bot}\|_2
  \le  \sup_{0\le t\le n_1}\|\widehat{T}_t\|_2\le  e^{Kn_1}.
\end{equation}
For any $t>0$, there exist $l\in\mathbb{N}$ and $r\in[0,n_1)$ such that $t=n_1l+r$.
 By \eqref{2.3} and \eqref{2.2}, we have
 $$
 \| \widehat{T}_t|_{(\mathcal{M}_k)^\bot}\|_2
 \le \|\widehat{T}_{n_1}|_{(\mathcal{M}_k)^\bot}\|^l_2\|\widehat{T}_r|_{(\mathcal{M}_k)^\bot}\|_2
 \le e^{-an_1l}e^{Kn_1}\le e^{Kn_1}\left(\sup_{0\le r\le n_1}e^{ar}\right)e^{-at}.
 $$
 Thus we can find $c(k,a)>1$ such that $\|\widehat{T}_t|_{(\mathcal{M}_k)^\bot}\|_2\le c(k,a)e^{-at}$.
 Similarly, we can show that $\|T_t|_{(\widehat{\mathcal{M}}_k)^\bot}\|_2\le c(k,a)e^{-at}$.
\hfill$\Box$

\begin{lemma}\label{lemma2.2}
For each $k\in \I$ and $t_1>0$, if $a<\Re_{k+1}$,
there exists a constant $c(k,a,t_1)>0$
such  that for all $(t, x, y)\in (2t_1, \infty)\times E\times E$,
\begin{equation}\label{density}
  \left|q(t,x,y)-\sum_{j=1}^ke^{-\lambda_jt}(\Phi_j(x))^TD_j(t)\overline{\Psi_j(y)}\right|
   \le ce^{-at}b_{t_1}(x)^{1/2}\widehat{b}_{t_1}(y)^{1/2}.
\end{equation}
\end{lemma}
\textbf{Proof:}
Recall that for any $f\in L^2(E,m;\mathbb{C})$ and $k\in \I$,
$\widehat{f}_k$ is defined
in the paragraph containing \eqref{fk}.
Since $|\langle f,\phi_l^{(j)}\rangle_m|\le \|f\|_2$, we have
$|\widehat{f}_k(x)|\le \|f\|_2\sum_{j=1}^k\sum_{l=1}^{n_j}|\psi_l^{(j)}(x)|$.
Thus, we get $\|\widehat{f}_k\|_2\le c_1(k)\|f\|_2$.
By Lemma \ref{lemma2.1}, for any $a<\Re_{k+1}$, there exists a constant $c_2=c_2(k,a)>0$ such that for all $t>0$,
\begin{equation}\label{2.7}
  \|\widehat{T}_t(f-\widehat{f}_k)\|_2\le c_2e^{-at}\|f-\widehat{f}_k\|_2\le c_3e^{-at}\|f\|_2,
\end{equation}
where $c_3=c_2(1+c_1(k)).$
For $t>t_1$, we have
$$q
(t,x,y)=\int_{E}q(t_1,x,z)q(t-t_1,z,y)\,m(dz)=\widehat{T}_{t-t_1}(h_x)(y),
$$
where $h_x(z)=q(t_1,x,z)\in L^2(E,m).$
It is easy to see that
$$
\langle h_x,\phi_l^{(j)}\rangle_m=\int_E q(t_1,x,z)\overline{\phi_l^{(j)}}(z)\,m(dz)
=\overline{T_{t_1}(\phi_l^{(j)})(x)}.
$$
Let
$$
h_{x,k}(z):=\sum_{j=1}^k\sum_{l=1}^{n_j}\langle h_x,\phi_l^{(j)}\rangle_m\psi_l^{(j)}(z)
=\sum_{j=1}^k\overline{T_{t_1}((\Phi_j)^T)(x)}\Psi_j(z).
$$
By \eqref{T-Jordan} and \eqref{T^*}, we have
\begin{eqnarray*}
  \widehat{T}_{t-t_1}(h_{x,k})(y)&=&\sum_{j=1}^k\overline{T_{t_1}(\Phi_j)^T(x)}\widehat{T}_{t-t_1}(\Psi_j)(y)
  =\sum_{j=1}^ke^{-\overline{\lambda_j}t}(\overline{\Phi_j(x)})^TD_j(t_1)D_j(t-t_1)\Psi_j(y)\\
  &=&\sum_{j=1}^ke^{-\overline{\lambda_j}t}(\overline{\Phi_j(x)})^TD_j(t)\Psi_j(y).
\end{eqnarray*}
Thus, by \eqref{2.7}, we have
$$
\int_E|q(t,x,y)-\sum_{j=1}^ke^{-\overline{\lambda_j}t}(\overline{\Phi_j(x)})^TD_j(t)\Psi_j(y)|^2\,m(dy)
\le (c_3)^2e^{-2a(t-t_1)}\|h_x\|_2^2=c_4e^{-2at}b_{t_1}(x),
$$
where $c_4=c_4(k,a, t_1)=c^2_3e^{-2at_1}$.
Since $q(t,x,y)$ is a real-valued function, we have, for $t> t_1$,
\begin{equation}\label{2.11}
  \int_E|q(t,x,y)-\sum_{j=1}^ke^{-\lambda_jt}(\Phi_j(x))^TD_j(t)\overline{\Psi_j(y)}|^2\,m(dy)\le c_4e^{-2at}b_{t_1}(x).
\end{equation}

Repeating the above argument with $T_t$, we get that there exists $c_5=c_5(k,a,t_1)>0$ such that
for $t>t_1$,
\begin{equation}\label{2.15}
  \int_E|q(t,z,y)-\sum_{j=1}^ke^{-\lambda_jt}(\Phi_j(z))^TD_j(t)\overline{\Psi_j(y)}|^2\,m(dz)
\le c_5e^{-2at}\widehat{b}_{t_1}(y).
\end{equation}
Since $D_j(t)=D_{j}(t/2)D_j(t/2)$, we get
\begin{equation}\label{2.12}
  e^{-\lambda_jt}(\Phi_j(x))^TD_j(t)\overline{\Psi_j(y)}=e^{-\lambda_jt/2}\int_E q(t/2,x,z)(\Phi_j(z))^TD_j(t/2)\overline{\Psi_j(y)}\,m(dz),
\end{equation}
\begin{equation}\label{2.13}
  e^{-\lambda_jt}(\Phi_j(x))^TD_j(t)\overline{\Psi_j(y)}=e^{-\lambda_jt/2}\int_E q(t/2,z,y)(\Phi_j(x))^TD_j(t/2)\overline{\Psi_j(z)}\,m(dz),
\end{equation}
and by \eqref{A}, we have
\begin{eqnarray}
  &&\int_E \left(\sum_{j=1}^ke^{-\lambda_jt/2}(\Phi_j(x))^TD_j(t/2)\overline{\Psi_j(z)}\right)
  \left(\sum_{j=1}^ke^{-\lambda_jt/2}(\Phi_j(z))^TD_j(t/2)\overline{\Psi_j(y)}\right)\,m(dz)\nonumber  \\
   &=&\sum_{j=1}^ke^{-\lambda_jt} (\Phi_j(x))^TD_j(t/2)D_j(t/2)\overline{\Psi_j(y)}
   =\sum_{j=1}^ke^{-\lambda_jt}(\Phi_j(x))^TD_j(t)\overline{\Psi_j(y)}.\label{2.14}
\end{eqnarray}
Thus, by the semigroup property of $T_t$ and \eqref{2.12}--\eqref{2.14}, we obtain
\begin{eqnarray*}
  &&q(t,x,y)-\sum_{j=1}^ke^{-\lambda_jt}(\Phi_j(x))^TD_j(t)\overline{\Psi_j(y)}\\
  & =& \int_Eq(t/2,x,z)q(t/2,z,y)\,m(dz)
  -\sum_{j=1}^ke^{-\lambda_jt/2}\int_E q(t/2,x,z)(\Phi_j(z))^TD_j(t/2)\overline{\Psi_j(y)}\,m(dz)\\
  &&-\sum_{j=1}^ke^{-\lambda_jt/2}\int_E q(t/2,z,y)(\Phi_j(x))^TD_j(t/2)\overline{\Psi_j(z)}\,m(dz)\\
&&+ \int_E \left(\sum_{j=1}^ke^{-\lambda_jt/2}(\Phi_j(x))^TD_j(t/2)\overline{\Psi_j(z)}\right)
  \left(\sum_{j=1}^ke^{-\lambda_jt/2}(\Phi_j(z))^TD_j(t/2)\overline{\Psi_j(y)}\right)\,m(dz)\\
   &=& \int_E\left(q(t/2,x,z)-\left(\sum_{j=1}^ke^{-\lambda_jt/2}(\Phi_j(x))^TD_j(t/2)\overline{\Psi_j(z)}\right)\right)\\
    && \quad\qquad\left(q(t/2,z,y)-\left(\sum_{j=1}^ke^{-\lambda_jt/2}(\Phi_j(z))^TD_j(t/2)\overline{\Psi_j(y)}\right)\right)\,m(dz).
\end{eqnarray*}
Therefore, by H\"{o}lder's inequality, \eqref{2.11} and \eqref{2.15}, we get, for $t>2t_1$,
\begin{eqnarray*}
  &&\left|q(t,x,y)-\sum_{j=1}^ke^{-\lambda_jt}(\Phi_j(x))^TD_j(t)\overline{\Psi_j(y)}\right|
   \le \sqrt{c_4c_5}e^{-at}b_{t_1}(x)^{1/2}\widehat{b}_{t_1}(y)^{1/2}.
\end{eqnarray*}
\hfill$\Box$
\smallskip

\begin{cor}\label{lemma2.3}
Assume $f\in L^2(E,m;\mathbb{C})$.
If $\gamma(f)<\infty$,
then, for any $t_1>0$, there exists a constant $c(f,t_1)>0$
such that for all $(t,x)\in (2t_1,\infty)\times E$,
\begin{equation}\label{1.31}
   \left|t^{-\tau(f)}e^{\Re_{\gamma(f)}t}T_tf(x)-\sum_{j=\gamma(f)}^{\zeta(f)}e^{-i\Im_jt}(\Phi_j(x))^TF_{f,j}\right|
   \le c(f,t_1) t^{-1}b_{t_1}(x)^{1/2}.
\end{equation}
Moreover, we have, for $(t,x)\in(2t_1,\infty)\times E$,
\begin{equation}\label{1.23}
  |T_tf(x)|\lesssim t^{\tau(f)}e^{-\Re_{\gamma(f)}t}b_{t_1}(x)^{1/2}.
\end{equation}
If $\gamma(f)=\infty$, for any $t_1>0$,
we have, for $(t,x)\in(2t_1,\infty)\times E$,
\begin{equation}\label{1.23'}
  |T_tf(x)|\lesssim b_{t_1}(x)^{1/2}.
\end{equation}
\end{cor}
\textbf{Proof:}
First, we consider the case $\gamma(f)<\infty$, which implies $\gamma(f)\in\I$.
By the definition of $\zeta(f)$, we have $\Re_{\gamma(f)}<\Re_{\zeta(f)+1}$.
Since $\langle f, (\widehat{b}_{t_1})^{1/2}\rangle_m\le \|\widehat{b}_{t_1}^{1/2}\|_2\|f\|_2$, applying Lemma \ref{lemma2.2} with
$k=\zeta(f)$ and
a fixed $a$ with $\Re_{\gamma(f)}<a<\Re_{\zeta(f)+1}$,
we get that there exists $c_1=c_1(f,t_1)>0$
such that for $(t,x)\in(2t_1,\infty)\times E$,
\begin{equation}\label{1.21}
  \left|T_tf(x)-e^{-\Re_{\gamma(f)}t}\sum_{j=\gamma(f)}^{\zeta(f)}e^{-i\Im_jt}(\Phi_j(x))^TD_j(t)\langle f,\Psi_j\rangle_m\right|
   \le  c_1e^{-at}b_{t_1}(x)^{1/2}.
\end{equation}
If $\tau(f)\ge 1$,
the degree of each component of $D_j(t)\langle f,\Psi_j\rangle_m-t^{\tau(f)}F_{f,j}$ is no larger than $\tau(f)-1$.
Thus, for $t>2t_1$,
\begin{equation}\label{1.22}
 |D_j(t)\langle f,\Psi_j\rangle_m-t^{\tau(f)}F_{f,j}|_\infty\lesssim t^{\tau(f)-1}.
\end{equation}
If $\tau(f)=0$, $D_j(t)\langle f,\Psi_j\rangle_m-t^{\tau(f)}F_{f,j}=0.$
By \eqref{phi}, we get, for $(t,x)\in(2t_1,\infty)\times E$,
\begin{eqnarray}\label{1.24}
  &&\left|\sum_{j=\gamma(f)}^{\zeta(f)}e^{-i\Im_jt}(\Phi_j(x))^TD_j(t)\langle f,\Psi_j(y)\rangle_m-t^{\tau(f)}\sum_{j=\gamma(f)}^{\zeta(f)}e^{-i\Im_jt}(\Phi_j(x))^TF_{f,j}\right|\nonumber\\
  &\lesssim& t^{\tau(f)-1}|\Phi_j(x)|_\infty\lesssim t^{\tau(f)-1}b_{t_1}(x)^{1/2}.
\end{eqnarray}
Now  \eqref{1.31} follows easily from \eqref{1.21} and  \eqref{1.24}.
By \eqref{1.31} and \eqref{phi}, we get \eqref{1.23} immediately.

Now, we deal with the case $\gamma(f)=\infty$.
Let $k_0:=\sup\{j: \Re_j\le 0\}$. Thus, we have $k_0\in \I$ and $\Re_{k_0+1}>0$.
Since $\gamma(f)=\infty$, so for any $k\in\I$, we have $\langle f,\Psi_k\rangle_m=0$.
Now, applying Lemma \ref{lemma2.2} with $k=k_0$ and $a=0$, we get \eqref{1.23'} immediately.
\hfill$\Box$

\begin{remark}\label{rek3}
Since $D_1(t)\equiv1$,
using \eqref{density} with $k=1$ and
$\lambda_1<a<\Re_2$,
we get that, for any $t_1>0$, there exists $c_1(t_1,a)>0$ such that for any
$f\in L^2(E,m)$ and $(t,x)\in(2t_1,\infty)\times E$,
\begin{equation}\label{2.6}
  |e^{\lambda_1t}T_tf(x)-\langle f,\psi_1\rangle_m\phi_1(x)|\le c_1(t_1,a)e^{-(a-\lambda_1)t}\|f\|_2b_{t_1}(x)^{1/2},
\end{equation}
and hence there exists $c_2(t_1,a)>0$ such that
\begin{equation}\label{2.8}
  e^{\lambda_1t}|T_tf(x)|\le c_2\|f\|_2b_{t_1}(x)^{1/2}.
\end{equation}
\end{remark}

\subsection{Estimates on the second moment of $X$}

We first recall the formula for the second moment of the branching Markov process $\{X_t: t\ge 0\}$
(see, for example, \cite[Lemma 3.3]{Sh}):
for $f\in \mathcal{B}_b(E)$, we have for any $(t, x)\in (0, \infty)\times E$,
\begin{equation}\label{1.19}
   \P_{\delta_x}\langle f,X_t\rangle^2=\int_0^tT_{s}[A|T_{t-s}f|^2](x)\,ds+ T_t(f^2)(x).
\end{equation}
For any $f\in L^2(E,m)\cap L^4(E,m)$ and $x\in E$, since $(T_{t-s}f)^2(x)\le e^{K(t-s)}T_{t-s}(f^2)(x)$, we have
\begin{equation*}
  \int_0^tT_{s}[A(T_{t-s}f)^2](x)\,ds\le Ke^{Kt}T_t(f^2)(x)<\infty,
\end{equation*}
which implies
\begin{equation}\label{2.27}
  \int_0^tT_{s}[A(T_{t-s}f)^2](x)\,ds+ T_t(f^2)(x)\le (1+Ke^{Kt})T_t(f^2)(x)<\infty.
\end{equation}
Thus, using a routine limit argument, one can easily check
that \eqref{1.19} also holds for $f\in L^2(E,m)\cap L^4(E,m)$.
Thus, for $f\in L^2(E,m;\mathbb{C})\cap L^4(E,m;\mathbb{C})$, we have
\begin{eqnarray}\label{1.13}
   \P_{\delta_x}|\langle f,X_t\rangle|^2 =  \P_{\delta_x}\langle \Re(f),X_t\rangle^2 +\P_{\delta_x}\langle \Im(f),X_t\rangle^2=  \int_0^tT_{s}[A|T_{t-s}f|^2](x)\,ds+ T_t(|f|^2)(x).
\end{eqnarray}
Let ${\V}{\rm ar}_{\nu}$ be the variance under $\P_{\nu}$.
Then by the branching property, we have
${\V}{\rm ar}_{\nu}\langle f,X_t\rangle=\langle{\V}{\rm ar}_{\delta_\cdot}\langle f,X_t\rangle,\nu\rangle.$
By \eqref{2.27} and \eqref{2.8}, we get, for $t>2t_0$,
\begin{eqnarray}
  {\V}{\rm ar}_{\delta_x}\langle f,X_t\rangle&\le& {\P}_{\delta_x}
  |\langle f,X_t\rangle|^2\le (1+Ke^{Kt})T_t(|f|^2)(x)\nonumber\\
  &\le& (1+Ke^{Kt})e^{-\lambda_1t}b_t(x)^{1/2}\||f|^2\|_2\in L^2(E,m)\cap L^4(E,m). \label{var}
\end{eqnarray}

Recall that $t_0$ is the constant in Assumption 1(c).

\begin{lemma}
Assume that $f\in L^2(E,m;\mathbb{C})\cap L^4(E,m;\mathbb{C})$.
If $\lambda_1>2\Re_{\gamma(f)}$, then for any $(t, x)\in(10t_0,\infty)\times E$ we have,
\begin{equation}\label{1.51}
  \sup_{t>10t_0}t^{-2\tau(f)}e^{2\Re_{\gamma(f)}t}\P_{\delta_x}|\langle f,X_t\rangle|^2\lesssim b_{t_0}(x)^{1/2}.
\end{equation}
\end{lemma}
\textbf{Proof:}
In this proof, we always assume $t>10 t_0$.
For $s\le 2t_0$, we have $T_{t-s}[A|T_{s}f|^2](x)\le Ke^{Ks}T_t(|f|^2)(x)\lesssim T_t(|f|^2)(x)$.
Thus, by \eqref{1.23}, we have for $t>10t_0$,
\begin{equation}\label{1.52}
  \int_0^{2t_0}T_{t-s}[A|T_{s}f|^2](x)\,ds\lesssim T_t(|f|^2)(x)\lesssim e^{-\lambda_1t}b_{t_0}(x)^{1/2}.
\end{equation}
It follows from \eqref{1.23} again that, for $(s,x)\in (8t_0,\infty)\times E$, $|T_s f(x)|\lesssim s^{\tau(f)}e^{-\Re_{\gamma(f)}s}b_{4t_0}(x)^{1/2}$ .
 Thus, for $(t,x)\in(10t_0,\infty)\times E$,
\begin{eqnarray}\label{2.46}
 &&\int_{t-2t_0}^{t}T_{t-s}[A|T_{s}f|^2](x)\,ds\lesssim t^{2\tau(f)}\int_{t-2t_0}^t e^{-2\Re_{\gamma(f)}s}T_{t-s}(b_{4t_0})(x)\,ds\nonumber\\
  &=& t^{2\tau(f)}e^{-2\Re_{\gamma(f)}t}\int_{0}^{2t_0} e^{2\Re_{\gamma(f)}s}T_s(b_{4t_0})(x)\,ds
  \lesssim t^{2\tau(f)}e^{-2\Re_{\gamma(f)}t}\int_{0}^{2t_0}T_s(b_{4t_0})(x)\,ds.
\end{eqnarray}
  We now show that for any $x\in E$, $\int_{0}^{2t_0}T_s(b_{4t_0})(x)\,ds<\infty$.
By \eqref{8.9}, we get
$$
b_{4t_0}(x)\le e^{8Kt_0}a_{4t_0}(x)\le e^{10Kt_0}T_{2t_0}(a_{2t_0})(x).
$$
Thus, by \eqref{2.8}, we have
\begin{equation}\label{1.37}
  \int_{0}^{2t_0}T_s(b_{4t_0})(x)\,ds\le e^{10Kt_0}\int_{0}^{2t_0}T_{s+2t_0}(a_{2t_0})(x)\,ds
  \lesssim \int_{0}^{2t_0}e^{-\lambda_1(s+2t_0)}\,ds b_{t_0}(x)^{1/2}\lesssim b_{t_0}(x)^{1/2}.
\end{equation}
By \eqref{2.46} and \eqref{1.37}, we get
\begin{equation}\label{5.2}
 \int_{t-2t_0}^{t}T_{t-s}[A|T_{s}f|^2](x)\,ds\lesssim t^{2\tau(f)}e^{-2\Re_{\gamma(f)}t}b_{t_0}(x)^{1/2}.
\end{equation}
For $s\in[2t_0, t-2t_0],$ by \eqref{1.23},
we have
$|T_sf(x)|^2\lesssim s^{2\tau(f)}e^{-2\Re_{\gamma(f)}s}b_{t_0}(x).$
By \eqref{2.8}, we get
$T_{t-s}[A(T_{s}f)^2](x)\lesssim s^{2\tau(f)}e^{-2\Re_{\gamma(f)}s}e^{-\lambda_1(t-s)}b_{t_0}(x)^{1/2}$.
So,
for $(t,x)\in(10t_0,\infty)\times E$,
\begin{eqnarray}
 \int_{2t_0}^{t-2t_0}T_{t-s}[A|T_{s}f|^2](x)\,ds&\lesssim& t^{2\tau(f)}e^{-\lambda_1t}\int_0^te^{(\lambda_1-2\Re_{\gamma(f)})s}\,dsb_{t_0}(x)^{1/2}\label{5.1}\\
 &\lesssim& t^{2\tau(f)}e^{-2\Re_{\gamma(f)}t}b_{t_0}(x)^{1/2}.\label{1.54}
\end{eqnarray}
Combining \eqref{1.52}, \eqref{5.2} and \eqref{1.54}, when $\lambda_1>2\Re_{\gamma(f)}$, we get
$$
\int_0^tT_{t-s}[A|T_{s}f|^2](x)\,ds\lesssim t^{2\tau(f)}e^{-2\Re_{\gamma(f)}t}b_{t_0}(x)^{1/2}.
$$
Since $\lambda_1>2\Re_{\gamma(f)}$, by \eqref{2.8}, we have, for $(t,x)\in(10t_0,\infty)\times E$,
$$
T_t(|f|^2)(x)\lesssim e^{-\lambda_1t}b_{t_0}(x)^{1/2}\lesssim t^{2\tau(f)}e^{-2\Re_{\gamma(f)}t}b_{t_0}(x)^{1/2}.
$$
Now \eqref{1.51} follows easily.\hfill$\Box$

\begin{lemma} \label{lem:2.2}
Assume that $f\in L^2(E,m)\cap L^4(E,m)$.
If $\lambda_1<2\Re_{\gamma(f)}$, then for  $(t,x)\in(10t_0,\infty)\times E$,
\begin{equation}\label{small}
   \left|e^{\lambda_1t}{\V}{\rm ar}_{\delta_x}\langle f,X_t\rangle
   - \sigma^2_{f}\phi_1(x)\right|
   \lesssim c_t(b_{t_0}(x)^{1/2}+b_{t_0}(x)),
\end{equation}
where $c_t$ is  independent of $x$ with
$\lim_{t\to\infty}c_t= 0$ and
$\sigma_{f}^2$ is defined in \eqref{e:sigma}.
\end{lemma}
\textbf{Proof:}
First, we consider the case $\gamma(f)<\infty$.
In this proof, we always assume $t>10t_0$ and $f\in L^2(E,m)\cap L^4(E,m)$.
By \eqref{1.23}, we have
\begin{equation}\label{1.6}
  e^{\lambda_1t/2}|\P_{\delta_x}\langle f, X_t\rangle|\lesssim t^{\tau(f)}e^{-(2\Re_{\gamma(f)}-\lambda_1)t/2}b_{t_0}(x)^{1/2}.
\end{equation}
We first show that $\sigma_f^2<\infty$.
For $s\le 2t_0$, by \eqref{Lp}, we have
\begin{equation}\label{1.32}
  \|A|T_s f|^2\|_2\le K\|T_s f\|_4^2\le Ke^{2Ks}\|f\|_4^2.
\end{equation}
For $s>2t_0$, by \eqref{1.23}, $|T_sf(x)|\lesssim e^{-\Re_{\gamma(f)}s}s^{\tau(f)}b_{t_0}(x)^{1/2}$.
Thus, we have
\begin{eqnarray}\label{1.34}
  &&\int_0^\infty e^{\lambda_1 s}\langle A|T_s f|^2,\psi_1\rangle_m \,ds\le K\|\psi_1\|_2\int_0^\infty e^{\lambda_1 s}\||T_s f|^2\|_2\,ds\nonumber\\
  &\lesssim& \int_0^{2t_0}e^{\lambda_1 s}\,ds+\int_{2t_0}^\infty e^{(\lambda_1-2\Re_{\gamma(f)})s}s^{2\tau(f)}\,ds<\infty,
\end{eqnarray}
from which we easily see that $\sigma_f^2<\infty$.
By \eqref{1.13}, we have
\begin{eqnarray}
  &&\left|e^{\lambda_1t}\P_{\delta_x}\langle f,X_t\rangle^2- \sigma_{f}^2 \phi_1(x)\right| \nonumber\\
  &\le& e^{\lambda_1t}\int^{t-2t_0}_0\left|T_{t-s}[A|T_sf|^2](x)-e^{-\lambda_1(t-s)}
 \langle A|T_s f|^2,\psi_1\rangle_m\phi_1(x)\right|\,ds\nonumber\\
  &&+e^{\lambda_1t}\int_{t-2t_0}^tT_{t-s}[A|T_sf|^2](x)\,ds+\int_{t-2t_0}^\infty e^{\lambda_1 s}
  \langle A|T_s f|^2,\psi_1\rangle_m \,ds \phi_1(x)\nonumber\\
   &&+|e^{\lambda_1t}T_t(|f|^2)(x)-\langle |f|^2,\psi_1\rangle_m\phi_1(x)|\nonumber\\
  &=:& V_1(t,x)+V_2(t,x)+V_3(t,x)+V_4(t,x).
\end{eqnarray}

First, we consider $V_1(t,x)$.
By \eqref{2.6}, for $t-s>2t_0$, there exists $a\in (\lambda_1, \Re_2)$ such that
\begin{eqnarray*}
  \left|T_{t-s}[A|T_sf|^2](x)-e^{-\lambda_1(t-s)}
  \langle A|T_s f|^2,\psi_1\rangle_m\phi_1(x)\right|
  &\lesssim&
  e^{-a(t-s)}\|A(T_s f)^2\|_2b_{t_0}(x)^{1/2}.
\end{eqnarray*}
Therefore, by  \eqref{1.23} and \eqref{1.32}, we have
\begin{eqnarray}\label{V1}
  V_1(t,x) &\lesssim&
 e^{\lambda_1t}t^{2\tau(f)}\int_{2t_0}^{t-2t_0}
 e^{-a(t-s)}e^{-2\Re_{\gamma(f)}s}\,ds
 \,b_{t_0}(x)^{1/2}
  +e^{\lambda_1t}\int_{0}^{2t_0}
  e^{-a(t-s)}\,ds
  \,b_{t_0}(x)^{1/2}\nonumber\\
  &\lesssim&e^{-(a-\lambda_1)t}t^{2\tau(f)}\int_{0}^{t} e^{(a-2\Re_{\gamma(f)})s}\,ds\,b_{t_0}(x)^{1/2}+e^{-(a-\lambda_1)t}b_{t_0}(x)^{1/2}\nonumber\\
   &\lesssim& t^{2\tau(f)}\left(e^{(\lambda_1-2\Re_{\gamma(f)})t}+e^{-(a-\lambda_1)t}\right)b_{t_0}(x)^{1/2}.
\end{eqnarray}

Now we deal with $V_2(t,x)$. By \eqref{5.2}, we have
\begin{equation}\label{V2}
  V_2(t,x)\lesssim t^{2\tau(f)}e^{(\lambda_1-2\Re_{\gamma(f)})t}b_{t_0}(x)^{1/2}.
\end{equation}

For $V_3(t,x)$, by \eqref{1.34}, we get
$\int_{t-2t_0}^\infty e^{\lambda_1 s}\langle A|T_s f|^2,\psi_1\rangle_m \,ds\to 0,$ as $t\to\infty$.
By \eqref{phi}, we have $\phi_1(x)\lesssim b_{t_0}(x)^{1/2}$.

Finally, we consider $V_4(t,x)$.
By \eqref{2.6}, we have
\begin{equation}\label{V4}
  V_4(t,x)\lesssim e^{-(a-\lambda_1)t}b_{t_0}(x)^{1/2}.
\end{equation}
Thus, by \eqref{V1}--\eqref{V4}, we have, for $(t,x)\in(10t_0,\infty)\times E$,
\begin{equation}\label{2.22}
  \left|e^{\lambda_1t}\P_{\delta_x}\langle f,X_t\rangle^2 - \sigma_{f}^2 \phi_1(x)\right|
 \lesssim c_tb_{t_0}(x)^{1/2},
\end{equation}
with $\lim_{t\to\infty}c_t=0$.
Now \eqref{small} follows immediately from \eqref{1.6} and \eqref{2.22}.

Now, we consider the case $\gamma(f)=\infty$. The proof is similar to that of the
case $\gamma(f)<\infty$, the only difference being that we now use \eqref{1.23'} instead of \eqref{1.23}.
\hfill$\Box$

\begin{lemma}\label{critical}
Assume that $f,h\in L^2(E,m)\cap L^4(E,m)$.
If $\lambda_1=2\Re_{\gamma(f)}=2\Re_{\gamma(h)}$, then for  $(t,x)\in(10t_0,\infty)\times E$,
\begin{equation}\label{7.49}
\left|t^{-(1+\tau(f)+\tau(h))}e^{\lambda_1 t}\mathbb{C}{\rm ov}_{\delta_x}(\langle f,X_t\rangle,\langle h,X_t\rangle)-\rho(f,h)\phi_1(x)\right|
\lesssim t^{-1}\left(b_{t_0}(x)^{1/2}+b_{t_0}(x)\right),
\end{equation}
where $\mathbb{C}{\rm ov}_{\delta_x}$ is the covariance under $\P_{\delta_x}$ and $\rho(f,h)$ is defined by
\eqref{rho2} with $f$ and $h$ in place of $h_1$ and $h_2$ respectively.
In particular, we have,  for  $(t,x)\in(10t_0,\infty)\times E$,
\begin{equation}\label{1.49}
\left|t^{-(1+2\tau(f))}e^{\lambda_1 t}{\V}{\rm ar}_{\delta_x}\langle f,X_t\rangle-\rho_{f}^2\phi_1(x)\right|
\lesssim t^{-1}\left(b_{t_0}(x)^{1/2}+b_{t_0}(x)\right),
\end{equation}
where $\rho^2_{f}$ is defined by \eqref{e:rho}.
Moreover, we have,  for  $(t,x)\in(10t_0,\infty)\times E$,
\begin{equation}\label{3.33}
  t^{-(1+2\tau(f))}e^{\lambda_1 t}{\V}{\rm ar}_{\delta_x}\langle f,X_t\rangle\lesssim \left(b_{t_0}(x)^{1/2}+b_{t_0}(x)\right).
\end{equation}

\end{lemma}
\textbf{Proof:}
 In this proof we always assume $t>10t_0$ and $f,h\in L^2(E,m)\cap L^4(E,m)$.
By \eqref{1.13}, we have
\begin{eqnarray}\label{7.1}
  &&\mathbb{C}{\rm ov}_{\delta_x}(\langle f,X_t\rangle, \langle h,X_t\rangle) \nonumber\\
  &=& \frac{1}{4}\left(\mathbb{V}{\rm ar}_{\delta_x}\langle (f+h),X_t\rangle- \mathbb{V}{\rm ar}_{\delta_x}\langle (f-h),X_t\rangle\right)\nonumber\\
   &=& \int_0^t T_{t-s}\left[A(T_sf)(T_sh)\right](x)\,ds+T_t(fh)(x)-T_t(f)(x)T_t(h)(x).
\end{eqnarray}
Let
$$C_f(s,x):=\sum_{j:\lambda_1=2\Re_{j}}\left(e^{-i\Im_j s}(\Phi_j(x))^TF_{f,j}\right),\quad C_h(s,x):=\sum_{j:\lambda_1=2\Re_{j}}\left(e^{-i\Im_j s}(\Phi_j(x))^TF_{h,j}\right).$$
Define
\begin{eqnarray*}
  V_5(t,x)&:=&e^{\lambda_1t}\int^{t-2t_0}_{2t_0}T_{t-s}[A(T_{s}f)(T_sh))](x)\,ds,\\
  V_6(t,x)&:=& e^{\lambda_1t}\int^{t-2t_0}_{2t_0}s^{\tau(f)+\tau(h)}e^{-\lambda_1s}T_{t-s}[AC_f(s,\cdot)C_h(s,\cdot)](x)\,ds, \\
  V_7(t,x)&:=& \int^{t-2t_0}_{2t_0}s^{\tau(f)+\tau(h)}\langle AC_f(s,\cdot)C_h(s,\cdot),\psi_1\rangle_m\,ds\phi_1(x) \\
  \mbox{and}\qquad\qquad&&\\
 V_8(t,x)&:=&\int_{2t_0}^{t-2t_0} s^{\tau(f)+\tau(h)}\langle AF_{f,h},\psi_1\rangle_m\,ds\phi_1(x),
\end{eqnarray*}
where $F_{f, h}$ is defined in \eqref{e:Ffg} with $f$ and $h$ in place of $h_1$ and $h_2$ respectively.
It is easy to see from the definition of $\rho(f,h)$ that
$$\rho(f,h)=t^{-(1+\tau(f)+\tau(h))}\int_{0}^{t} s^{\tau(f)+\tau(h)}\langle AF_{f,h},\psi_1\rangle_m\,ds.$$
Thus we have
\begin{eqnarray}\label{1.45}
  &&\left|e^{\lambda_1t}\int_0^tT_{t-s}[A(T_{s}f)(T_sh)](x)\,ds- t^{1+\tau(f)+\tau(h)}\rho(f,h)\phi_1(x)\right| \nonumber\\
  &\le&e^{\lambda_1t}\left(\int_0^{2t_0}+\int_{t-2t_0}^t\right)T_{t-s}[A|T_{s}f||T_sh|](x)\,ds+|V_5(t,x)-V_6(t,x)|+|V_6(t,x)-V_7(t,x)|\nonumber\\
  &&+|V_7(t,x)-V_8(t,x)|+\left(\int_0^{2t_0}+\int_{t-2t_0}^t\right) s^{\tau(f)+\tau(h)}\,ds\langle AF_{f,h},\psi_1\rangle_m\phi_1(x).
\end{eqnarray}
By \eqref{2.8}, for $s\le t-2t_0$, we have
$$T_{t-s}[A|T_sf||T_sh|](x)\lesssim e^{-\lambda_1(t-s)}\left\|A|T_sf||T_sh|\right\|_2(b_{t_0}(x))^{1/2}.$$
By \eqref{Lp}, it is easy to see that
$\|A|T_sf||T_sh|\|_2\le K\|T_sf\|_4\|T_sh\|_4\le Ke^{2Ks}\|f\|_4\|h\|_4.$
Thus,
\begin{equation}\label{7.3}
  e^{\lambda_1t}\int_0^{2t_0}T_{t-s}[A|T_{s}f||T_sh|](x)\,ds\lesssim \int_0^{2t_0} e^{\lambda_1s}\,ds(b_{t_0}(x))^{1/2}\lesssim (b_{t_0}(x))^{1/2}.
\end{equation}
For $s>t-2t_0$, using arguments similar to those leading to \eqref{5.2}, we get
\begin{eqnarray}
  &&e^{\lambda_1t}\int_{t-2t_0}^t T_{t-s}\left[A|T_sf||T_sh|\right](x)\,ds
   \lesssim t^{\tau(f)+\tau(h)}e^{\lambda_1t}e^{-(\Re_\gamma(h)+\Re_{\gamma(f)})t}(b_{t_0}(x))^{1/2}\label{7.5}\\
   &=& t^{\tau(f)+\tau(h)}(b_{t_0}(x))^{1/2}.\label{7.6}
\end{eqnarray}
By \eqref{phi}, it is easy to see that
\begin{equation}\label{1.42}
  \left(\int_0^{2t_0}+\int_{t-2t_0}^t\right) s^{\tau(f)+\tau(h)}\,ds\langle AF_{f,h},\psi_1\rangle_m\phi_1(x)\lesssim t^{\tau(f)+\tau(h)}b_{t_0}(x)^{1/2}.
\end{equation}

Next we consider $|V_5(t,x)-V_6(t,x)|$.
By \eqref{1.31}, we have, for $(s,x)\in (2t_0,\infty)\times E$,
\begin{eqnarray*}
  |T_{s}f(x)-s^{\tau(f)}e^{-\lambda_1s/2}C_f(s,x)|&\lesssim& s^{\tau(f)-1}e^{-\lambda_1s/2}b_{t_0}(x)^{1/2} .
\end{eqnarray*}
The same is also true for $h$.
Thus by \eqref{1.23} and \eqref{phi}, we get, for $(s,x)\in(2t_0,\infty)\times E$,
\begin{eqnarray}\label{1.36}
  &&\left||T_{s}f(x)T_sh(x)|-s^{\tau(f)+\tau(h)}e^{-\lambda_1s}C_f(s,x)C_h(s,x)\right|\nonumber\\
  &\lesssim&\left|T_{s}f(x)-s^{\tau(f)}e^{-\lambda_1s/2}C_f(s,x)\right|\left|T_{s}h(x)-s^{\tau(h)}e^{-\lambda_1s/2}C_h(s,x)\right|\nonumber\\
  &&+s^{\tau(h)}e^{-\lambda_1s/2}\left|T_{s}f(x)-s^{\tau(f)}e^{-\lambda_1s/2}C_f(s,x)\right||C_h(s,x)|\nonumber\\
  &&+s^{\tau(f)}e^{-\lambda_1s/2}\left|T_{s}h(x)-s^{\tau(h)}e^{-\lambda_1s/2}C_h(s,x)\right||C_f(s,x)|\nonumber\\
  &\lesssim& s^{\tau(f)+\tau(h)-1}e^{-\lambda_1s}b_{t_0}(x).
\end{eqnarray}
Therefore, by \eqref{2.8},  we have, for $(t,x)\in(10t_0,\infty)\times E$,
\begin{eqnarray}\label{1.44}
&&|V_5(t,x)-V_6(t,x)|
   \lesssim   \int^{t-2t_0}_{2t_0}s^{\tau(f)+\tau(h)-1}e^{\lambda_1(t-s)}T_{t-s}(b_{t_0})(x)\,ds\nonumber\\
   &\lesssim &\int^{t-2t_0}_{2t_0}s^{\tau(f)+\tau(h)-1}\,ds b_{t_0}(x)^{1/2}
   \lesssim t^{\tau(f)+\tau(h)}b_{t_0}(x)^{1/2}.
\end{eqnarray}

For $|V_6(t,x)-V_7(t,x)|$, by \eqref{2.6}, there exists $\lambda_1<a<\Re_2$, such that, for $t-s>2t_0$,
\begin{eqnarray*}
  &&\left|e^{\lambda_1(t-s)}T_{t-s}\left[AC_f(s,\cdot)C_h(s,\cdot)\right](x)-\langle AC_f(s,\cdot)C_h(s,\cdot),\psi_1\rangle_m\phi_1(x)\right| \\
 &\lesssim& e^{-(a-\lambda_1)(t-s)}\|C_f(s,\cdot)C_h(s,\cdot)\|_2b_{t_0}(x)^{1/2}.
\end{eqnarray*}
By \eqref{phi}, we get, for $s>2t_0$,
$|C_f(s,x)C_h(s,x)|\lesssim b_{t_0}(x).$
Thus, we get
\begin{eqnarray}\label{1.38}
  &&|V_6(t,x)-V_7(t,x)|\lesssim \int_{2t_0}^{t-2t_0}s^{\tau(f)+\tau(h)}
  e^{-(a-\lambda_1)(t-s)}\,dsb_{t_0}(x)^{1/2}\nonumber\\
  &\lesssim&  t^{\tau(f)+\tau(h)}\int_{2t_0}^{t-2t_0}e^{-(a-\lambda_1)(t-s)}\,dsb_{t_0}(x)^{1/2}\lesssim t^{\tau(f)+\tau(h)}b_{t_0}(x)^{1/2}.
\end{eqnarray}

Now we deal with  $|V_7(t,x)-V_8(t,x)|$.
We can check that $C_h(s,x)$ is real.
In fact, for each $j$ with $\lambda_1=2\Re_j$, we also have $\lambda_1=2\Re_{j'}$
and
$e^{-i\Im_{j'} s}(\Phi_{j'}(x))^TF_{h,j'}=\overline{e^{-i\Im_j s}(\Phi_j(x))^TF_{h,j}}$.
Thus, we have $C_h(s,x)=\overline{C_h(s,x)}=\sum_{j:\lambda_1=2\Re_{j}}\left(e^{i\Im_j s}\overline{(\Phi_j(x))^TF_{h,j}}\right)$.
Therefore,
\begin{eqnarray*}
  C_f(s,x)C_h(s,x)&=&\sum_{j:\lambda_1=2\Re_j}\Phi_j(x))^TF_{f,j}\overline{(\Phi_j(x))^TF_{h,j}}\\
&&+\sum_{\gamma(f)\le j\ne l\le \zeta(f)}\left(e^{-i(\Im_j-\Im_l) s}(\Phi_j(x))^TF_{f,j}(\overline{\Phi_l(x))^TF_{h,l}}\right).
\end{eqnarray*}
When $j\ne l$, since $\lambda_j\ne\lambda_l$ and $\Re_j=\Re_l$, we have $\Im_j\ne\Im_l$.

We claim that for any non-zero $\theta\in\mathbb{R}$ and $n\ge 0$, we have for $t>2t_0$,
\begin{equation}\label{2.26}
 \left|\int^{t-2t_0}_{2t_0}s^{n}e^{i\theta s}\,ds\right|\lesssim t^{n}.
\end{equation}
Then, using \eqref{e:Ffg}, we get
\begin{eqnarray}\label{1.41}
  &&|V_7(t,x)-V_8(t,x)|\nonumber\\
  &\lesssim&  \sum_{\gamma(f)\le j\ne l\le\zeta(f)} \left|\int_{2t_0}^{t-2t_0}s^{\tau(f)+\tau(h)}e^{-i(\Im_j-\Im_l) s}\,ds\right|\left|\langle(\Phi_j(x))^TF_{f,j}(\overline{\Phi_l(x)})^T\overline{F_{h,l}},\psi_1\rangle_m\right|\phi_1(x)\nonumber\\
  &\lesssim& t^{\tau(f)+\tau(h)}b_{t_0}(x)^{1/2}.
\end{eqnarray}
Now we prove \eqref{2.26}.  Using integration by parts, for $n\ge 1$, we get
\begin{eqnarray*}
  &&\left|\int_{2t_0}^{t-2t_0}s^{n}e^{i\theta s}\,ds \right|
  =\left|  \frac{s^{n}e^{i\theta s}|_{2t_0}^{t-2t_0}-\int_{2t_0}^{t-2t_0}ns^{n-1}e^{i\theta s}\,ds}{i\theta}\right|
   \lesssim t^{n}+\int_{2t_0}^{t-2t_0}s^{n-1}\,ds\lesssim t^{n}.
\end{eqnarray*}
For $n=0$, we have
$$\left|\int_{2t_0}^{t-2t_0}e^{i\theta s}\,ds \right|=\left|\frac{e^{i\theta(t-2t_0)}-e^{i2\theta t_0}}{i\theta}\right|\le 2/|\theta|.$$
Therefore, \eqref{2.26} follows immediately.

Combining \eqref{7.3}, \eqref{7.6},  \eqref{1.42}, \eqref{1.44}, \eqref{1.38} and \eqref{1.41}, we get $(t,x)\in(10t_0,\infty)\times E$,
\begin{equation}\label{1.46}
  \left|e^{\lambda_1t}\int_0^tT_{t-s}[A(T_{s}f)(T_sh)](x)\,ds- t^{1+\tau(f)+\tau(h)}\rho(f,h)\phi_1(x)\right|  \lesssim t^{\tau(f)+\tau(h)}b_{t_0}(x)^{1/2}.
\end{equation}
By \eqref{2.8}, we have, for $(t,x)\in(10t_0,\infty)\times E$,
$$e^{\lambda_1t}T_t(|fh|)(x)\lesssim b_{t_0}(x)^{1/2}.$$
And by \eqref{2.6} and $\lambda_1=2\Re_{\gamma(f)}$,
$$e^{\lambda_1t}|T_tf(x)||T_th(x)|\lesssim t^{\tau(f)+\tau(h)}b_{t_0}(x).$$
Now \eqref{7.49} follows immediately.
\hfill$\Box$

\begin{lemma}\label{small and critical}
Assume that $f\in L^2(E,m)\cap L^4(E,m)$ with $\lambda_1<2\Re_{\gamma(f)}$
and $h\in L^2(E,m)\cap L^4(E,m)$ with $\lambda_1=2\Re_{\gamma(h)}$. Then, for any $(t, x)\in (10t_0, \infty)\times E$,
\begin{equation}\label{cov:sc}
 e^{\lambda_1t}\mathbb{C}{\rm ov}_{\delta_x}(\langle f,X_t\rangle, \langle h,X_t\rangle)\lesssim ((b_{t_0}(x))^{1/2}+b_{t_0}(x)).
\end{equation}
\end{lemma}
\textbf{Proof:}
In this proof, we always assume that $t>10t_0$, $f\in L^2(E,m)\cap L^4(E,m)$ with $\lambda_1<2\Re_{\gamma(f)}$
and $h\in L^2(E,m)\cap L^4(E,m)$ with $\lambda_1=2\Re_{\gamma(h)}$.
First, we assume $\gamma(f)<\infty$.
By \eqref{7.1}, we have
\begin{eqnarray*}
  &&\mathbb{C}{\rm ov}_{\delta_x}(\langle f,X_t\rangle, \langle h,X_t\rangle)
  =\int_0^t T_{t-s}\left[A(T_sf)(T_sh)\right](x)\,ds+T_t(fh)(x)-T_t(f)(x)T_t(h)(x).
\end{eqnarray*}
By \eqref{7.3} and \eqref{7.5}, we have, for $(t, x)\in (10t_0, \infty)\times E$,
\begin{eqnarray*}
  &&e^{\lambda_1t}\left(\int_{0}^{2t_0}+\int_{t-2t_0}^t\right) T_{t-s}\left[A|T_sf||T_sh|\right](x)\,ds\\
   &\lesssim& b_{t_0}(x)^{1/2}+t^{\tau(f)+\tau(h)}e^{(\lambda_1/2-\Re_{\gamma(f)})t}(b_{t_0}(x))^{1/2}\lesssim (b_{t_0}(x))^{1/2}.
\end{eqnarray*}
By \eqref{1.23}, we have
\begin{eqnarray}\label{1.57}
 &&e^{\lambda_1t}\int_{2t_0}^{t-2t_0} T_{t-s}\left[A|T_sf||T_sh|\right](x)\,ds
  \lesssim e^{\lambda_1t}\int_{2t_0}^{t-2t_0} s^{\tau(f)+\tau(h)}e^{-(\lambda_1/2+\Re_{\gamma(f)})s}T_{t-s}(b_{t_0})(x)\,ds\nonumber\\
  &\lesssim & \left(\int_{2t_0}^{t-2t_0}s^{\tau(f)+\tau(h)}e^{(\lambda_1/2-\Re_{\gamma(f)})s}\,ds\right)\,b_{t_0}(x)^{1/2}
   \lesssim b_{t_0}(x)^{1/2}.
\end{eqnarray}
Thus, we have
\begin{equation}\label{1.61}
  e^{\lambda_1t}\left|\int_{0}^t T_{t-s}\left[A(T_sf)(T_sh)\right](x)\,ds\right|\lesssim (b_{t_0}(x))^{1/2}.
\end{equation}
By \eqref{2.8}, we get
$$e^{\lambda_1t}|T_t(fh)(x)|\le e^{\lambda_1t}T_t(|fh|)(x)\lesssim b_{t_0}(x)^{1/2}.$$
By \eqref{1.23}, for $(t,x)\in(10t_0,\infty)\times E$, we have
$$e^{\lambda_1t}|T_tf(x)T_th(x)|\lesssim t^{\tau(f)+\tau(h)}e^{(\lambda_1/2-\Re_{\gamma(f)})t}b_{t_0}(x)\lesssim b_{t_0}(x).$$
Now \eqref{cov:sc} follows immediately.

Repeating the proof above by using \eqref{1.23'} instead of \eqref{1.23},
 we get \eqref{cov:sc} also holds when $\gamma(f)=\infty.$
\hfill$\Box$

\section{Proofs of Main Results}

In this section, we will prove the main results of this paper.
When referring to individuals in $X$ we will use the classical Ulam-Harris notation
so that every individual in
$X$ has a unique label, see \cite{HH}.
For each individual
$u\in{\cal T}$ we shall write
$b_u$ and $d_u$ for its birth and death times respectively and
$\{z_u(r): r\in [b_u,d_u]\}$ for its spatial trajectory.
Define
$${\cal L}_t=\{ u\in {\cal T}, b_u\le t<d_u\},\quad t\ge 0.$$
Thus, $X_{s+t}$ has the following decomposition:
\begin{equation}\label{3.22}
  X_{s+t}=\sum_{u\in \cL_t}X^{u,t}_s,
\end{equation}
where given $\mathcal{F}_t$, $X^{u,t}_s$, $u\in \cL_t$, are independent and $X^{u,t}_s$ has the same law as $X_s$ under $\P_{\delta_{z_u(t)}}$.

\subsection{A basic law of large numbers}
Recall that
$$
H_{t}^{(k)}:=e^{\lambda_kt}(\langle \phi_1^{(k)}, X_t\rangle, \cdots,\langle\phi_{n_k}^{(k)},X_t\rangle) (D_k(t))^{-1}.
$$
\begin{lemma}\label{thrm1}
Assume that $b$ is an $n_k$-dimensional vector. If $\lambda_1>2\Re_k$, then, for any $\nu\in {\cal M}_a(E)$, $H_t^{(k)}b$ is a martingale under $\P_{\nu}$.
Moreover, the limit
\begin{equation}\label{1.48}
H_\infty^{(k)}:=\lim_{t\to\infty}H_t^{(k)}
\end{equation}
exists $\P_{\nu}$-a.s. and in $L^2(\P_{\nu})$.
\end{lemma}

\textbf{Proof:}
By the branching property, it suffices to prove the lemma for $\nu=\delta_x$ with $x\in E$.
By \eqref{T-Jordan}, we have
$$\P_{\delta_x}H_t^{(k)}b=e^{\lambda_kt}T_t((\Phi_k)^T)(x)(D_k(t))^{-1}b=(\Phi_k(x))^Tb.$$
Thus, by the Markov property, we get that $H_t^{(k)}b$ is a martingale under $\P_{\delta_x}$.
We claim that, for $(t,x)\in(2t_0,\infty)\times E$,
\begin{equation}\label{1.47}
  \P_{\delta_x}|H_t^{(k)}b|^2\lesssim |b|_\infty^2b_{t_0}(x)^{1/2},
\end{equation}
from which \eqref{1.48} follows immediately.

Now we prove the claim.
Let $f_t(x)=e^{\lambda_kt}b^T(D_k(t)^{-1})^T\Phi_k(x)$.
Then $H_t^{(k)}b=\langle f_t,X_t\rangle$, and
by \eqref{T-Jordan}, for $s<t$, we have
$$T_s(f_t)(x)=e^{\lambda_k(t-s)}b^T(D_k(t-s)^{-1})^T\Phi_k=f_{t-s}(x).$$
By \eqref{1.13}, we have
\begin{equation}\label{1.50}
  \P_{\delta_x}|H_t^{(k)}b|^2=\P_{\delta_x}|\langle f_t,X_t\rangle|^2
  =\int_0^tT_{s}[A|f_{s}|^2](x)\,ds+T_t(|f_t|^2)(x).
\end{equation}
Since each component of $D_k(s)^{-1}=D_k(-s)$ is a polynomial of $s$ with degree
no larger than $\nu_k$, we get  $|D_k(s)^{-1}|_\infty\lesssim (1+s^{\nu_k})$.
Thus, for all $s>0$, we have
\begin{equation}\label{7.8}
  |f_s|\lesssim e^{\Re_ks}|b|_\infty|D_k(s)|_\infty|\Phi_k(x)|_\infty\lesssim |b|_\infty(1+s^{\nu_k})e^{\Re_ks}b_{4t_0}(x)^{1/2}.
\end{equation}
By \eqref{2.8}, we have, for $(s,x)\in(2t_0,\infty)\times E$,
\begin{equation}\label{1.56}
  T_s(|f_s|^2)(x)\lesssim e^{-\lambda_1s}\||f_s|^2\|_2b_{t_0}(x)^{1/2}\lesssim |b|_\infty^2(1+s^{2\nu_k})e^{-(\lambda_1-2\Re_k)s}b_{t_0}(x)^{1/2}.
\end{equation}
Thus, we have
\begin{equation}\label{7.7}
   \int_{2t_0}^tT_{s}[A|f_{s}|^2](x)\,ds\lesssim|b|_\infty^2 b_{t_0}(x)^{1/2}.
\end{equation}
By \eqref{7.8} and \eqref{1.37}, we get
\begin{equation}\label{1.59}
  \int_0^{2t_0}T_{s}[A|f_{s}|^2](x)\,ds\lesssim |b|_\infty^2 \int_0^{2t_0}T_{s}b_{4t_0}(x)\,ds\lesssim |b|_\infty^2 b_{t_0}(x)^{1/2}.
\end{equation}
Thus, by \eqref{7.7} and \eqref{1.59}, we have
\begin{equation}\label{1.55}
  \int_{0}^tT_{s}[A|f_{s}|^2](x)\,ds\lesssim|b|_\infty^2 b_{t_0}(x)^{1/2}.
\end{equation}
Since $\lambda_1>2\Re_k$,  we have $\sup_{s>2t_0}(1+s^{2\nu_k})e^{-(\lambda_1-2\Re_k)s}<\infty$.
Thus, by \eqref{1.56}, we get $$T_t(|f_t|^2)(x)\lesssim |b|_\infty^2b_{t_0}(x)^{1/2},$$
from which \eqref{1.47} follows immediately.
\hfill$\Box$

Now, we present the proof of Theorem \ref{thrm2}.

\textbf{Proof of Theorem \ref{thrm2}:}
By the branching property, it suffices to prove the theorem for $\nu=\delta_x$ with $x\in E$.
Put
$$f^*(x):=\sum_{j=\gamma(f)}^{\zeta(f)}\Phi_j(x)^T\langle f,\Psi_j\rangle,\quad \widetilde{f}(x):=f(x)-f^*(x)$$
and
$f_t(x):=\sum_{j=\gamma(f)}^{\zeta(f)}\Phi_j(x)^TD_j(t)^{-1}F_{f,j}.$
Then
$$t^{-\tau(f)}f^*(x)-f_t(x)=\sum_{j=\gamma(f)}^{\zeta(f)}\Phi_j(x)^TD_j(t)^{-1}\left(t^{-\tau(f)}D_j(t)\langle f,\Psi_j\rangle-F_{f,j}\right).$$
By \eqref{1.47} and \eqref{1.22}, we have, for $(t,x)\in(2t_0,\infty)\times E$,
\begin{eqnarray}\label{2.17}
  &&\P_{\delta_x}\left|t^{-\tau(f)}e^{\Re_{\gamma(f)}t}\langle f^*,X_t\rangle-e^{\Re_{\gamma(f)}t}\langle f_t ,X_t\rangle\right|^2\nonumber\\
  &\lesssim& \sum_{j=\gamma(f)}^{\zeta(f)} |t^{-\tau(f)}D_j(t)\langle f,\Psi_j\rangle-F_{f,j}| _\infty^2 b_{t_0}(x)^{1/2}
  \lesssim t^{-2}b_{t_0}(x)^{1/2}.
\end{eqnarray}
By the definition of $H_t^{(j)}$ and \eqref{1.48}, we have, as $t\to\infty$,
\begin{equation}\label{2.16}
  e^{\Re_{\gamma(f)}t}\langle f_t ,X_t\rangle-\sum_{j=\gamma(f)}^{\zeta(f)}\left(e^{-i\Im_j t}H_\infty^{(j)}F_{f,j}\right)=\sum_{j=\gamma(f)}^{\zeta(f)}\left(e^{-i\Im_j t}(H_t^{(j)}-H_\infty^{(j)})F_{f,j}\right)\to0,
\end{equation}
in $L^2(\P_{\delta_x})$.
Thus, by \eqref{2.17} and \eqref{2.16}, we obtain that, as $t\to\infty$,
\begin{equation}\label{2.18}
  t^{-\tau(f)}e^{\Re_{\gamma(f)}t}\langle f^*,X_t\rangle-\sum_{j=\gamma(f)}^{\zeta(f)}\left(e^{-i\Im_j t}H_\infty^{(j)}F_{f,j}\right)\to 0,
  \quad\mbox{in }L^2(\P_{\delta_x}).
\end{equation}

Now, to complete the proof, we only need to show that, as $t\to\infty$,
\begin{equation}\label{2.21}
  t^{-2\tau(f)}e^{2\Re_{\gamma(f)}t}\P_{\delta_x}|\langle\widetilde{ f},X_t\rangle|^2\to 0.
\end{equation}

(1) If $\lambda_1>2\Re_{\gamma(\widetilde{f})}$, then by \eqref{1.51}, we get, for $(t,x)\in(2t_0,\infty)\times E$, as $t\to\infty$,
$$t^{-2\tau(f)}e^{2\Re_{\gamma(f)}t}\P_{\delta_x}|\langle\widetilde{ f},X_t\rangle|^2\lesssim t^{-2\tau(f)}t^{2\tau(\widetilde{f})}e^{2(\Re_{\gamma(f)}-\Re_{\gamma(\widetilde{f})})t}b_{t_0}(x)^{1/2}\to 0.$$

(2) If $\lambda_1=2\Re_{\gamma(\widetilde{f})}$, then by \eqref{1.49}, we get, as $t\to\infty$,
$$t^{-2\tau(f)}e^{2\Re_{\gamma(f)}t}\P_{\delta_x}|\langle\widetilde{ f},X_t\rangle|^2
=t^{-2\tau(f)}t^{(1+2\tau(\widetilde{f}))t}e^{(2\Re_{\gamma(f)}-\lambda_1)t}t^{-(1+2\tau(wide{f}))t}e^{\lambda_1t}
\P_{\delta_x}|\langle\widetilde{ f},X_t\rangle|^2\to 0. $$

(3) If $\lambda_1<2\Re_{\gamma(\widetilde{f})}$, then by \eqref{small}, we get, as $t\to\infty$,
$$t^{-2\tau(f)}e^{2\Re_{\gamma(f)}t}\P_{\delta_x}|\langle\widetilde{ f},X_t\rangle|^2=t^{-2\tau(f)}e^{(2\Re_{\gamma(f)}-\lambda_1)t}e^{\lambda_1t}\P_{\delta_x}|\langle\widetilde{ f},X_t\rangle|^2\to 0. $$
Combining the three cases above, we get \eqref{2.21}.
The proof is now complete.\hfill$\Box$

\subsection{Proof of the main theorem}\label{s:3}
First, we recall a metric on the space of distributions on $\mathbb{R}^d$.
For $f:\mathbb{R}^d\to\mathbb{R}$, define
$$
\|f\|_{BL}:=\|f\|_{\infty}+ \sup_{x\ne y}\frac{|f(x)-f(y)|}{|x-y|}.
$$
For any distributions $\nu_1$ and $\nu_2$ on $\mathbb{R}^d$, define
\begin{equation*}
  \beta(\nu_1,\nu_2):=\sup\left\{\left|\int f\,d\nu_1-\int f\,d\nu_2\right|~:~\|f\|_{BL}\leq1\right\}.
\end{equation*}
Then $\beta$ is a metric. It follows from \cite[Theorem 11.3.3]{Dudley} that the topology generated by this metric
is equivalent to the weak convergence topology.
 From the definition, we can easily see that, if $\nu_1$ and $\nu_2$ are the distributions of
 two $\mathbb{R}^d$-valued random variables $X$ and $Y$ respectively, then
\begin{equation}\label{5.20}
  \beta(\nu_1,\nu_2)\leq \E\|X-Y\|\leq\sqrt{ \E\|X-Y\|^2}.
\end{equation}

\begin{lemma}\label{lem:small}
If $f\in \C_s$,
then $\sigma_{f}^2\in (0, \infty)$ and,
for any nonzero $\nu\in {\cal M}_a(E)$, it holds under $\P_{\nu}$ that
$$
  \left(e^{\lambda_1 t}\langle \phi_1, X_t\rangle, ~e^{\lambda_1t/2}\langle f , X_t\rangle \right)\stackrel{d}{\rightarrow}\left(W_\infty,~G_1(f)\sqrt{W_\infty}\right), \quad t\to\infty,
$$
where $G_1(f)\sim \mathcal{N}(0,\sigma_{f}^2)$. Moreover, $W_\infty$ and $G_1(f)$ are independent.
\end{lemma}

\textbf{Proof:} The proof is similar that of \cite[Theorem 1.8 ]{RSZ2}.
We define an ${\mathbb R}^2$-valued random variable $U_1(t)$ by
\begin{equation}\label{6.4}
   U_1(t):=\left(e^{\lambda_1 t}\langle \phi_1,X_t\rangle, e^{\lambda_1t/2}\langle f, X_t\rangle\right).
\end{equation}
To prove this lemma, it suffices to show that,
for any $x\in E$, under $\P_{\delta_x}$,
\begin{equation}\label{6.5}
   U_1(t)\stackrel{d}{\to}\left(W_\infty, \sqrt{W_\infty}G_1(f)\right),
\end{equation}
where $G_1(f)\sim\mathcal{N}(0,\sigma_f^2)$ is independent of $W_\infty$.
In fact, if
$\nu=\sum_{j=1}^n\delta_{x_j}, n=1, 2, \dots,
\{x_j; j=1,\cdots, n\}\subset E$, then
$$
X_t=\sum_{j=1}^nX_t^j,
$$
where $X^j_t$ is a branching Markov process starting from
$\delta_{x_j}, j=1, \dots, n$, and $X^j, j=1,\cdots, n$, are independent.
If \eqref{6.5} is valid, we put
$W^j_\infty:=\lim_{t\to\infty}e^{\lambda_1t}\langle \phi_1, X_t^j\rangle$.
Then we easily get that, under $\P_{\nu}$, $W_\infty=\sum_{j=1}^n W^j_\infty$.
Since $\lambda_1<2\Re_{\gamma(f)}$,
\begin{eqnarray*}
  &&\P_{\nu}\exp\left\{i\theta_1e^{\lambda_1 t}\langle \phi_1, X_t\rangle+i\theta_2e^{(\lambda_1/2) t}\langle f, X_t\rangle\right) \\
  &=&\prod_{j=1}^n  \P_{\delta_{x_j}}\exp\left\{i\theta_1e^{\lambda_1 t}\langle \phi_1, X_t^j\rangle+i\theta_2e^{(\lambda_1/2) t}\langle f, X_t^j\rangle\right) \\
  &\to& \prod_{j=1}^n\P_{\delta_{x_j}}\exp\left\{i\theta_1W^j_\infty-\frac{1}{2}\theta_2^2\sigma_f^2W^j_\infty\right) \\
  &=& \P_{\nu}\exp\left\{i\theta_1W_\infty-\frac{1}{2}\theta_2^2\sigma_f^2W_\infty\right),
\end{eqnarray*}
which implies that \eqref{6.5} is valid for $\P_\nu$.

Now we show that \eqref{6.5} is valid.
In the remainder of this proof, we assume $s, t >10t_0$ and write
\begin{equation*}
  U_1(s+t)=\left(e^{\lambda_1 (s+t)}\langle \phi_1,X_{t+s}\rangle, e^{(\lambda_1/2)(s+t)}\langle f,X_{s+t}\rangle\right).
\end{equation*}
Recall the decomposition of $X_{s+t}$ in \eqref{3.22}.
Define
\begin{equation}\label{e:new}
Y_1^{u,t}(s):=e^{\lambda_1 s/2}\langle f,X_{s}^{u,t}\rangle\quad\mbox{and}\quad y_1^{u,t}(s):=\P_{\delta_x}(Y^{u,t}_1(s)|\mathcal{F}_t).
\end{equation}
Given $\mathcal{F}_t$, $Y_1^{u,t}(s)$ has the same law as $Y_1(s):=e^{\lambda_1 s/2}\langle f,X_{s}\rangle$ under $\P_{\delta_{z_u(t)}}$.
Then we have
\begin{eqnarray}\label{6.14}
  &&e^{(\lambda_1/2) (s+t)}\langle f,X_{s+t}\rangle= e^{(\lambda_1/2) t}\sum_{u\in\cL_t}Y_1^{u,t}(s)\nonumber\\
  &=& e^{(\lambda_1/2) t}\sum_{u\in\cL_t}(Y_1^{u,t}(s)-y^{u,t}_s)
      + e^{(\lambda_1/2) (t+s)}\P_{\delta_x}(\langle f,X_{s+t}\rangle|\mathcal{F}_t)\nonumber\\
  &=:& J_1(s,t)+J_2(s,t).
\end{eqnarray}

We first consider $J_{2}(s,t)$.
By the Markov property, we have
\begin{equation*}
  J_{2}(s,t)=e^{(\lambda_1/2)(s+t)}\langle T_sf, X_t\rangle.
\end{equation*}
By \eqref{1.13}, we get
\begin{eqnarray*}
      \P_{\delta_x}\langle T_{s}f,X_t\rangle^2=\int_0^tT_{t-u}[A(T_{u+s}(f))^2](x)\,du+T_t(T_{s}f)^2(x).
\end{eqnarray*}
First, we consider the case $\gamma(f)<\infty$.
Since $u+s\ge s>10t_0$, by \eqref{1.23}, we get
\begin{equation}\label{9.1}
  |T_{u+s}f(x)|^2\lesssim (u+s)^{2\tau(f)}e^{-2\Re_{\gamma(f)}(u+s)}b_{4t_0}(x).
\end{equation}
Thus, for $t>10t_0$, we have
\begin{eqnarray}
 &&\int_{0}^{t-2t_0}T_{t-u}[A(T_{s+u}f)^2](x)\,du\nonumber\\
 &\lesssim& e^{-2s\Re_{\gamma(f)}}\int_0^{t-2t_0}\left(u+s\right)^{2\tau(f)}e^{-2\Re_{\gamma(f)}u}T_{t-u}(b_{4t_0})(x)\,du\nonumber\\
 &\lesssim & e^{-2s\Re_{\gamma(f)}}\int_0^{t-2t_0}\left(u+s\right)^{2\tau(f)}e^{-2\Re_{\gamma(f)}u}e^{-\lambda_1(t-u)}\,dub_{t_0}(x)^{1/2}.\label{4.4}\\
 &\lesssim & e^{-\lambda_1 t}e^{-2\Re_{\gamma(f)}s}
 \left(\int_0^{t-2t_0}u^{2\tau(f)}e^{(\lambda_1-2\Re_{\gamma(f)})u}\,du
 +s^{2\tau(f)}\int_0^{t-2t_0}e^{(\lambda_1-2\Re_{\gamma(f)})u}\,du\right)b_{t_0}(x)^{1/2}\nonumber\\
 &\lesssim &s^{2\tau(f)}e^{-\lambda_1 t}e^{-2\Re_{\gamma(f)}s}b_{t_0}(x)^{1/2}.\nonumber
\end{eqnarray}
The second inequality above follows from  \eqref{2.8}.
And by \eqref{9.1} and \eqref{1.37}, we have
\begin{eqnarray}\label{9.7}
  &&\int_{t-2t_0}^tT_{t-u}[A(T_{s+u}f)^2](x)\,du\nonumber\\
  &\lesssim & (t+s)^{2\tau(f)}e^{-2\Re_{\gamma(f)}(t+s-2t_0)}\int_{t-2t_0}^tT_{t-u}(b_{4t_0})(x)\,du\nonumber\\
  &\lesssim &(t+s)^{2\tau(f)}e^{-2\Re_{\gamma(f)}(t+s)}b_{t_0}(x)^{1/2}.
\end{eqnarray}
By \eqref{1.23}, we get that
$|T_{s}f(x)|^2\lesssim s^{2\tau(f)}e^{-2\Re_{\gamma(f)}s}b_{t_0}(x).$
Thus, we have
\begin{equation}\label{4.3}
  T_t(T_{s}f)^2(x)\lesssim s^{2\tau(f)}e^{-\lambda_1t}e^{-2\Re_{\gamma(f)}s}b_{t_0}(x)^{1/2}.
\end{equation}
Consequently, we have
\begin{equation}\label{4.49}
  \P_{\delta_x}\langle T_sf,X_t\rangle^2\lesssim (t+s)^{2\tau(f)}e^{-2\Re_{\gamma(f)}(t+s)}b_{t_0}(x)^{1/2}+s^{2\tau(f)}e^{-\lambda_1t}e^{-2\Re_{\gamma(f)}s}b_{t_0}(x)^{1/2}.
\end{equation}
Therefore, we have
\begin{equation}\label{6.7}
  \limsup_{t\to\infty}\P_{\delta_x}J_2(s,t)^2=
 \limsup_{t\to\infty}e^{\lambda_1(t+s)}\P_{\delta_x}\langle T_sf,X_t\rangle^2
  \lesssim s^{2\tau(f)}e^{(\lambda_1-2\Re_{\gamma(f)})s}b_{t_0}(x)^{1/2}.
\end{equation}
 Similarly, for the case $\gamma(f)=\infty,$ we have
\begin{equation}\label{4.49'}
  \P_{\delta_x}\langle T_sf,X_t\rangle^2\lesssim b_{t_0}(x)^{1/2}+e^{-\lambda_1t}b_{t_0}(x)^{1/2}.
\end{equation}
Thus,
\begin{equation}\label{6.7'}
  \limsup_{t\to\infty}\P_{\delta_x}J_2(s,t)^2=
 \limsup_{t\to\infty}e^{\lambda_1(t+s)}\P_{\delta_x}\langle T_sf,X_t\rangle^2
  \lesssim e^{\lambda_1s}b_{t_0}(x)^{1/2}.
\end{equation}
Combining \eqref{6.7} and \eqref{6.7'}, we get
\begin{equation}\label{6.7''}
   \limsup_{s\to\infty}\limsup_{t\to\infty}\P_{\delta_x}J_2(s,t)^2=0.
\end{equation}

Next we consider $J_1(s,t)$.
We define an ${\mathbb R}^2$-valued random variable $U_2(s,t)$ by
\begin{eqnarray*}
U_2(s,t):=\left(e^{\lambda_1 t}\langle \phi_1, X_t\rangle, J_1(s,t) \right).
\end{eqnarray*}
Let $V_s(x):={\V}{\rm ar}_{\delta_x}Y_1(s)$. We claim that, for any $x\in E$, under $\P_{\delta_x}$,
\begin{equation}\label{6.1}
  U_2(s,t)\stackrel{d}{\to}\left(W_\infty, \sqrt{W_\infty}G_1(s)\right), \quad \mbox{ as } t\to\infty,
\end{equation}
where $G_1(s)\sim\mathcal{N}(0,\sigma^2_f(s))$ is independent of $W_\infty$ and $\sigma^2_f(s)=\langle V_s,\phi_1\rangle$.
Denote the characteristic function of $U_2(s,t)$ under $\P_{\delta_x}$ by
$\kappa(\theta_1,\theta_2,s,t)$:
\begin{eqnarray}
  \kappa(\theta_1,\theta_2,s,t)
  &=&\P_{\delta_x}\left(\exp\left\{i\theta_1e^{\lambda_1 t}\langle \phi_1, X_t\rangle+i\theta_2e^{(\lambda_1/2) t}\sum_{u\in\cL_t}(Y_1^{u,t}(s)-y_1^{u,t}(s))\right\}\right)\nonumber\\
  &=& \P_{\delta_x}\left(\exp\{i\theta_1e^{\lambda_1 t}\langle \phi_1, X_t\rangle\}\prod_{u\in\cL_t}h_s(z_u(t),e^{(\lambda_1/2) t}\theta_2)\right),\label{8.3}
\end{eqnarray}
where
$$
h_s(x,\theta):=\P_{\delta_x}e^{i\theta(Y_1(s)-\P_{\delta_x}Y_1(s))}.
$$
Let $t_k,m_k\to\infty$, as $k\to\infty$, and $a_{k,j}\in E$, $j=1,2,\cdots m_k$. Now we consider
\begin{equation}\label{6.16}
  S_k:=e^{\lambda_1t_k/2}\sum_{j=1}^{m_k}(Y_{k,j}-y_{k,j}),
\end{equation}
where $Y_{k,j}$ has the same law as $Y_1(s)$ under $\P_{\delta_{a_{k,j}}}$ and $y_{k,j}=\P_{\delta_{a_{k,j}}}Y_1(s)$. Further, $Y_{k,j},j=1,2,\dots$ are independent. Suppose the following Lindeberg conditions hold:
\begin{description}
\item{(i)} as $k\to\infty$,
    $$e^{\lambda_1t_k}\sum_{j=1}^{m_k}\E(Y_{k,j}-y_{k,j})^2=e^{\lambda_1t_k}\sum_{j=1}^{m_k}V_s(a_{k,j})\to\sigma^2;$$
  \item{(ii)} for any $\epsilon>0$,
   \begin{eqnarray*}
     && e^{\lambda_1t_k}\sum_{j=1}^{m_k}\E\left(|Y_{k,j}-y_{k,j}|^2,|Y_{k,j}-y_{k,j}|>\epsilon e^{-\lambda_1 t_k/2}\right) \\
     &=& e^{\lambda_1t_k}\sum_{j=1}^{m_k}g(a_{k,j},s,t_k)\to 0,\quad \mbox{as }k\to\infty
  \end{eqnarray*}
  where
  $g(x,s,t)=\P_{\delta_x}\left(|Y_1(s)-\P_{\delta_x}Y_1(s)|^2,|Y_1(s)-\P_{\delta_x}Y_1(s)|>\epsilon e^{-\lambda_1 t/2}\right).$
\end{description}
Then using the Lindeberg-Feller theorem,
we have $S_k\stackrel{d}{\to}\mathcal{N}(0,\sigma^2)$, which implies
\begin{equation}\label{6.17}
  \prod_{j=1}^{m_k}h_s(a_{k,j},e^{\lambda_1t_k/2}\theta)\to e^{-\frac{1}{2}\sigma^2\theta^2}.
\end{equation}
By \eqref{var}, we get
$V_s\in L^2(E,m)\cap L^4(E,m)$.
So using Remark \ref{rem:large}, we have
\begin{equation}\label{6.18}
  e^{\lambda_1t}\sum_{u\in \cL_t}V_s(z_u(t))=e^{\lambda_1t}\langle V_s,X_t\rangle
  \to\langle V_s,\phi_1\rangle W_\infty,\quad \mbox{  in probability}, \mbox{as }t\to\infty.
\end{equation}

We note that $g(x,s,t)\downarrow0$ as  $t\uparrow\infty$ and $g(x,s,t)\leq V_s(x)$ for any $x\in E$.
Thus by \eqref{2.8} we have for any $x\in E$,
\begin{equation*}
 e^{\lambda_1t}\P_{\delta_x}\langle g(\cdot,s,t),X_t\rangle\lesssim \|g(\cdot,s,t)\|_2b_{t_0}(x)^{1/2}\to 0,
 \quad \mbox{as}\quad t\to\infty,
\end{equation*}
which implies
\begin{equation}\label{6.19}
  e^{\lambda_1t}\sum_{u\in \cL_t}g(z_u(t),s,t)\to 0,\quad \mbox{as}\quad t\to\infty,
\end{equation}
in $\P_{\delta_x}$-probability.
Therefore, for any sequence $s_k\to\infty$, there exists a subsequence $s_k'$ such that,
if we let $t_k=s_k'$, $m_k=|X_{s_k'}|$ and
$\{a_{k,j},j=1,2\cdots m_k\}=\{z_u(s_k'),u\in \cL_{s_k'}\}$,
then the Lindeberg conditions hold $\P_{\delta_x}$-a.s. for any $x\in E$, which implies
\begin{equation}\label{6.25}
  \lim_{k\to\infty}\prod_{u\in\cL_{s_k'}}h_s(z_u(s_k'),e^{(\lambda_1/2) s_k'}\theta_2)
  =\exp\left\{-\frac{1}{2}\theta_2^2\langle V_s,\phi_1\rangle W_\infty\right\},\quad\P_{\delta_x}\mbox{-a.s.}
\end{equation}
Consequently, we have
\begin{equation}\label{6.20}
  \lim_{t\to\infty}\prod_{u\in\cL_t}h_s(z_u(t),e^{(\lambda_1/2) t}\theta_2)
  =\exp\left\{-\frac{1}{2}\theta_2^2\langle V_s,\phi_1\rangle W_\infty\right\},\quad\mbox{in probability}.
\end{equation}
Hence by the dominated convergence theorem, we get
\begin{equation}\label{6.21}
 \lim_{t\to\infty}\kappa(\theta_1,\theta_2,s,t)= \P_{\delta_x}\exp\left\{i\theta_1W_\infty\right\}
  \exp\left\{-\frac{1}{2}\theta_2^2\langle V_s,\psi_1\rangle_m W_\infty\right\},
\end{equation}
which implies our claim \eqref{6.1}.
Thus, we easily get that, for any $x\in E$, under $\P_{\delta_x}$,
\begin{eqnarray*}
U_3(s,t):=\left(e^{\lambda_1 (t+s)}\langle \phi_1, X_{t+s}\rangle,J_1(s,t) \right)
\stackrel{d}{\to}\left(W_\infty, \sqrt{W_\infty}G_1(s)\right), \quad \mbox{ as } t\to\infty.
\end{eqnarray*}
By \eqref{small}, we have
$
\lim_{s\to\infty}\langle V_s,\psi_1\rangle_m=\sigma^2_f.
$
Let $G_1(f)$ be a  $\mathcal{N}(0,\sigma_f^2)$ random variable independent of $W_\infty$.
Then
\begin{equation}\label{6.22}
\lim_{s\to\infty}\beta(G_1(s),G_1(f))=0.
\end{equation}
Let $\mathcal{D}(s+t)$ and $\widetilde{\mathcal{D}}(s,t)$ be the distributions of $U_1(s+t)$ and $U_3(s,t)$
respectively, and let $\mathcal{D}(s)$ and $\mathcal{D}$ be the distributions of $(W_\infty, \sqrt{W_\infty}G_1(s))$
and $(W_\infty, \sqrt{W_\infty}G_1(f))$ respectively. Then, using \eqref{5.20}, we have
\begin{eqnarray}\label{6.23}
  \limsup_{t\to\infty}\beta(\mathcal{D}(s+t),\mathcal{D})&\leq&
  \limsup_{t\to\infty}[\beta(\mathcal{D}(s+t),\widetilde{\mathcal{D}}(s,t))+\beta(\widetilde{\mathcal{D}}(s,t),\mathcal{D}(s))
  +\beta(\mathcal{D}(s),\mathcal{D})]\nonumber\\
 &\leq &\limsup_{t\to\infty}(\P_{\delta_x}J_2(s,t)^2)^{1/2}+0+\beta(\mathcal{D}(s),\mathcal{D}).
\end{eqnarray}
Using this and the definition of $\limsup_{t\to\infty}$, we easily get that
$$
\limsup_{t\to\infty}\beta (\mathcal{D}(t),\mathcal{D})=\limsup_{t\to\infty}\beta (\mathcal{D}(s+t),\mathcal{D})
\le \limsup_{t\to\infty}(\P_{\delta_x}J_2(s,t)^2)^{1/2}+\beta(\mathcal{D}(s),\mathcal{D}).
$$
Letting $s\to\infty$, we get $ \limsup_{t\to\infty}\beta(\mathcal{D}(t),\mathcal{D})=0$.
The proof is now complete. \hfill$\Box$

\begin{lemma}\label{lem:5.5}
 Assume $f(x)=\sum_{j: \lambda_1=2\Re_j}(\Phi_j(x))^Tb_j\in\C_c$, where $b_j\in \mathbb{C}^{n_j}$.
Define
$$
  S_tf(x):=t^{-(1+2\tau(f))/2}e^{(\lambda_1/2) t}(\langle f,X_t\rangle-{T}_tf(x)), \qquad (t, x)\in (0, \infty)\times E.
$$
Then for any $c>0$, $\delta>0$ and $x\in E$, we have
\begin{equation}\label{4.5}
 \lim_{t\to\infty}\mathbb{P}_{\delta_x}\left(|S_tf(x)|^2;|S_tf(x)|>ce^{\delta t }\right)=0.
\end{equation}
\end{lemma}

\textbf{Proof:}\quad
In this proof, we always assume $t>10t_0$.
For each $j$, define
$$S_{j,t}(x):=t^{-(1+2\tau(f))/2}e^{\lambda_1t/2}\left(\langle \Phi_j^T,X_t\rangle-e^{-\lambda_jt}(\Phi_j(x))^TD_j(t)\right).$$
Thus, $S_tf(x)=\sum_{j: \lambda_1=2\Re_j}S_{j,t}(x)b_j$.
Using the fact that for every $n\ge 1$,
\begin{equation}\label{e:elem}
\left|\sum_{l=1}^n x_l\right|^2{\bf 1}_{|\sum_{l=1}^n x_l|^2>M}\le n\sum_{l=1}^n |x_l|^2{\bf 1}_{|x_l|^2>M/n},
\end{equation}
we see that, to prove \eqref{4.5}, it suffices to show that, as $t\to\infty$,
 $$
 F(t, x,b_j):=\P_{\delta_x}\left(|S_{j,t}(x)b_j|^2;|S_{j,t}(x)b_j|>ce^{\delta t }\right)\to 0.
 $$
Choose an integer $n_0>2t_0$. We write $t=l_tn_0+\epsilon_t$, where $l_t\in \mathbb{N}$ and $0\le\epsilon_t<n_0$.
By \eqref{T-Jordan}, we easily get $T_u(\Phi_j^T)(x)=e^{-\lambda_ju}\Phi_j(x)^TD_j(t)$.
Since $\lambda_1=2\Re_j$,
for any $(t, x)\in (0, \infty)\times E$, we have
\begin{eqnarray}\label{4.8}
S_{j,t+n_0}(x)&=&\left(\frac{1}{t+n_0}\right)^{1/2+\tau(f)}e^{\lambda_1 (t+n_0)/2}\left(\langle \Phi_j^T, X_{t+n_0}\rangle-
  \langle e^{-\lambda_jn_0}\Phi_j^T, X_t\rangle D_j(n_0)\right)\nonumber\\
  &&+\left(\frac{1}{t+n_0}\right)^{1/2+\tau(f)}e^{-i\Im_jn_0}e^{\lambda_1t/2}\left(\langle \Phi_j^T, X_t\rangle-e^{-\lambda_jt}(\Phi_j(x))^TD_j(t)\right)D_j(n_0)\nonumber\\
   &=&\left(\frac{1}{t+n_0}\right)^{1/2+\tau(f)}R_j(t)+e^{-i\Im_jn_0}\left(\frac{t}{t+n_0}\right)^{1/2+\tau(f)}S_{j,t}(x)D_j(n_0),
\end{eqnarray}
where
$$
R_j(t):=e^{(\lambda_1/2) (t+n_0)}\left(\langle \Phi_j^T,X_{t+n_0}\rangle-\langle e^{-\lambda_j}\Phi_j^T, X_t\rangle D_j(n_0)\right).
$$
Hence, for any $(t, x)\in (0, \infty)\times E$, we have
 \begin{eqnarray*}
 &&F(t+n_0, x,b_j)\\
 &\leq& \P_{\delta_x}\left(|S_{j,t+n_0}(x)b_j|^2;|S_{j,t}(x)D_j(n_0)b_j|>ce^{\delta t }\right)\\
 &&+\P_{\delta_x}\left(|S_{j,t+n_0}(x)b_j|^2;|S_{j,t}(x)D_j(n_0)b_j|\leq ce^{\delta t},|S_{j,t+n_0}(x)b_j|>ce^{\delta (t+n_0)}\right)\\
 &=:&M_1(t, x)+M_2(t, x).
 \end{eqnarray*}
 Put
 \begin{eqnarray*}
    A_1(t,x,b_j) &:=& \{|S_{j,t}(x)D_j(n_0)b_j|>ce^{\delta t }\},\\
    A_2(t,x,b_j) &:=& \{|S_{j,t}(x)D_j(n_0)b_j|\leq ce^{\delta t},|S_{j,t+n_0}(x)b_j|>ce^{\delta (t+n_0)}\}\\
    \mbox{and}\qquad \qquad \qquad && \\
    A(t,x,b_j)&:=&A_1(t,x,b_j)\cup A_2(t,x,b_j).
 \end{eqnarray*}
 Since
$A_1(t, x,b_j)\in\mathcal{F}_t$ and $\P_{\delta_x}(R_j(t) |\mathcal{F}_t)$=0 for any $(t, x)\in (0, \infty)\times E$,
we have by \eqref{4.8} that
 \begin{eqnarray*}
 M_1(t, x)&=&\left(\frac{1}{t+n_0}\right)^{1+2\tau(f)}\P_{\delta_x} \left(|R_j(t)b_j|^2;A_1(t,x,b_j)\right)
  +\left(\frac{t}{t+n_0}\right)^{1+2\tau(f)}F(t, x,D_j(n_0)b_j)
 \end{eqnarray*}
 and
 \begin{eqnarray*}
   M_2(t, x)&\leq&2\left(\frac{1}{t+n_0}\right)^{1+2\tau(f)}\P_{\delta_x}\left(|R_j(t)b_j|^2;A_2(t,x,b_j)\right)\\
   &&+2\left(\frac{t}{t+n_0}\right)^{1+2\tau(f)}\P_{\delta_x}\left(|S_{j,t}(x)D_j(n_0)b_j|^2;A_2(t, x.b_j)\right).
 \end{eqnarray*}
 Thus,  for any $(t, x)\in (0, \infty)\times E$, we have
 \begin{eqnarray}\label{4.9}
F(t+n_0, x,b_j)&\leq&\left(\frac{t}{t+n_0}\right)^{1+2\tau(f)}F(t, x,D_j(n_0)b_j)\nonumber\\
&&\quad +\left(\frac{1}{t+n_0}\right)^{1+2\tau(f)}(F_1(t, x,b_j)+F_2(t, x,b_j)),
 \end{eqnarray}
 where
 \begin{eqnarray*}
   F_1(t, x,b_j) &:=& 2\P_{\delta_x}\left(|R_j(t)b_j|^2;A_1(t, x,b_j)\cup A_2(t, x,b_j)\right),\\
   F_2(t, x,b_j) &:=& 2t^{1+2\tau(f)}\P_{\delta_x}\left(|S_{j,t}(x)D_j(n_0)b_j|^2;A_2(t, x,b_j)\right).
 \end{eqnarray*}
 Iterating \eqref{4.9}, we get for $t$ large enough,
 \begin{eqnarray}\label{4.10}
   && F(t+n_0,x,b_j)\nonumber\\
   &\leq&\left(\frac{1}{t+n_0}\right)^{1+2\tau(f)}
    \sum_{m=5}^{l_t}\left(F_1(mn_0+\epsilon_t,x,D_j((l_t-m)n_0)b_j)\right)\nonumber\\
    &&\left(\frac{1}{t+n_0}\right)^{1+2\tau(f)}
    \sum_{m=5}^{l_t}\left(F_2(mn_0+\epsilon_t,x,D_j((l_t-m)n_0)b_j)\right)\nonumber\\
   && +\left(\frac{5n_0+\epsilon_t}{t+n_0}\right)^{1+2\tau(f)}
   F(5n_0+\epsilon_t,x,D_j((l_t-4)n_0)b_j)\nonumber\\
     &=:&L_1(t, x)+L_2(t, x)+\left(\frac{5n_0+\epsilon_t}{t+1}\right)^{1+2\tau(f)}
   F(5n_0+\epsilon_t,x,D_j((l_t-4)n_0)b_j).
 \end{eqnarray}

First, we  consider $L_1(t, x)$.
By the definition of $\tau(f)$, we have for $s>0$,
\begin{equation}\label{2.33}
  |D_j(s)b_j|_2\lesssim |D_j(s)b_j|_\infty\lesssim 1+s^{\tau(f)}.
\end{equation}
Thus, we have for $0\le s\le t$ and $t\ge 2t_0$,
\begin{equation}\label{R2}
  |R_j(s)D_j(t-s)b_j|^2\le |R_j(s)|_2^2|D_j(t-s)b_j|_2^2\lesssim t^{2\tau(f)}|R_j(s)|_2^2.
\end{equation}
It follows that for any $\epsilon\in (0, 1)$,
\begin{eqnarray}\label{L1}
  L_1(t,x)&\le& \frac{2}{t+n_0}\sum_{5\le m\le \epsilon l_t}\P_{\delta_x}\left(|R_j(mn_0+\epsilon_t)|_2^2\right)\nonumber\\
  &&+\frac{2}{t+n_0}\sum_{l_t\epsilon\le m\le l_t}\P_{\delta_x}\left(|R_j(mn_0+\epsilon_t)|_2^2; A(mn_0+\epsilon_t,x,D_j((l_t-m)n_0)b_j)\right)\nonumber\\
  &=:&L_{1,1}(t,x)+L_{1,2}(t,x).
\end{eqnarray}
By the definition of $R_j(s)$, we have
\begin{equation}\label{2.34}
  |R_j(s)|_2^2=e^{\lambda_1(s+n_0)}\sum_{l=1}^{n_j}|\langle \phi_l^{(j)},X_{s+n_0}\rangle-\langle T_{n_0}(\phi_l^{(j)}),X_s\rangle|^2.
\end{equation}
Note that
\begin{eqnarray*}
  |\langle \phi_l^{(j)},X_{s+n_0}\rangle-\langle T_{n_0}(\phi_l^{(j)}),X_s\rangle|^2&=&|\langle \Re(\phi_l^{(j)}),X_{s+n_0}\rangle-\langle T_{n_0}(\Re(\phi_l^{(j)})),X_s\rangle|^2 \\
   &+& |\langle \Im(\phi_l^{(j)}),X_{s+n_0}\rangle-\langle T_{n_0}(\Im(\phi_l^{(j)})),X_s\rangle|^2.
\end{eqnarray*}
Thus, we have
\begin{equation}\label{2.35}
  \P_{\delta_x}|\langle \phi_l^{(j)},X_{s+n_0}\rangle-\langle T_{n_0}(\phi_l^{(j)}),X_s\rangle|^2
=T_s({\V}{\rm ar}_{\delta_{\cdot}}\langle \Re(\phi_l^{(j)}),X_{n_0} \rangle)(x)+T_s({\V}{\rm ar}_{\delta_{\cdot}}\langle \Im(\phi_l^{(j)}),X_{n_0} \rangle)(x).
\end{equation}
Hence, by \eqref{2.8}, we get, for $s\ge 5n_0>2t_0$,
\begin{equation}\label{2.37}
  \P_{\delta_x}|R_j(s)|_2^2
=e^{\lambda_1(s+n_0)}\sum_{l=1}^{n_j}\P_{\delta_x}|\langle \phi_l^{(j)},X_{t+n_0}\rangle-\langle T_{n_0}(\phi_l^{(j)}),X_t\rangle|^2
\lesssim b_{t_0}(x)^{1/2}.
\end{equation}
Therefore, we have, for $(t,x)\in(5n_0,\infty)\times E$,
\begin{equation}\label{L11}
  L_{1,1}(t,x)\lesssim \epsilon b_{t_0}(x)^{1/2}.
\end{equation}
We claim that, for any $x\in E$,

 (i) \begin{equation}\label{2.44}
  \lim_{M\to\infty}\limsup_{s\to\infty}\P_{\delta_x}(|R_j(s)|_2^2; |R_j(s)|_2^2>M)=0
  ,\mbox{ and}
  \end{equation}

  (ii) \begin{equation}\label{Ato0}
   \sup_{t\epsilon\le s\le t}\P_{\delta_x}(A_1(s, x,D_j(t-s)b_j)\cup A_2(s, x,D_j(t-s)b_j))\to 0.
  \end{equation}
Using these two claims we get that, as $t\to\infty$,
\begin{eqnarray}
  &&L_{1,2}(t,x) \nonumber \\
  &\le& \frac{2}{t+n_0}\sum_{\epsilon l_t\le m\le l_t}\left(\P_{\delta_x}\left(|R_j(mn_0+\epsilon_t)|_2^2; |R_j(mn_0+\epsilon_t)|_2^2>M\right)\right.\nonumber\\
  &&\left.+M\P_{\delta_x}\left(A(mn_0+\epsilon_t,x,D_j((l_t-m)n_0)b_j)\right)\right)\nonumber \\
   &\lesssim& \sup_{s\ge t\epsilon} \P_{\delta_x}(|R_j(s)|_2^2; |R_j(s)|_2^2>M)
+M\sup_{t\epsilon\le s\le t}\P_{\delta_x}(A(s, x,D_j(t-s)b_j))\nonumber \\
   &\to&\limsup_{s\to\infty}\P_{\delta_x}(|R_j(s)|_2^2; |R_j(s)|_2^2>M).
\end{eqnarray}
Letting $M\to\infty$, we get
\begin{equation}\label{L12}
  \lim_{t\to\infty}L_{1,2}(t,x)=0.
\end{equation}
Now we prove the two claims.

(i) For $l=1,2,\cdots,n_j,$ define
$$R_{j,l,1}(s):=e^{\lambda_1 (s+n_0)/2}\langle \Re(\phi_l^j),X_{s+n_0}\rangle-\langle T_{n_0}(\Re(\phi_l^j)),X_s\rangle $$
and
$$R_{j,l,2}(s):=e^{\lambda_1 (s+n_0)/2}\langle \Im(\phi_l^j),X_{s+n_0}\rangle-\langle T_{n_0}(\Im(\phi_l^j)),X_s\rangle. $$
Using \eqref{e:elem} and \eqref{2.34}, we easily see that,
to prove \eqref{2.44},
we only need to show that, for $k=1,2$,
\begin{equation}\label{to-prove}
\lim_{M\to\infty}\limsup_{s\to\infty}\P_{\delta_x}(|R_{j,l,k}(s)|^2,|R_{j,l,k}(s)|^2>M)=0.
\end{equation}
Repeating the proof of \eqref{6.1} with $s=n_0$, we see that \eqref{6.1} is valid for
$f\in L^2(E,m)\cap L^4(E,m)$.
Thus, for $l=1,2,\cdots,n_j,$
as $s\to\infty$,
$$
R_{j,l,1}(s)\stackrel{d}{\rightarrow}\sqrt{W_\infty}G,
$$
where $G\sim \mathcal{N}(0,e^{\lambda_1n_0}\langle {\V}{\rm ar}_{\delta_\cdot}\langle \Re(\phi_l^j),X_{n_0} \rangle, \psi_1\rangle_m$.
And by \eqref{2.6}, we get, as $s\to\infty$,
\begin{equation}\label{2.39}
  \P_{\delta_x}(|R_{j,l,1}(s)|^2)
=e^{\lambda_1 (s+n_0)}T_s({\V}{\rm ar}_{\delta_{\cdot}}\langle \Re(\phi_l^j),X_{n_0} \rangle)(x)
\to e^{\lambda_1n_0}\langle ({\V}{\rm ar}_{\delta_{\cdot}}\langle \Re(\phi_l^j),X_{n_0} \rangle,\psi_1\rangle_m\phi_1(x).
\end{equation}
Let $h_{M}(r)=r$ on $[0,M-1]$, $h_M(r)=0$ on $[M,\infty]$, and let $h_M$ be linear on $[M-1,M]$.
By \eqref{2.39}, we have that for any $x\in E$,
\begin{eqnarray*}
  &&\limsup_{s\to\infty}\P_{\delta_x}(|R_{j,l,1}(s)|^2,|R_{j,l,1}(s)|^2>M)\le \limsup_{t\to\infty}\P_{\delta_x}(|R_{j,l,1}(s)|^2)-\P_{\delta_x}(h_M(|R_{j,l,1}(s)|^2))\\
  &=& e^{\lambda_1n_0}\langle ({\V}{\rm ar}_{\delta_{\cdot}}\langle \Re(\phi_l^j),X_{n_0} \rangle,\psi_1\rangle_m\phi_1(x)
-\P_{\delta_x}(h_M(W_\infty G^2)).
\end{eqnarray*}
By the monotone convergence theorem, we have that for any $x\in E$,
$$\lim_{M\to\infty}\P_{\delta_x}(h_M(W_\infty G^2))=\P_{\delta_x}(W_\infty G^2)=\P_{\delta_x}(W_\infty)\P_{\delta_x}(G^2)
=e^{\lambda_1n_0}\langle ({\V}{\rm ar}_{\delta_{\cdot}}\langle \Re(\phi_l^j),X_{n_0} \rangle,\psi_1\rangle_m\phi_1(x),$$
which implies
$$
  \lim_{M\to\infty}\limsup_{s\to\infty}\P_{\delta_x}(|R_{j,l,1}(s)|^2,|R_{j,l,1}(s)|^2>M)=0,
  $$
which says \eqref{to-prove} holds for $k=1$. Using similar arguments,
we get \eqref{to-prove} holds for $k=2$.

(ii)
Since $\tau(\phi_l^j)\le \nu_j$, by \eqref{3.33}, we get for $10t_0\le s$,
\begin{equation}\label{vars}
  \P_{\delta_x}|S_{j,s}(x)|_2^2\lesssim s^{1+2\nu_j}s^{-(1+2\tau(f))}\le s^{2\nu_j}.
\end{equation}
By \eqref{2.33}, we get, for $10t_0\le s\le t$,
\begin{equation}\label{3.34}
  \P_{\delta_x}|S_{j,s}(x)D_j(t+1-s)b_j|^2\lesssim s^{2\nu_j}(1+t^{2\tau(f)}).
\end{equation}
By Chebyshev's inequality and \eqref{3.34}, we have that, for any $x\in E$, as $t\to\infty$
\begin{eqnarray*}
  &&\sup_{t\epsilon\le s\le t}\P_{\delta_x}( A_1(s, x,D_{j}(t-s)))
\leq \sup_{t\epsilon\le s\le t}c^{-2}e^{-2\delta s }\P_{\delta_x}|S_{j,s}(x)D_j(t+1-s)b_j|^2\\
  &\lesssim &e^{-2\delta \epsilon t} t^{2\nu_j}(1+t^{2\tau(f)})\to 0.
\end{eqnarray*}
It is easy to see that, under $\P_{\delta_x}$, for any $t>0$,
\begin{equation}\label{4.24}
A_2(s,x, D_j(t-s)b_j)\subset
\left\{|R_j(s)D_j(t-s)b_j|>ce^{\delta s}\left(e^{\delta n_0}-1\right)s^{(2\tau(f)+1)/2}\right\}.
\end{equation}
By \eqref{R2} and \eqref{2.37}, we get
$$\P_{\delta_x}|R_j(s)D_j(t-s)b_j|^2\lesssim t^{2\tau(f)}b_{t_0}(x)^{1/2}.$$
Similarly, by Chebyshev's inequality, we have that, for any $x\in E$, as $t\to\infty$,
\begin{eqnarray*}
 &&\sup_{t\epsilon\le s\le t}\P_{\delta_x}A_2(s,x, D_j(t-s)b_j)\\
 &\leq& \sup_{t\epsilon\le s\le t}c^{-2}(e^{\delta n_0}-1)^{-2}e^{-2\delta s}s^{-(1+2\tau(f))}
\P_{\delta_x}|R_j(s)D_j(t-s)b_j|^2\\
 &\lesssim& e^{-2\delta \epsilon t}(t\epsilon)^{-(1+2\tau(f))}t^{2\tau(f)}\to 0.
\end{eqnarray*}
Thus we have finished proving the two claims.
Therefore, by \eqref{L11} and \eqref{L12}, we get
$$\limsup_{t\to\infty}L_1(t,x)\lesssim \epsilon   b_{t_0}(x)^{1/2}.$$
Letting $\epsilon\to0,$ we get
\begin{equation}\label{L1to0}
  \lim_{t\to\infty}L_1(t,x)=0.
\end{equation}

Now we consider $L_2(t, x)$.
By \eqref{4.24}, we have that for any $x\in E$,
\begin{eqnarray*}
&&F_2(s, x, D_j(t-s)b_j)\\
&=&2s^{(1+2\tau(f))}\P_{\delta_x}\left(|S_{j,s}(x)D_j(t+n_0-s)b_j|^2;A_2(s,x,D_j(t-s)b_j)\right)\\
 &\leq& 2s^{(1+2\tau(f))}ce^{\delta s}
\P_{\delta_x}\left(|S_{j,s}(x)D_j(t+n_0-s)b_j|;|R_j(s)D_j(t-s)b_j|>ce^{\delta s}(e^{\delta n_0}-1)s^{(2\tau(f)+1)/2}\right)\\
    &\leq& 2c^{-1}(e^{\delta n_0}-1)e^{-\delta s}\P_{\delta_x}\left(|S_{j,s}(x)D_j(t+n_0-s)b_j|\cdot|R_j(s)D_j(t-s)b_j|^2\right)\\
   &\lesssim& e^{-\delta s}e^{\lambda_1(s+n_0)}t^{\tau(f)}\P_{\delta_x}\left(|S_{j,s}(x)|_2
\langle {\V}{\rm ar}_{\delta_\cdot}(\langle \Phi_j^TD_j(t-s)b_j,X_{n_0}\rangle),X_s\rangle\right)\\
   &\lesssim& e^{-\delta s}t^{\tau(f)}\sqrt{\P_{\delta_x}|S_{j,s}(x)|_2^2}\sqrt{e^{2\lambda_1s}
\P_{\delta_x}\left(\langle {\V}{\rm ar}_{\delta_\cdot}(\langle \Phi_j^TD_j(t-s)b_j,X_{n_0}\rangle),X_s\rangle^2\right)}.
\end{eqnarray*}
By \eqref{3.33} and \eqref{phi}, we get for $s\le t$,
\begin{eqnarray*}
  &&{\V}{\rm ar}_{\delta_x}(\langle \Phi_j^TD_j(t-s)b_j,X_{n_0}\rangle\le \P_{\delta_x}|\langle \Phi_j^TD_j(t-s)b_j,X_{n_0}\rangle|^2
\lesssim t^{2\tau(f)}\P_{\delta_x}\langle b_{t_0}(x)^{1/2},X_{n_0}\rangle^2.
\end{eqnarray*}
Thus by \eqref{vars} and \eqref{1.51}, we have for $5n_0\le s\le t$
$$
F_2(s, x, D_j(t-s)b_j)\lesssim e^{-\delta s}t^{2\tau(f)}s^{\nu_j}\sqrt{e^{2\lambda_1s}
\P_{\delta_x}\left(\langle b_{t_0}(x)^{1/2},X_s\rangle^2\right)}
\lesssim  e^{-\delta s}t^{2\tau(f)}s^{\nu_j}.
$$
Thus, we get, as $t\to\infty$,
\begin{equation}\label{L2to0}
  L_2(t,x)\lesssim \frac{1}{t+n_0}\sum_{m=5}^{l_t}e^{-\delta (mn_0+\epsilon_t)}(mn_0+\epsilon_t)^{(1+2\nu_j)/2}
\le \frac{1}{t+n_0}\sum_{m=5}^{l_t}e^{-\delta mn_0}((m+1)n_0)^{(1+2\nu_j)/2}\to 0.
\end{equation}
To finish the proof,  we only need to show that for any $x\in E$,
\begin{equation}\label{L3}
\lim_{t\to\infty}\left(\frac{5n_0+\epsilon_t}{t+n_0}\right)^{1+2\tau(f)}F(5n_0+\epsilon_t,x,D_j((l_t-4)n_0)b_j)=0.
\end{equation}
By \eqref{2.33} and \eqref{vars}, we get that for any $x\in E$,
\begin{eqnarray*}
  &&(5n_0+\epsilon_t)^{1+2\tau(f)}F(5n_0+\epsilon_t,x,D_j((l_t-4)n_0)b_j)\\
    &\le&
  (6n_0)^{1+2\tau(f)}\sup_{5n_0\le s\le 6n_0}\P_{\delta_x}|S_{j,s}(x)D_j((l_t-4)n_0)b_j|^2\lesssim t^{2\tau(f)}(6n_0)^{2\nu_j},
\end{eqnarray*}
which implies \eqref{L3}.

The proof is now complete. \hfill$\Box$

\begin{lemma}\label{lem:5.6}
Assume that $f\in \C_s$ and $h\in\C_c$.
Define
$$
  Y_1(t):=e^{\lambda_1t/2}\left(\langle f,X_t\rangle-T_tf(x)\right), \quad
  Y_2(t):=t^{-(1+\tau(h)/2)}e^{\lambda_1t/2}\left(\langle h,X_t\rangle-T_th(x)\right),\quad t>0,
$$
and $Y_t:=Y_1(t)+Y_2(t)$, $ t>0.$
Then for any $c>0$, $\delta>0$ and $x\in E$, we have
\begin{equation}\label{4.6}
 \lim_{t\to\infty}\mathbb{P}_{\delta_x}\left(|Y_t|^2;|Y_t|>ce^{\delta t }\right)=0.
\end{equation}
\end{lemma}
\textbf{Proof:}
By \eqref{e:elem} and Lemma \ref{lem:5.5}, it suffices to show that
\begin{equation}\label{3.55}
  \lim_{t\to\infty}\mathbb{P}_{\delta_x}\left(|Y_1(t)|^2;|Y_1(t)|>ce^{\delta t }\right)=0.
\end{equation}
If $\gamma(f)<\infty$, by \eqref{1.23}, we get, as $t\to\infty$,
$$e^{\lambda_1t/2}|T_tf(x)|\lesssim t^{\tau(f)}e^{(\lambda_1/2-\Re_{\gamma(f)})t}b_{t_0}(x)^{1/2}\to 0.$$
If $\gamma(f)=\infty$, by \eqref{1.23'}, we get, as $t\to\infty$, $e^{\lambda_1t/2}|T_tf(x)|\lesssim e^{\lambda_1t/2}b_{t_0}(x)^{1/2}\to 0.$
Thus, by Lemma \ref{lem:small}, $Y_1(t)\stackrel{d}{\to}\sqrt{W_\infty}G_1(f)$.
By Lemma \ref{lem:2.2}, we have
$$\lim_{t\to\infty}\mathbb{P}_{\delta_x}\left(|Y_1(t)|^2\right)=\sigma^2_f\phi_1(x).$$
Thus, for any $M>0$, we have
\begin{eqnarray*}
  \mathbb{P}_{\delta_x}\left(|Y_1(t)|^2;|Y_1(t)|>ce^{\delta t }\right)
  &\le& \mathbb{P}_{\delta_x}\left(|Y_1(t)|^2;|Y_1(t)|>M\right)
+M^2\mathbb{P}_{\delta_x}\left(|Y_1(t)|>ce^{\delta t }\right)\\
&=:& I_1(t, x, M)+I_2(t, x, M).
\end{eqnarray*}
Let $h_{M}(r)=r$ on $[0,M-1]$,
$h_M(r)=0$ on $[M,\infty]$, and let $h_M$ be linear on $[M-1,M]$.
Then
$$
\limsup_{t\to\infty}I_1(t, x, M)
\le \limsup_{t\to\infty}\mathbb{P}_{\delta_x}\left(|Y_1(t)|^2\right)-\mathbb{P}_{\delta_x}(h_M(|Y_1(t)|)^2)
=\sigma_{f}^2\phi_1(x)-\mathbb{P}_{\delta_x}(h_M(|G_1(f)\sqrt{W_\infty}|)^2.
$$
By Chebyshev's inequality,
we have, as $t\to\infty$,
$$
I_2(t, x, M)\le M^2c^{-2}e^{-2\delta t}\mathbb{P}_{\delta_x}\left(|Y_1(t)|^2\right)\to 0.
$$
Thus, we have
$$
\limsup_{t\to\infty}\mathbb{P}_{\delta_x}\left(|Y_1(t)|^2;|Y_1(t)|>ce^{\delta t }\right)
\le \sigma_{f}^2\phi_1(x)-\mathbb{P}_{\delta_x}(h_M(|G_1(f)\sqrt{W_\infty}|)^2.
$$
Letting $M\to\infty$, by the monotone convergence theorem, we have that for any $x\in E$,
$$
\lim_{M\to\infty}\mathbb{P}_{\delta_x}(h_M(|G_1(f)\sqrt{W_\infty}|)^2
=\mathbb{P}_{\delta_x}(G_1(f)^2W_\infty)=\sigma_{f}^2\phi_1(x),
$$
which implies \eqref{3.55}.
The proof is now complete.\hfill$\Box$

\begin{lemma}\label{lem:cs}
Assume that $f\in \C_s$ and $h\in\C_c$.
Then
\begin{equation}\label{cs}
  \left(e^{\lambda_1t}\langle\phi_1,X_t\rangle, t^{-(1+2\tau(h))/2}e^{\lambda_1t/2}\langle h,X_t\rangle,
  e^{\lambda_1t/2}\langle f,X_t\rangle\right)
  \stackrel{d}{\to}\left(W_\infty, \sqrt{W_\infty}G_2(h),\sqrt{W_\infty}G_1(f)\right),
\end{equation}
where $G_2(h)\sim\mathcal{N}(0,\rho^2_h)$ and $G_1(f)\sim\mathcal{N}(0,\sigma^2_f)$.
Moreover, $W_\infty$, $G_2(h)$ and $G_1(f)$ are independent.
\end{lemma}
\textbf{Proof:}\quad
In this proof, we always assume $t>10t_0$, $f\in \C_s$ and $h\in\C_c$.
We define an ${\mathbb R}^3$-valued random variable by
\begin{equation*}
  U_1(t):=\left(e^{\lambda_1t}\langle\phi_1,X_t\rangle, t^{-(1+2\tau(h))/2}e^{\lambda_1t/2}\langle h,X_t\rangle,
  e^{\lambda_1t/2}\langle f,X_t\rangle\right).
\end{equation*}
For $n>2$, we define
\begin{equation*}
  U_1(nt)=\left(e^{\lambda_1nt}\langle\phi_1,X_{nt}\rangle, (nt)^{-(1+2\tau(h))/2}e^{\lambda_1nt/2}\langle h,X_{nt}\rangle,
  e^{\lambda_1nt/2}\langle f,X_{nt}\rangle\right).
\end{equation*}
Now we define another ${\mathbb R}^3$-valued random variable $U_2(n,t)$ by
\begin{eqnarray*}
&&U_2(n,t)\\
&:=&\left(e^{\lambda_1t}\langle\phi_1,X_t\rangle,
\frac{e^{\lambda_1nt/2}(\langle h,X_{nt}\rangle-\langle T_{(n-1)t}h,X_t\rangle)}{((n-1)t)^{(1+2\tau(h))/2}},
e^{\lambda_1nt/2}(\langle f,X_{nt}\rangle-\langle T_{(n-1)t}f,X_t\rangle)\right).
\end{eqnarray*}
We claim that
\begin{equation}\label{9.5}
  U_2(n,t)\stackrel{d}{\to}\left(W_\infty, \sqrt{W_\infty}G_2(h),\sqrt{W_\infty}G_1(f) \right), \quad \mbox{ as } t\to\infty.
\end{equation}
Denote the characteristic function of $U_2(n,t)$ under $\P_\mu$ by
$\kappa_2(\theta_1,\theta_2,\theta_3,n,t)$.
Define
$$
  Y_1^{u,t}(s):=e^{\lambda_1s/2}\langle f,X_{s}^{u,t}\rangle,\quad
  Y_2^{u,t}(s):=s^{-(1+2\tau(h))/2}e^{\lambda_1s/2}\langle h,X_{s}^{u,t}\rangle,\quad s,t>0.
$$
We also define
$$
  Y_1(s):=e^{\lambda_1s/2}\langle f,X_{s}\rangle,\quad Y_2(s):=s^{-(1+2\tau(h))/2}e^{\lambda_1s/2}\langle h,X_{s}\rangle
$$
and
$$Y_s(\theta_2,\theta_3):=\theta_2Y_2(s)+\theta_3Y_1(s).$$
Given $\mathcal{F}_t$, for $k=1,2$, $Y_k^{u, t}(s)$ has the same distribution as $Y_k(s)$ under $\P_{\delta_{z_u(t)}}$.
Thus, for $k=1,2$,
$$y_k^{u, t}(s):=\P_{\delta_x}(Y^{u, t}_k(s)|\mathcal{F}_t)=\P_{\delta_{z_u(t)}}Y_k(s).$$
Thus, by \eqref{3.22}, we have
\begin{eqnarray}\label{3.11}
  U_2(n,t)
=\left(e^{\lambda_1t}\langle\phi_1,X_t\rangle,
e^{\lambda_1t/2}\sum_{u\in\cL_t}(Y_2^{u,t}((n-1)t)-y_2^{u,t}((n-1)t)),\right.&&\nonumber\\
\left.e^{\lambda_1t/2}\sum_{u\in\cL_t}(Y_1^{u,t}((n-1)t)-y_1^{u,t}((n-1)t))\right)&&.
\end{eqnarray}
Let $h(s,x,\theta,\theta_2,\theta_3)
=\P_{\delta_x}(\exp\{i\theta(Y_s(\theta_2,\theta_3)-\P_{\delta_x}Y_s(\theta_2,\theta_3))\}).$
Thus, we get
\begin{equation}\label{4.1}
  \kappa_2(\theta_1,\theta_2,\theta_3,n,t)
  =\P_{\delta_x}\left(\exp\{i\theta_1e^{\lambda_1 t}\langle \phi_1, X_t\rangle\}
  \prod_{u\in\cL_t}h\left((n-1)t,z_u(t),e^{\lambda_1t/2}, \theta_2,\theta_3\right)\right).
\end{equation}
Let $t_k,m_k\to\infty$, as $k\to\infty$. Now we consider
\begin{equation}\label{3.16}
  S_k:=e^{\lambda_1t_k/2}\sum_{j=1}^{m_k}(Y_{k,j}-y_{k,j}),
\end{equation}
where $Y_{k,j}$ has the same law as $Y_{(n-1)t_k}(\theta_2,\theta_3)$ under $\P_{\delta_{a_{k,j}}}$ and $y_{k,j}=\P_{\delta_{a_{k,j}}}Y_{(n-1)t_k}(\theta_2,\theta_3)$
with $a_{k, j}\in E$.
Further, for each positive integer $k$, $Y_{k,j},j=1,2,\dots$ are independent.
Denote $V_t^n(x):={\V}ar_{\delta_x}Y_{(n-1)t}(\theta_2,\theta_3)$.
Suppose the following Lindeberg conditions hold:
\begin{description}
  \item{(i)} as $k\to\infty$, $$e^{\lambda_1t_k}\sum_{j=1}^{m_k}\E(Y_{k,j}-y_{k,j})^2
  =e^{\lambda_1t_k}\sum_{j=1}^{m_k}V_{t_k}^n(a_{k,j})\to\sigma^2;$$
  \item{(ii)} for every $c>0$,
  \begin{eqnarray*}
     && e^{\lambda_1t_k}\sum_{j=1}^{m_k}\E\left(|Y_{k,j}-y_{k,j}|^2,|Y_{k,j}-y_{k,j}|>c e^{-\lambda_1 t_k/2}\right) \\
     &=& e^{\lambda_1t_k}\sum_{j=1}^{m_k}g_{(n-1)t_k}(a_{k,j},\theta_2,\theta_3)\to 0,\quad k\to\infty,
  \end{eqnarray*}
  where
  $$g_s(x,\theta_2,\theta_3)=\P_{\delta_{x}}\left(|Y_{s}(\theta_2,\theta_3)-\P_{\delta_{x}}Y_{s}(\theta_2,\theta_3)|^2,
     |Y_{s}(\theta_2,\theta_3)-\P_{\delta_{x}}Y_{s}(\theta_2,\theta_3)|>ce^{-\lambda_1 s/(2(n-1))}\right).$$
\end{description}
Then $S_k\stackrel{d}{\to}\mathcal{N}(0,\sigma^2)$, which implies
\begin{equation}\label{3.17}
  \prod_{j=1}^{m_k}h((n-1)t_k,a_{k,j},e^{\lambda_1t_k/2},\theta_2,\theta_3)\to e^{-\frac{1}{2}\sigma^2\theta^2},
  \quad\mbox{as }k\to\infty.
\end{equation}
By the definition of $Y_s$, we get
\begin{eqnarray*}
  V_t^n(x)&:=&{\V}ar_{\delta_x}Y_{(n-1)t}(\theta_2,\theta_3)
  =\theta_2^2 {\V}ar_{\delta_x}Y_2((n-1)t)+\theta_3^2 {\V}ar_{\delta_x}Y_1((n-1)t)\\
 &&+2\theta_2\theta_3((n-1)t)^{-(1+2\tau(h))/2}e^{\lambda_1(n-1)t}
\mathbb{C}{\rm ov}_{\delta_x}(\langle f,X_{(n-1)t}\rangle, \langle h,X_{(n-1)t}\rangle).
\end{eqnarray*}
Thus, by \eqref{small}, \eqref{1.49} and \eqref{cov:sc}, we easily get
\begin{eqnarray}\label{3.20}
 &&\left| V_t^n(x)-(\theta_2^2\rho_{h}^2+\theta_3^2\sigma_{f}^2)\phi_1(x)\right|
 \lesssim \left(c_{(n-1)t}+{t}^{-1}+{t}^{-(1+2\tau(h))/2}\right)(b_{t_0}(x)^{1/2}+b_{t_0}(x)),
\end{eqnarray}
where $c_{t}\to 0$ as $t\to \infty$.
By \eqref{2.8}, we get, as  $t\to\infty,$
$$
e^{\lambda_1t}T_t\left|V_t^n(x)-(\theta_2^2\rho_{h}^2+\theta_3^2\sigma_{f}^2)\phi_1(x)\right|(x)
\lesssim \left(c_{(n-1)t}+t^{-1}+t^{-(1+2\tau(h))/2}\right)e^{\lambda_1t}T_t(\sqrt{b_{t_0}}+b_{t_0})(x)
\to 0,
$$
which implies
\begin{equation}\label{3.38}
\lim_{t\to\infty}e^{\lambda_1t}\sum_{u\in \cL_t}V_t^n(z_u(t))
=\lim_{t\to\infty}e^{\lambda_1t}(\theta_2^2\rho_h^2+\theta_3^2\sigma_f^2)\langle\phi_1,X_t\rangle
=(\theta_2^2\rho_h^2+\theta_3^2\sigma_f^2) W_\infty,
\end{equation}
in probability.

By Lemma \ref{lem:5.6}, we get, as $s\to\infty$, $g_s(x,\theta_2,\theta_3)\to 0$.
Since
$$
g_{(n-1)t}(x,\theta_2,\theta_3)\leq V_t^n(x)\lesssim b_{t_0}(x)^{1/2}+b_{t_0}(x)\in L^2(E,m),
$$
by the dominated convergence theorem, we have that for any $x\in E$,
$$
\lim_{t\to\infty}\|g_{(n-1)t}(x,\theta_2,\theta_3)\|_2=0.
$$
By Lemma \ref{2.8}, we have that for any $x\in E$,
$$
e^{\lambda_1t}\P_{\delta_x}\langle g_{(n-1)t}(\cdot,\theta_2,\theta_3),X_t\rangle
\lesssim \|g_{(n-1)t}(\cdot,\theta_2,\theta_3)\|_2b_{t_0}(x)^{1/2}\to 0,\quad \mbox{as}\quad t\to\infty,
$$
which implies
\begin{equation}\label{3.21}
  e^{\lambda_1t}\sum_{u\in \cL_t}g_{(n-1)t}(z_u(t),\theta_2,\theta_3)
  =e^{\lambda_1t}\langle g_{(n-1)t}(x,\theta_2,\theta_3),X_t\rangle\to0,
\end{equation}
in probability.
Thus, for any sequence $s_k\to\infty$, there exists a subsequence $s'_k$ such that,
if we let $t_k=s'_k$, $m_k=|X_{t_k}|$ and $\{a_{k, j}, j=1, \dots, m_k\}=
\{z_u(t_k), u\in {\cal L}_{t_k}\}$, then the Lindeberg conditions hold
$\P_{\delta_x}$-a.s.
Therefore, by \eqref{3.17}, we have
\begin{equation}\label{3.24}
  \lim_{t\to\infty}\prod_{u\in\cL_t}h\left((n-1)t,z_u(t),e^{\lambda_1t/2}, \theta_2,\theta_3\right)
  =\exp\left\{-\frac{1}{2}\left(\theta_2^2\rho_h^2+\theta_3^2\sigma_f^2\right)W_\infty\right\},\quad\mbox{in probability}.
\end{equation}
Hence by the dominated convergence theorem, we get
\begin{equation}\label{3.25}
 \lim_{t\to\infty}\kappa_2(\theta_1,\theta_2,\theta_3,n,t)
 =\P_{\delta_x}\left(\exp\left\{i\theta_1W_\infty\right\}
 \exp\left\{-\frac{1}{2}\left(\theta_2^2\rho_h^2+\theta_3^2\sigma_f^2\right)W_\infty\right\}\right),
\end{equation}
which implies our claim \eqref{9.5}.

By \eqref{9.5} and the fact that
$e^{\lambda_1 nt}\langle \phi_1,X_{nt}\rangle-e^{\lambda_1 t}\langle \phi_1,X_{t}\rangle\to 0$, in probability, as $t\to\infty$ , we easily get that
\begin{eqnarray*}
&&U_3(n,t)\\
&&:=\left(e^{\lambda_1 nt}\langle \phi_1,X_{nt}\rangle,
\frac{e^{\lambda_1nt/2}(\langle h,X_{nt}\rangle-\langle T_{(n-1)t}h,X_t\rangle)}{(nt)^{(1+2\tau(h))/2}},
e^{\lambda_1nt/2}(\langle f,X_{nt}\rangle-\langle T_{(n-1)t}f,X_t\rangle)\right)\\
&&\stackrel{d}{\to}\left(W_\infty, \left(\frac{n-1}{n}\right)^{(1+2\tau(h))/2}\sqrt{W_\infty}G_2(h),\sqrt{W_\infty}G_1(f)\right).
\end{eqnarray*}
Using \eqref{4.49} with $s=(n-1)t$, we get that, if $\gamma(f)<\infty$,
$$
\P_{\delta_x}\langle T_{(n-1)t}f,X_t\rangle^2
\lesssim (nt)^{2\tau(f)}e^{-2nt\Re_{\gamma(f)}}b_{t_0}(x)^{1/2}
+((n-1)t)^{2\tau(f)}e^{-\lambda_1t}e^{-2\Re_{\gamma(f)}(n-1)t}b_{t_0}(x)^{1/2}.
$$
If $\gamma(f)=\infty$, using \eqref{4.49'} with $s=(n-1)t$, we get
$$
\P_{\delta_x}\langle T_{(n-1)t}f,X_t\rangle^2
\lesssim b_{t_0}(x)^{1/2}+e^{-\lambda_1t}b_{t_0}(x)^{1/2}.
$$
Therefore, we have
\begin{equation}\label{4.7}
\lim_{t\to\infty}e^{\lambda_1nt}\P_{\delta_x}\langle T_{(n-1)t}f,X_t\rangle^2=0.
\end{equation}
By \eqref{4.4}, when $\lambda_1=2\Re_{\gamma(h)}$, we get
\begin{eqnarray}\label{9.3}
   && \int_0^{t-2t_0}T_{t-u}[A(T_{u+(n-1)t}h)^2](x)\,du \nonumber\\
   &\lesssim &e^{-\lambda_1nt}\int_0^{t-2t_0}\left(u+(n-1)t\right)^{2\tau(h)}\,dub_{t_0}(x)^{1/2}\lesssim n^{2\tau(h)} t^{1+2\tau(h)}e^{-\lambda_1nt}b_{t_0}(x)^{1/2}.
\end{eqnarray}
By \eqref{9.3}, \eqref{9.7} and \eqref{4.3}, when $\lambda_1=2\Re_{\gamma(h)}$, we have
$$
\P_{\delta_x}\langle T_{(n-1)t}h,X_t\rangle^2
\lesssim n^{2\tau(h)}t^{1+2\tau(f)}e^{-\lambda_1nt}b_{t_0}(x)^{1/2}
+(nt)^{2\tau(h)}e^{-\lambda_1nt}b_{t_0}(x)^{1/2}.
$$
Therefore, we have
\begin{equation}\label{4.2}
  \lim_{n\to\infty}\limsup_{t\to\infty}(nt)^{-(1+2\tau(h))}e^{\lambda_1nt}\P_{\delta_x}\langle T_{(n-1)t}h,X_t\rangle^2=0.
\end{equation}
Let $\mathcal{D}(nt)$ and $\widetilde{\mathcal{D}}^n(t)$ be the distributions of
$U_1(nt)$ and $U_3(n,t)$ respectively,
and let $\mathcal{D}^n$ and $\mathcal{D}$ be those of
$\left(W_\infty, \left(\frac{n-1}{n}\right)^{(1+2\tau(h))/2}\sqrt{W_\infty}G_2(h),\sqrt{W_\infty}G_1(f)\right)$
and $\left(W_\infty, \sqrt{W_\infty}G_2(h),\sqrt{W_\infty}G_1(f)\right)$ respectively.
Then, using \eqref{5.20}, we have
\begin{eqnarray}\label{4.12}
 &&\limsup_{t\to\infty}\beta(\mathcal{D}(nt),\mathcal{D})\leq
  \limsup_{t\to\infty}[\beta(\mathcal{D}(nt),\widetilde{\mathcal{D}}^n(t))
  +\beta(\widetilde{\mathcal{D}}^n(t),\mathcal{D}^n)+\beta(\mathcal{D}^n,\mathcal{D})]\nonumber\\
 &\leq &\limsup_{t\to\infty}\left((nt)^{-(1+2\tau(h))}e^{\lambda_1 nt}{\P}_{\mu}\langle T_{(n-1)t}h,X_t\rangle^2
 +e^{\lambda_1 nt}{\P}_{\mu}\langle T_{(n-1)t}f,X_t\rangle^2\right)^{1/2}
 +0+\beta(\mathcal{D}^n,\mathcal{D}).\nonumber\\
\end{eqnarray}
Using  the definition of
 $\limsup_{t\to\infty}$, \eqref{4.7} and \eqref{4.2}, we easily get that
 \begin{eqnarray*}
 &&\limsup_{t\to\infty}\beta(\mathcal{D}(t),\mathcal{D})=
 \limsup_{t\to\infty}\beta(\mathcal{D}(nt),\mathcal{D})\\
 & \le& \limsup_{t\to\infty}(nt)^{-(1+2\tau(h))}e^{nt\lambda_1t}\P_{\delta_x}\langle T_{(n-1)t}h,X_t\rangle^2
 +\beta(\mathcal{D}^n,\mathcal{D}).
 \end{eqnarray*}
Letting $n\to\infty$, we get $ \limsup_{t\to\infty}\beta(\mathcal{D}(t),\mathcal{D})=0$.
The proof is now complete. \hfill$\Box$

\textbf{Proof of Corollary \ref{cor:2}:}
Define
$$
  Y_1(s):=s^{-(1+2\tau(h_1))/2}e^{\lambda_1s/2}\langle h_1,X_{s}\rangle,\quad Y_2(s):=s^{-(1+2\tau(h_2))/2}
e^{\lambda_1s/2}\langle h_2,X_{s}\rangle
$$
and
$$Y_s(\theta_2,\theta_3):=\theta_2Y_1(s)+\theta_3Y_2(s).$$
Thus, we have
\begin{eqnarray}\label{7.10}
  {\V}{\rm ar}_{\delta_x}Y_{(n-1)t}(\theta_2,\theta_3)&=&\theta_2^2{\V}{\rm ar}_{\delta_x}Y_1((n-1)t)
  +\theta_3^2{\V}{\rm ar}_{\delta_x}Y_2((n-1)t)\nonumber\\
 &&\quad +2\theta_2\theta_3{\mathbb{C}}{\rm ov}_{\delta_x}(Y_1((n-1)t),Y_2((n-1)t)).
\end{eqnarray}
By \eqref{7.49} and \eqref{1.49}, we get
$$\left|{\V}{\rm ar}_{\delta_x}Y_{(n-1)t}(\theta_2,\theta_3)-(\theta_2^2\rho_{h_1}^2
+\theta_3^2\rho_{h_2}^2+2\theta_2\theta_3\rho(h_1,h_2))\phi_1(x)\right|
\lesssim t^{-1}\left(b_{t_0}(x)^{1/2}+b_{t_0}(x)\right).$$
Using arguments similar to those leading to Lemma \ref{lem:cs}, we get
\begin{eqnarray}\label{7.11}
  &&\lim_{t\to\infty}\P_{\delta_x}\exp\left\{i\theta_1e^{\lambda_1t}\langle \phi_1,X_t\rangle+i\theta_2Y_1(t)+i\theta_3Y_2(t)\right\}\nonumber\\
&=&\P_{\delta_x}\exp\left\{i\theta_1W_\infty-\frac{1}{2}\left(\theta_2^2\rho_{h_1}^2+\theta_3^2\rho_{h_2}^2+2\theta_2\theta_3\rho(h_1,h_2)\right)W_\infty\right\}.
\end{eqnarray}
The proof of Corollary \ref{cor:2} is now complete.\hfill$\Box$

\smallskip

Recall that
$$g(x)=\sum_{k: \lambda_1>2\Re_k}\Phi_k(x)^Tb_k\in \C_c\quad \mbox{and }\quad I_sg(x)
=\sum_{k: \lambda_1>2\Re_k}e^{\lambda_ks}\Phi_k(x)^TD_k(s)^{-1}b_k.$$
We can show that $I_sg$ is real.
In fact, for $k$ with $\lambda_1>2\Re_k$, we have $\lambda_1>2\Re_{k'}$.
And
$$e^{\lambda_{k'}s}\Phi_{k'}(x)^TD_{k'}(s)^{-1}b_{k'}
=e^{\overline{\lambda_k}s}\overline{\Phi_k(x)^T}D_k(s)^{-1}\overline{b_k}
=\overline{e^{\lambda_ks}\Phi_k(x)^TD_k(s)^{-1}b_k},$$
which implies that $I_sg(x)$ is real.
Define
$$
H_\infty:=\sum_{k: \lambda_1>2\Re_k}H_\infty^{(k)}b_k.
$$
By Lemma \ref{thrm1}, we have, as $s\to\infty$
$$
\langle I_sg,X_s\rangle\to H_\infty,\quad \P_{\delta_x}\mbox{-a.s.}
\quad \mbox{and in } L^2(\P_{\delta_x}).
$$
Since $\P_{\delta_x}\langle I_sg,X_s\rangle=g(x)$,
we get
\begin{equation}\label{L1H}
  \P_{\delta_x}(H_\infty)=g(x).
\end{equation}
By \eqref{1.13},
we have
\begin{eqnarray}\label{var:Iu}
  \P_{\delta_x}\langle I_sg,X_s\rangle^2
  =\int_0^sT_u\left[ A\left|I_ug\right|^2\right](x)\,du+T_s[(I_sg)^2](x).
\end{eqnarray}
It is easy to see that,
$$|I_sg(x)|^2\lesssim \sum_{k: \lambda_1>2\Re_k}e^{2\Re_ks}s^{2\nu_k}b_{4t_0}(x).$$
Thus, by \eqref{2.8}, we have, for $s>2t_0$,
\begin{equation}\label{5.5}
  T_s|I_sg|^2(x)\lesssim \sum_{k: \lambda_1>2\Re_k}e^{2\Re_ks}s^{2\nu_k}T_s(b_{t_0})(x)
  \lesssim \sum_{k: 2\Re_k<\lambda_1}s^{2\nu_k}e^{(2\Re_k-\lambda_1)s}b_{t_0}(x)^{1/2}.
\end{equation}
By \eqref{1.37}, we get
\begin{eqnarray*}
  &&\int_0^\infty T_u\left[ A\left|I_ug\right|^2\right](x)\,du \\
  &\lesssim & \sum_{k: \lambda_1>2\Re_k}\left(\int_0^{2t_0} e^{2\Re_ku}u^{2\nu_k}T_u(b_{4t_0})(x)\,du
  +\int_{2t_0}^\infty u^{2\nu_k}e^{(2\Re_k-\lambda_1)u}\,dub_{t_0}(x)^{1/2}\right)\\
   &\lesssim &b_{t_0}(x)^{1/2}  \in L^2(E,m)\cap L^4(E,m).
\end{eqnarray*}
Therefore, by \eqref{var:Iu} and \eqref{5.5}, we get
\begin{equation}\label{L2H}
  \P_{\delta_x}(H_\infty)^2=\lim_{s\to\infty}\P_{\delta_x}|\langle I_sg,X_s\rangle|^2
=\int_{0}^\infty T_u\left[ A\left|I_ug\right|^2\right](x)\,du\in L^2(E,m)\cap L^4(E,m).
\end{equation}
Hence, we have
\begin{equation}\label{LH}
  {\V}{\rm ar}_{\delta_x}H_\infty=\P_{\delta_x}(H_\infty)^2-(\P_{\delta_x}H_\infty)^2
=\int_{0}^\infty T_u\left( A\left|I_ug\right|^2\right)(x)\,du-g(x)^2.
\end{equation}

\textbf{Proof of Theorem \ref{The:1.3}:}
Recall that
$$E_t(g)=\left(\sum_{k: 2\lambda_k<\lambda_1}e^{-\lambda_kt}H^{(k)}_\infty D_k(t)b_k\right)$$
and
$$
  Y_1(t):=e^{\lambda_1t/2}\langle f,X_{t}\rangle,\quad  Y_2(t):=t^{-(1+2\tau(h))/2}e^{\lambda_1t/2}\langle h,X_{t}\rangle.
$$
Consider an $\mathbb{R}^4$-valued random variable $U_4(t)$ defined by:
 \begin{eqnarray*}\label{8.5}
   &&U_4(t):= \left(e^{\lambda_1 t}\langle\phi_1,X_t\rangle,
   ~e^{\lambda_1t/2}\left(\langle g,X_t\rangle-E_t(g)\right),
   Y_2(t), Y_1(t)\right).
 \end{eqnarray*}
 To get the conclusion of Theorem \ref{The:1.3}, it suffices to show that, under $\P_{\delta_x}$,
\begin{equation}\label{2.5a}
   U_4(t)\stackrel{d}{\to}\left(W_\infty, \sqrt{W_\infty}G_3(g), \sqrt{W_\infty}G_2(h),\sqrt{W_\infty}G_1(f)\right),
\end{equation}
where $W_\infty$, $G_3(g)$, $G_2(h)$ and $G_1(f)$ are independent.
Denote the characteristic function of $U_4(t)$  under  $\P_{\delta_x}$ by $\kappa_3(\theta_1,\theta_2,\theta_3,\theta_4,t)$.
Then, we only need to prove
 \begin{equation}\label{8.7}
   \lim_{t\to\infty}\kappa_3(\theta_1,\theta_2,\theta_3,\theta_4,t)
   =\P_{\mu}\left(\exp\{i\theta_1W_\infty\}
   \exp\left\{-\frac{1}{2}(\theta_2^2\beta_g^2 +\theta_3^2\rho_h^2 +\theta_4^2\sigma^2_f)W_\infty\right\}\right).
 \end{equation}
Note that, by Lemma \ref{thrm1}, we get
 $$
 E_t(g)=\lim_{s\to\infty}\langle I_sg,X_{t+s}\rangle=\sum_{u\in\cL_t}\lim_{s\to\infty}\langle I_sg,X_{s}^{u,t}\rangle.
 $$
Since $X_{s}^{u,t}$ has the same law as $X_s$ under $\P_{\delta_{z_u(t)}}$,
$H^{u,t}_\infty:=\lim_{s\to\infty}\langle I_sg,X_{s}^{u,t}\rangle$ exists
and has the same law as $H_\infty$ under $\P_{\delta_{z_u(t)}}$.
Thus, we get $E_t(g)=\sum_{u\in\cL_t}H^{u,t}_\infty$.
Let $ h(x,\theta)=\P_{\delta_x}\exp\left\{i\theta(H_\infty-g(x))\right\}$.
Therefore, we obtain that
\begin{eqnarray}\label{8.6}
  &&\kappa_3(\theta_1,\theta_2,\theta_3,\theta_4,t)\nonumber\\
  &=&
   \P_{\delta_x}\left(\exp\left\{i\theta_1e^{\lambda_1 t}\langle\phi_1,X_t\rangle+i\theta_3Y_2(t)+i\theta_4Y_1(t)\right\}\prod_{u\in\cL_t}h\left(z_u(t),-\theta_2e^{\lambda_1t/2}\right)\right).
\end{eqnarray}
Let $V(x)={\V}{\rm ar}_{\delta_x}H_\infty$.
 We claim that
\begin{description}
  \item{(i)}
  as $t\to\infty$,
  \begin{equation}\label{8.1}
    e^{\lambda_1t}\sum_{u\in{\cL}_t}\P_{\delta_x}|H^{u,t}_\infty-g(z_u(t))|^2
  =e^{\lambda_1t}\langle V,X_t\rangle\to\langle V,\psi_1\rangle_m W_\infty, \mbox{ in probability};
  \end{equation}
  \item{(ii)} for any $\epsilon>0$,
  as $t\to\infty$,
  \begin{eqnarray}\label{8.2}
   &&e^{\lambda_1t}\sum_{u\in{\cL}_t}
    \P_{\delta_x}(|H^{u,t}_\infty-g(z_u(t))|^2,|H^{u,t}_\infty-g(z_u(t))|>\epsilon e^{-\lambda_1t/2})\nonumber\\
    &=&e^{\lambda_1t}\langle k(\cdot,t),X_t\rangle\to 0,
    \mbox{ in probability,}
  \end{eqnarray}
  where
  $k(x,t):=\P_{\delta_x}(|H_\infty-g(x)|^2,|H_\infty-g(x)|>\epsilon e^{-\lambda_1t/2}).$
\end{description}
Then using arguments similar to those in the proof Lemma \ref{lem:cs}, we have
\begin{equation}\label{8.4}
  \prod_{u\in{\cL}_t}h\left(z_u(t),-\theta_2e^{(\lambda_1/2)t}\right)
  \to \exp\left\{-\frac{1}{2}\theta_2^2\langle V,\psi_1\rangle_m W_\infty\right\},\mbox{ in probability.}
\end{equation}

Now we prove the claims.

(i)
By \eqref{L2H}, we have $V(x)\in L^2(E,m)\cap L^4(E,m)$.
By Remark \ref{rem:large}, \eqref{8.1} follows immediately.

(ii)
We easily see that $k(x,t)\downarrow0$ as  $t\uparrow\infty$ and $k(x,t)\leq V(x)\in L^2(E,m)$ for any $x\in E$.
Thus, $\lim_{t\to\infty}\|k(\cdot,t)\|_2=0.$
So, by \ref{2.8}, we have that for any $x\in E$,
\begin{equation*}
  e^{\lambda_1t}\P_{\delta_x}\langle k(\cdot,t),X_t\rangle
  \lesssim \|k(\cdot,t)\|_2b_{t_0}(x)^{1/2}\to 0,\quad\mbox{ as }t\to\infty,
\end{equation*}
which implies \eqref{8.2}.

By \eqref{8.6}, \eqref{8.4} and the dominated convergence theorem, we get that as $t\to\infty$,
\begin{eqnarray}\label{8.30}
  &&\left|\kappa_3(\theta_1,\theta_2,\theta_3,\theta_4,t)
  -\P_{\delta_x}\left(\exp\left\{(i\theta_1-\frac{1}{2}\theta_2^2\langle V,\psi_1\rangle_m)
e^{\lambda_1 t}\langle\phi_1,X_t\rangle+i\theta_3Y_2(t)+i\theta_4Y_1(t)\right\}\right)\right|\nonumber\\
  &\le&\P_{\delta_x}\left|\prod_{u\in{\cL}_t}h\left(z_u(t),-\theta_2e^{(\lambda_1/2)t}\right)
  -\exp\left\{-\frac{1}{2}\theta_2^2\langle V,\psi_1\rangle_m e^{\lambda_1t}\langle\phi_1,X_t\rangle\right\}\right|
  \to 0.
\end{eqnarray}
By Lemma \ref{lem:cs}, we get
\begin{eqnarray*}
  &&\lim_{t\to\infty}\kappa_3(\theta_1,\theta_2,\theta_3,\theta_4,t) \\
  &=& \lim_{t\to\infty}\P_{\delta_x}
  \left(\exp\left\{(i\theta_1-\frac{1}{2}\theta_2^2\langle V,\psi_1\rangle_m)e^{\lambda_1 t}\langle\phi_1,X_t\rangle
  +i\theta_3Y_1(t)+i\theta_4Y_2(t)\right\}\right) \\
  &=& \P_{\delta_x}\left(\exp\{i\theta_1W_\infty\}
   \exp\left\{-\frac{1}{2}\left(\theta_2^2\langle V,\psi_1\rangle_m+\theta_3^2\rho_h^2 +\theta_4^2\sigma^2_f\right)W_\infty\right\}\right).
\end{eqnarray*}
By \eqref{LH}, we get
$$
\langle V,\psi_1\rangle_m=\int_{0}^\infty e^{-\lambda_1u}\langle A\left|I_ug\right|^2,\psi_1\rangle_m\,du-\langle g^2,\psi_1\rangle_m.
$$
The proof is now complete.\hfill$\Box$

\begin{singlespace}

\end{singlespace}

\vskip 0.2truein
\vskip 0.2truein

\noindent{\bf Yan-Xia Ren:} LMAM School of Mathematical Sciences \& Center for
Statistical Science, Peking
University,  Beijing, 100871, P.R. China. Email: {\texttt
yxren@math.pku.edu.cn}

\smallskip
\noindent {\bf Renming Song:} Department of Mathematics,
University of Illinois,
Urbana, IL 61801, U.S.A.
Email: {\texttt rsong@math.uiuc.edu}

\smallskip

\noindent{\bf Rui Zhang:} LMAM School of Mathematical Sciences, Peking
University,  Beijing, 100871, P.R. China. Email: {\texttt
ruizhang8197@gmail.com}

\end{doublespace}
\end{document}